\documentclass[11pt]{amsart}

\usepackage{amsfonts,amssymb}
\usepackage{stmaryrd} 
\usepackage{tikz,color}  
\usepackage[headings]{fullpage} 
\usepackage[all]{xy}

\usetikzlibrary{arrows}
\usetikzlibrary{decorations.pathmorphing}
\usepackage[colorlinks]{hyperref} 

\theoremstyle{plain}
\newtheorem{thm}{Theorem}[section]

\newtheorem{cor}[thm]{Corollary}

\theoremstyle{definition}
\newtheorem{definition}[thm]{Definition}
\newtheorem{example}[thm]{Example}

\newtheorem{construction}[thm]{Construction}

\def\eq{\addtocounter{thm}{1}\begin{equation}} 
\def\myalign{\addtocounter{thm}{1}\begin{align}} 

\newenvironment{mylist}
{\begin{list}{\text} 
{\setlength{\rightmargin}{0cm}
\setlength{\itemindent}{0cm}
\setlength{\labelwidth}{1cm}
\setlength{\leftmargin}{1.2cm}
\setlength{\labelsep}{.2cm}
}}
{\end{list}}

\makeatletter
\font \sectionfont = cmbx12 
\renewcommand\section{\@startsection{section}%
{1}{0pt}{-1.8\baselineskip}{.8\baselineskip}%
{\sectionfont \raggedright}}
\makeatother

\makeatletter
\def\theenumi{\@roman\c@enumi}
\def\@linkcolor{black}
\def\@citecolor{black}
\makeatother

\def\b{\begin}
\def\l{\label} 
\def\e{\end} 

\def\qed{\unskip\nobreak\hfil\penalty100\hfill\null\nobreak\hfil
{\vbox{\hrule\hbox{\vrule\kern3.1pt
\vbox{\kern6pt}\kern3pt\vrule}\hrule}}}
{\parfillskip=0pt\finalhyphendemerits=0\par}  

\font \prooffont = cmcsc10 
\def \proof{\vskip 0mm\noindent {\prooffont Proof:\ \ }} 

\font \sectionfont = cmbx12  

\def \myspace{\vglue .3cm}  

\def\to{\longrightarrow} 

\def \a{{\mathrel{\smash-}}{\mathrel{\mkern-8mu}}
{\mathrel{\smash-}}{\mathrel{\mkern-8mu}}
{\mathrel{\smash-}}{\mathrel{\mkern-8mu}}}

\def\ar #1{{
{\mathop{{\mathrel{\smash-}}{\mathrel{\mkern-8mu}}\a\a\a
{\mathrel{\smash-}}{\mathrel{\mkern-8mu}}\rightarrow}\limits^{#1}}
}} 

\def\ar#1{\mathop{{\mathrel{\smash-}}{\mathrel{\mkern-8mu}}
{\mathrel{\smash-}}{\mathrel{\mkern-8mu}}\rightarrow}\limits^{#1}}

\def\Id{{\rm Id}}

\def\A{{\bf A}}
\def\YA{{\bf Y}}

\def\B{{\bf B}}

\def\K{{\bf K}}
\def\T{{\bf T}}
\def\W{{\bf W}}
\def\U{{\bf U}}

\def\Q{{\bf Q}}
\def\P{{\bf P}}
\def\F{{\bf F}}

\def\G{{\bf G}}
\def\V{{\bf V}}
\def\Y{{\bf Y}}
\def\L{{\bf L}}
\def\M{{\bf M}}
\def\C{\mathcal C}

\def\coker{{\rm{Coker}}}

\def\deg{{\rm{deg}}}

\def\dim{{\rm{dim}}}
\def\ext{{\rm{Ext}}}

\def\hom{{\rm{Hom}}}
\def\id{{\rm{Id}}}
\def\im{{\rm{Im}}}

\def\ker{{\rm{Ker}}}

\def\rank{{\rm{rank}\,}}
\def\reg{{\rm{reg}}}

\def\syz{{\rm{Syz}}}
\def\tor{{\rm{Tor}}}
\def\coker{{\rm{Coker}}}
\def\MM{{\bf Cone}}

\def\to{\longrightarrow}
\def\Box{\rm{Box}}
\def\poin{{\mathcal P}}
\def\RR{{\mathcal R}}
\def\cosyz{{\rm Cosyz}}
\def\std{{\rm Sh}}
\def\Hom{{\rm Hom}}
\def\KK{{\bf K}}

\def\mm{{\bf m}}

\def\ff{{f_{1}, \dots, f_{c}}}
\def\ff#1{{f_{1},\dots, f_{#1}}}
\def\Id{{\rm Id}}
\def\tilde{\widetilde}
\def\MCM{\underline{\bf MCM}}
\def\mod{{{\rm mod}}}

\def\MF{{\bf MF}}

\def\Mor{{\rm Mor}}
\def\s{{\sigma}}

\def\extt{{\widehat{\ext}}}
\def\highsyz{{stable  syzygy\ }}

\begin{document}

\title{
MATRIX FACTORIZATIONS\\
for COMPLETE INTERSECTIONS\\
and MINIMAL FREE RESOLUTIONS}

\author[David Eisenbud]{David Eisenbud}
\address{Mathematics Department, University of California at Berkeley, Berkeley, CA 94720, USA}

\author[Irena Peeva]{Irena Peeva}
\address{Mathematics Department, Cornell University, Ithaca, NY 14853, USA}

\begin{abstract}
We describe the asymptotic structure of minimal free resolutions over complete intersections of arbitrary codimension.
To do this we define a \emph{higher matrix factorization of a regular sequence $f_{1}, \dots, f_{c}$} in a way that extends Eisenbud's definition of a matrix factorization of one element. Using this notion we can
describe the minimal free resolutions, both over a regular local ring $S$ and over the complete intersection ring $R = S/(f_{1}, \dots, f_{c})$, of modules that are high syzygies over $R$.
\end{abstract}


\maketitle

\section{Introduction} \label{Sintro} 

Let $S$ be a regular local ring and $f_{1}, \dots, f_{c}$ be a regular sequence.
If $N$ is a  finitely generated module over the complete intersection  $R:= S/(f_{1}, \dots, f_{c})$, then it
can be also considered as an $S$-module annihilated by $f_{1}, \dots, f_{c}$.
 In this paper we will describe the minimal free resolutions of $N$ as an $S$ module and as an $R$-module when $N$ is a high syzygy over $R$.

The case when $N$ is the residue field of $S$ is classical: its minimal free resolution over $S$ is the Koszul complex. Perhaps motivated by questions of group cohomology, Tate \cite{Ta}, in 1957, gave an elegant description of its minimal free resolution over $R$.  

The understanding of minimal resolutions of an arbitrary module $N$ over $R$  began with the 1974 paper  \cite{Gu} of Gulliksen, who showed that  $\ext_{R}(N,k)$ can be regarded as a finitely generated graded module over a polynomial ring $\RR = k[\chi_{1},\dots,\chi_{c}]$, where $c$ is the codimension of $R$. 
He used this to show
 that the Poincar\'e series $\sum_{i}b_{i}^{R}(N)x^{i}$, the generating function of the Betti numbers $b_{i}^{R}(N)$, is rational and that the denominator divides $(1-x^2)^c$.  

Gulliksen's finite generation result implies that the even Betti numbers 
$b_{2i}^{R}(N)$ are eventually given by a polynomial in $i$, and similarly for the odd Betti numbers. In 1989 Avramov~\cite{Av}  proved that the two polynomials have the same leading coefficient, and he also 
extended  constructions from group cohomology to the general case.
 In 1997 Avramov, Gasharov and Peeva \cite{AGP} gave further restrictions on the Betti numbers, 
establishing in particular that the Betti sequence $\{b_{i}^{R}(N)\}_{i\ge q}$ is either strictly increasing or constant for $q\gg 0$.
Examples in \cite{Ei1} and  \cite{AGP} show that, as with the Betti numbers,
minimal free resolutions over a complete intersection can have
intricate structure, but the examples exhibit stable patterns when sufficiently truncated.

The theory of matrix factorizations entered the picture in the 1980 paper \cite{Ei1} of Eisenbud, who introduced them to describe the minimal free resolutions of modules that are high syzygies over hypersurface rings---the case of codimension one. They have had many applications:

Starting with Kapustin and Li \cite{KL}, who followed an  idea of Kontsevich, physicists
discovered  amazing connections with string theory --- see  \cite{As} for a survey.
A major advance was made by Orlov \cite{Or1, Or3, Or4, Or5},  who showed  that matrix factorizations could be used to study
Kontsevich's homological mirror symmetry by giving a new description of singularity categories.
Matrix factorizations have also proven
useful for the study of  cluster tilting \cite{DH}, 
Cohen-Macaulay modules and singularity theory \cite{BGS, BHU, CH, Kn}, Hodge theory \cite{BFK},  Khovanov-Rozansky homology \cite{KR1, KR2}, moduli of curves \cite{PV2},  quiver and group representations \cite{AM, Av, KST, Re}, 
and other topics, for example, \cite{BDFIK, CM, DM, Dy, Ho, HW, Is,
PV1, Se, Sei, Sh}. 

Orlov \cite{Or2} and subsequent authors, for example   \cite{Bu, BW, PV2}, have studied modules over a complete intersection
$S/(f_1,\allowbreak\dots , f_c)$  by
reducing to families of codimension $1$ matrix factorizations over the hypersurace $\sum z_{i}f_{i} = 0$
in the projective space $\P^{c-1}_{S}$, where the $z_{i}$ are the homogeneous coordinates of $\P^{c-1}_S$.
By contrast, our theory is focused on understanding  minimal free resolutions.

Minimal free resolutions of high syzygies over a  codimension two complete intersection were constructed by 
Avramov and Buchweitz in \cite{AB} in 2000 using the classification of modules over the exterior algebra on two variables. In higher codimension,
 non-mimimal resolutions have been known  for over forty years from the work of Shamash \cite{Sh}, but minimal free resolutions, which carry much more information and exhibit much more varied behavior, have remained mysterious.
We introduce the concept of higher matrix factorization in order to describe the structure of minimal resolutions of high syzygies.

\myspace
\noindent{\bf What is a Matrix Factorization?} 
\vglue .15cm
\noindent  We briefly review the codimension 1 case. If $0\neq f\in S$ is an element in a commutative ring then a {\it matrix factorization} of $f$ is a pair $(d,h)$
of maps of finitely generated free modules
$$
A_{0}\ar{h} A_{1} \ar{d} A_{0}
$$
such that the diagram 

\begin{center}
\begin{tikzpicture}  [scale=0.9, every node/.style={scale=0.9}]
\node(20){$A_{1}$};
\node(21)[node distance=2.5cm, right of=20]{$A_{0}$};
\node(22)[node distance=2.5cm, right of=21]{$A_{1}$};
\node(23)[node distance=2.5cm, right of=22]{$A_{0}$};
\draw[->, >=stealth] (20) to node [below=2pt] {$d$} (21);
\draw[->] (21) to node [below=2pt]{$h$} (22);
\draw[->] (22) to node [below=2pt] {$d$} (23);
\draw[->, bend left=20, red] (21) to node [above=2pt] {$f$} (23);
\draw[->, bend left=20, red, >=stealth] (20)  to node [above=2pt]  {$f$} (22);
\end{tikzpicture}
\end{center}

\noindent
commutes or, equivalently:  
\b{align*}
dh &= f\cdot \Id_{A_{0}}\cr
hd &= f\cdot \Id_{A_{1}}\, .
\end{align*}
If $f$ is a non-zerodivisor 
and $S$ is local, then the matrix factorization describes  the minimal free resolutions
of $M := \coker (d)$ over the rings $S$ and  $R:= S/(f)$; if $M$ has no direct summand then the free resolutions are:
\myalign
0\to A_{1}\ar{d} A_{0} \to \,M\to 0\ \hbox{ over $S$; and }&\   \label{twores}\\
\nonumber \cdots\ar{h} R\otimes A_{1}\ar{d} R\otimes A_{0}\ar{h} R\otimes A_{1}\ar{d}R\otimes A_{0}\to &\,M\to 0\  \hbox{ over $R$.}
\end{align}
\noindent
Minimal free resolutions of all sufficiently high syzygies over a hypersurface ring are always of this form  by \cite{Ei1}. 

\myspace

\noindent{\bf What is a Higher Matrix Factorization?}
\vglue .15cm
\noindent To extend the theory to higher codimensions, we make a new definition.
After giving the definition and an example, we will outline the main results of this paper.

\begin{definition} \label{matrixfiltration}
Sometimes we will abbreviate ``higher matrix factorization''  to ``HMF''. 
Let $f_{1}, \dots, f_{c}\in S$ be elements of a commutative ring, and set $R = S/(f_{1}, \dots, f_{c})$.  A  {\it higher matrix factorization  $(d,h)$
with respect to  $f_{1}, \dots, f_{c}$}  is:
\begin{mylist}
\item[(1)] A pair of free finitely generated $S$-modules $A_{0}, A_{1}$ with filtrations
$$
0\subseteq A_{s}(1)\subseteq\cdots \subseteq A_{s}(c) = A_{s},\ \ \hbox{for}\ s=0,1,
$$
such that each $A_{s}(p-1)$ is a free  summand of $A_{s}(p)$;

\item[(2)] A pair of maps $d,\, h$ preserving filtrations,
$$
\bigoplus_{q=1}^{c} A_{0}(q) \ar{h} A_{1}\ar{d} A_{0},
$$
where we regard $\oplus_{q} A_{0}(q)$ as filtered by the submodules
$\oplus_{q\leq p}A_{0}(q)$;
\end{mylist}
such that, writing 
$$
A_{0}(p)\ar{h_p} A_{1}(p)\ar{d_p} A_{0}(p)
$$
for the induced maps, the diagrams

\begin{center}
\begin{tikzpicture}  
 [every node/.style={scale=1},  auto]
\node(20){$A_{1}(p)$};
\node(21)[node distance=2.5cm, right of=20]{$A_{0}(p)$};
\node(22)[node distance=2.5cm, right of=21]{$A_{1}(p)$};
\node(23)[node distance=2.5cm, right of=22]{$A_{0}(p)$};
\node(30)[node distance=2.5cm, below of=20]{$A_{1}(p)/A_{1}(p-1)$};
\node(32)[node distance=2.5cm, below of=22]{$A_{1}(p)/A_{1}(p-1)$};
\draw[->] (20) to node[swap] {} (30);
\draw[->] (22) to node {} (32);
\draw[->] (20) to node [below=1pt] {$d_p$} (21);
\draw[->] (21) to node [below=1pt]{$h_p$} (22);
\draw[->] (22) to node [below=1pt] {$d_p$} (23);
\draw[->, bend left=20, red] (21)  to node[above=0pt]  {$f_{p}$} (23);
\draw[->, bend left=20, red] (30)  to node[above=0pt] {$f_{p}$} (32);
\end{tikzpicture}
\end{center}
\myspace
\noindent commute modulo $(\ff{p-1})$ for all $p$; or, equivalently,
\begin{mylist}
\item[(a)] 
$d_ph_p \equiv f_{p}\, \Id_{A_{0}(p)} 
\hbox{ mod}(f_{1},\dots, f_{p-1})A_{0}(p);$
\item[(b)] 
$\pi_{p}  h_p d_p \equiv f_{p}\,\pi_{p} \hbox{ mod}(f_{1},\dots, f_{p-1})\big(A_1(p)/A_{1}(p-1)\big), $
where $\pi_{p}$ denotes the projection $A_{1}(p) \to A_{1}(p)/A_{1}(p-1)$.
\end{mylist}

\noindent We define  {\it the module of the higher matrix factorization} $(d,h)$  to be 
$$M:=\coker(R\otimes d)\, .$$
\vglue -.1cm \noindent
 We refer to modules of this form as {\it higher matrix factorization modules} or {\it HMF modules}. 

In Section \ref{SMorph}, we show that a homomorphism of HMF modules induces a morphism of the whole higher matrix factorization structure; see   Definition~\ref{mfhom} and
Theorem~\ref{mfmorph} for details. In Section~\ref{Sfunctorial} we show that our constructions yield functors to  stable categories of Cohen-Macaulay modules. In Section~\ref{SHMFstrong}
we give a stronger version of Definition~\ref{matrixfiltration}  requiring that the map $h$ is part of a homotopy, and we prove in Theorem~\ref{strongmf}  that an HMF module always has
such a strong matrix factorization.

For each $1\le p\le c$, we have a higher matrix factorization $(d_p,(h_{1}|\cdots|h_{p}))$ with respect to $f_1,\dots ,f_p$, where $(h_{1}|\cdots|h_{p}))$ denotes the concatenation of the matrices
$h_{1},\dots,h_{p}$ and thus an HMF module $$M(p)=\coker (S/(f_1,\dots ,f_p)\otimes d_p)\,.$$ This allows us to do induction on $p$.

If $S$ is local, then we call the higher matrix factorization  {\it minimal\/} if  $d$ and  $h$ are minimal (that is, the image of each map is contained in the maximal ideal times the target).
\end{definition}

\begin{example} \label{exintro} 
Let $S = k[a,b,x,y]$ over a field $k$, and consider the
complete intersection $R = S/(xa,yb)$. Let $N = R/(x,y)$. The module $N$ is a maximal Cohen-Macaulay $R$-module. The earliest  syzygy of
$N$ that is an HMF module  is the third syzygy $M$. We can describe the higher  matrix factorization
for $M$ as follows. After choosing a splitting $A_s(2) = A_s(1)\oplus B_s(2)$, we can represent the map $d$ as

\myspace
\begin{center}
\begin{tikzpicture}  
 [every node/.style={scale=0.8},  auto]
  \node (TL) {$A_0(1)=B_0(1)=S^2$};
    \node (TM) [node distance=6.5cm, left of=TL] {$A_1(1)=B_1(1)=S^2$};
     \node (M) [node distance=3cm, below of=TM] {$B_1(2) = S^2$};
       \node (L) [node distance=3cm, below of=TL] {$B_0(2)=S\  .$};
   \draw[->] (TM) to node  {$\begin{pmatrix}a&0\\  y&x\end{pmatrix}$} (TL);
     \draw[->] (M) to node [above=10pt, right=-15pt] {$\begin{pmatrix}y&x\end{pmatrix}$} (L);
      \draw[->] (M) to node [above=10pt, left=5pt]  {$\begin{pmatrix}0&-b\\  0&0\end{pmatrix}$} (TL);
\end{tikzpicture}
\end{center}
\myspace

\noindent
The pair of maps
$$d_1: \, A_{1}(1)
 {\mathop{{\mathrel{\smash-}}{\mathrel{\mkern-8mu}}\a\a
{\mathrel{\smash-}}{\mathrel{\mkern-8mu}}\rightarrow}\limits^{ \begin{pmatrix}a &  0   \\ 
                  y  & x\end{pmatrix} }}
 A_{0}(1)\quad\ \ \hbox{and}\quad \ \ h_1: \,  A_{0}(1) 
  {\mathop{{\mathrel{\smash-}}{\mathrel{\mkern-8mu}}\a\a
{\mathrel{\smash-}}{\mathrel{\mkern-8mu}}\rightarrow}\limits^{ \begin{pmatrix} x &  0   \\ 
                  -y  & a \end{pmatrix}
}}
 A_{1}(1)$$
is a  matrix factorization for the element $xa$ since $d_1h_1=h_1d_1=xa\,\id$.
The map 
\def\dashing{\leaders\hbox to 1em{\hfil-\hfil}\hfil}
$h_2: \,  A_{0}=A_{0}(2)\to A_{1}= A_{1}(2)$ is given by the matrix
$$
h_2= \begin{pmatrix}0&b&0\\ 
0&0&0\\ 
x&0&b\\ 
-y&a&0
 \end{pmatrix}\, ,\quad\hbox{and}\ \ \ \ d_2= \begin{pmatrix} a&0&0&-b\\ 
y& x& 0& 0\\ 
0&0&y&x \end{pmatrix}\, .$$
 \vglue -.3cm \noindent
Hence
$$
d_2h_2 = 
 \begin{pmatrix}
yb&0&0\\ 
0&yb&0\\ 
0&xa&yb \end{pmatrix}
\quad\ \ \hbox{and}\quad\ \ 
h_2d_2 =  \begin{pmatrix}
yb&xb&0&0\\ 
0&0&0&0\\ 
\multispan4\dashing\\ 
xa&0&yb&0\\ 
0&xa&0&yb \end{pmatrix}.
$$

Thus $d_2h_2$ is congruent, modulo $(xa)$, to $yb\,\id$.
Furthermore, condition (b) of Definition~\ref{matrixfiltration} is the statement that the two bottom rows in the latter matrix are congruent modulo $(xa)$
to $yb\pi_{2}$.
In the context of the diagram in the definition, with $p=2$, the fact that the lower left $(2\times 2)$-matrix is congruent to 0 modulo $f_{1}=xa$ is necessary for the map $d_{2}h_{2}:A_{1}(2) \to A_{1}(2)$ to induce
a map $A_{1}(2)/A_{1}(1) = B_{1}(2) \to A_{1}(2)/A_{1}(1) = B_{1}(2)$.
\end{example}

\myspace
In the rest of the introduction we focus on the case when $S$ is a regular local ring and
$R = S/(f_{1}, \dots, f_{c})$ is a complete intersection, although most of our results are proved in greater generality.
We will keep the notation of Definition~\ref{matrixfiltration} throughout the introduction.

\myspace
\noindent{\bf High Syzygies are Higher Matrix Factorization Modules}
\vglue .15cm
\noindent 
 The next result was the key motivation for our definition of a higher matrix factorization.
A more precise version of this result is proved in Corollary~\ref{highsyzy}.

\b{thm}\label{highsyzthm}
 Let $S$ be a regular local ring with infinite residue field, 
and let $I\subset S$ be an ideal generated by
a regular sequence of length $c$. Set $R = S/I$, and
suppose that $N$ is a finitely generated  $R$-module. Let
 $f_{1}, \dots, f_{c}$ be a generic choice of elements
minimally generating $I$.
If $M$ is a sufficiently high syzygy of $N$ over $R$, then $M$ is the HMF  module of a minimal higher matrix factorization
 $(d,h)$ with respect to  $f_{1}, \dots, f_{c}$. Moreover $d\otimes R$ and $h\otimes R$ are the first two differentials in the minimal free resolution
of $M$ over $R$.  
\end{thm}

The meaning of ``a sufficiently high syzygy" is explained in  Section~\ref{Sgeneric}, where we  introduce a class of $R$-modules that we call {\it pre-stable syzygies} and show that they have the property given in Theorem~\ref{highsyzthm}. Given an $R$-module $N$ we give in Corollary~\ref{highsyzy} a sufficient condition, in terms of $\ext_{R}(N,k)$, for the $r$-th syzygy module of $N$ to be pre-stable. We also explain more about the genericity condition. Over a local Gorenstein ring, we introduce the concept of a stable syzygy in Section~\ref{Sgeneric} and discuss it in  Section~\ref{SGor}.

\myspace
\noindent{\bf Minimal $R$-free and $S$-free Resolutions}
\vglue .15cm
\noindent Theorem~\ref{highsyzthm} shows that in order to understand the asymptotic behavior of minimal free resolutions over the complete intersection $R$ it suffices to construct the resolutions of HMF modules. This is accomplished by Construction~\ref{inres}  and 
Theorem~\ref{inresthm}.

The finite minimal free resolution over $S$ of an HMF module
is given by Construction~\ref{fres} and Theorem~\ref{fresthm}. 
 Here is an outline of the codimension $2$ case:
Let $(d,h)$ be a codimension $2$ higher matrix factorization.
We first choose   splittings $A_s(2) = B_s(1)\oplus B_s(2)$.  
Since $d(B_{1}(1))\subset B_{0}(1),$ we can represent the differential $d$ as 
\begin{center}
\begin{tikzpicture}  
 [every node/.style={scale=0.9},  auto]
\node(00){};
\node(prepre)[node distance = .5cm, left of=00]{};
\node(preprepre)[node distance=.5cm, left of=prepre]{$\B(1):$};
\node(preprepre2)[node distance=2cm, below of=preprepre]{$\B(2):$};
\node(01)[node distance=2cm, right of=00]{$B_{1}(1)$};
\node(02)[node distance=4cm, right of=01]{$B_{0}(1)\ $};
\node(11)[node distance=2cm, below of=01]{$B_{1}(2)$};
\node(12)[node distance=2cm, below of=02]{$B_{0}(2)\,,$};
\draw[->] (01) to node [ pos=.4]{$b_{1}$} (02);
\draw[->] (11) to node [swap,pos=.39]{$b_{2}$} (12);
\draw[->] (11) to node [above=3pt, left=3pt] {$\psi_2$} (02);
\end{tikzpicture}
\end{center}
\noindent which may be thought of as a map of two-term complexes
$\psi_{2}: \B(2)[-1] \to \B(1)$.
\noindent This extends to a map of complexes $\K(f_{1})\otimes \B(2)[-1]\to \B(1)$, as in the following
diagram:
\begin{center}
\begin{tikzpicture}   [every node/.style={scale=0.8},  auto]
  \node (TL) {$B_0(1)$};
    \node (TM) [node distance=4cm, left of=TL] {$B_1(1)$};
     \node (M) [node distance=2cm, below of=TM] {$B_1(2)$};
       \node (L) [node distance=2cm, below of=TL] {$B_0(2)$};
       \node (BM) [node distance=2cm, below of=M] {$B_0(2)$};
        \node (BL) [node distance=4cm, left of=BM] {$B_1(2)$};
   \draw[->] (TM) to node  {$b_1$} (TL);
     \draw[->] (M) to node [above=10pt, left=-10pt] {$b_2$} (L);
      \draw[->] (M) to node [above=5pt, right=-20pt]  {$\psi_2$} (TL);
     \draw[->] (BL) to node [swap] {$b_2$} (BM);
    \draw[->] (BL) to node [above=20pt, right=-15pt]  [swap] {$h_1\psi_2$} (TM);
 \draw[->] (BL) to node [above=-10pt, right=-8pt] [swap]{$-f_1$} (M);
  \draw[->] (BM) to node [above=-10pt, right=-5pt] [swap]{$f_1$} (L);
\end{tikzpicture}
\end{center}
\noindent Theorem~\ref{fresthm} asserts that this is the minimal $S$-free resolution of the HMF module
$M = \coker(S/(f_{1}, f_{2}) \otimes d)$. 

Strong restrictions on the finite minimal $S$-free resolution of a high syzygy $M$ over the complete intersection  $S/(\ff c)$  follow from our results:  for example, by Corollary~\ref{simplerelations}
the minimal presentation matrix of $M$ must include $c-1$ columns of the form
$$
\b{pmatrix}
f_{1}&\cdots&f_{c-1}\\ 
0 &\cdots&0\\ 
\vdots&&\vdots\\ 
0 &\cdots&0
\end{pmatrix}
$$
 for a generic choice of $f_1,\dots ,f_c$. For instance, in Example~\ref{exintro}, the presentation matrix of $M$ 
 is
 $$\b{pmatrix} a&0&0&-b&0\\ 
y& x& 0& 0&0\\ 
0&0&y&x&xa\end{pmatrix}\, ,$$
and the last column is of the desired type.
There are numerical restrictions as well; see Corollary~\ref{syzbound} and the  remark following it.

Every maximal Cohen-Macaulay
$ S/(f_{1})$-module is a pre-stable syzygy, but this is not true in higher codimension --- one must go further back in the
syzygy chain. This is not surprising, since {\it every} $S$-module of finite length is a maximal Cohen-Macaulay module over an artinian complete intersection, 
and it seems hopeless to characterize the minimal free resolutions of all such modules.

 In Corollary~\ref{bettifinite} and  Corollary~\ref{bettiinf} we get formulas for the Betti numbers of an HMF  module over $S$ and over $R$ respectively. Furthermore, the vector spaces
$\ext_{S}^{i}(M,k)$ and $\ext_{R}^{i}(M,k)$   can be expressed as follows.

\b{cor}\label{corzero}
Suppose that $f_{1}, \dots, f_{c}$ is a regular sequence in a
 regular local ring $S$ with infinite residue field $k$, so that
 $R= S/(f_{1}, \dots, f_{c})$ is a local complete intersection. Let $M$ be the HMF module  of a minimal higher matrix factorization $(d, h)$ with respect to $f_{1}, \dots, f_{c}$.
Using notation as in  Definition~{\rm \ref{matrixfiltration}}, for $s=0,1$, 
choose splittings $A_{s}(p) = A_{s}(p-1)\oplus B_{s}(p)$ for $s=0,1$, 
so $$A_{s}(p) = \oplus_{1\le q\le p}  B_{s}(q)\,.$$ 
Set
 $B(p)= B_{1}(p)\oplus B_{0}(p)$, where we think of $B_{s}(p)$ as placed in homological degree $s$. There are decompositions
 \b{align*}
\ext_{S}(M,k) &\cong \bigoplus_{p=1}^{c} k\langle e_{1}, \dots, e_{p-1}\rangle \otimes \hom_S(B(p),k)\\ 
\ext_{R}(M,k) &\cong \bigoplus_{p=1}^{c} k[ \chi_{p}, \dots, \chi_{c}] \otimes \hom_S(B(p),k),
\end{align*}
\noindent as vector spaces, where $k\langle e_{1}, \dots, e_{p-1}\rangle$ denotes the exterior algebra on variables of degree $1$ and
$k[ \chi_{p}, \dots, \chi_{c}]$ denotes the polynomial ring on variables of degree $2$.  
\end{cor}
The former formula in \ref{corzero} follows from Remark~\ref{koszulexpression}
 and the latter  from Corollary~\ref{action on ext}.   We  explain in  \cite{EPS1} and Corollary~\ref{action on ext} how the given decompositions  reflect certain natural actions of the exterior and symmetric algebras on the graded modules $\ext_{S}(M,k)$ and $\ext_{R}(M,k)$.

\myspace
\noindent{\bf Syzygies over intermediate quotient rings}
\vglue .15cm
\noindent
For each $0\le p\le c$ set  $R(p):=S/(f_{1}, \dots, f_{p})$. 
In the case of a
codimension $1$ matrix factorization $(d,h)$, one can use the data of the matrix factorization to describe two
minimal free resolutions, as explained in (\ref{twores}). In the case of a codimension $c$ higher matrix factorization we construct the minimal free resolutions of its HMF module  over all $c+1$ rings
$$S=R(0), \ S/(f_{1})=R(1), \,\dots, \,S/(f_{1}, \dots, f_{c}) = R(c)\,.$$ See Theorem~\ref{intermediate} for the 
intermediate cases.

By Definition~\ref{matrixfiltration} an HMF module $M$ with respect to
the regular sequence $f_{1}, \dots, f_{c}$ determines, for each $p\le c$, an HMF  $R(p)$-module $M(p)$
with respect to $f_{1}, \dots, f_{p}$.
In the notation and hypotheses as in Theorem~\ref{highsyzthm},
Corollary~\ref{cosyzcor} shows that 
$$ M(p-1)=\syz_{2}^{R(p-1)}\Big(\cosyz_{2}^{R(p)}\Big( M(p)\Big)\Big)\, ,$$
where  $\syz(-)$  and $\cosyz(-)$ denote syzygy and cosyzygy, respectively.
Furthermore, Corollary~\ref{takingsyz} says that  if we replace
$M$ by its first syzygy, then all the modules $M(p)$
are  replaced by their first syzygies:
$$
\Big(\syz_1^{R(p)}(M(p))\Big)(p-1)=
\syz_1^{R(p-1)}\Big(M(p-1)\Big)\, .
$$
Theorem~\ref{agreement} 
expresses the modules $M(p)$ as syzygies of $Y: = \cosyz_{c+1}^{R}(M)$ over the intermediate rings
$R(p)$ as follows:
$$
\syz_{c+1}^{R(p)}(Y) \cong M(p)\ \ \hbox{for}\ p\ge 0 \, .
$$

 The package CompleteIntersectionResolutions, available from the first author, 
 can compute  in Macaulay2  examples of many of the constructions in this paper. 
 
\section{Notation and  Conventions} \label{SNotation}

Unless otherwise stated, in the rest of the paper {\bf all rings are assumed commutative and Noetherian, and all modules are assumed finitely generated}.

 If $S$ is a local ring with maximal ideal $\mm$ then  a map  of $S$-modules is  called {\it minimal} if its image is contained in $\bf m$ times the target.

 To distinguish a matrix factorization for one element  from the general concept, sometimes we will refer to the former as a
{\it codimension $1$ matrix factorization} or a
 {\it hypersurface matrix factorization}. 

We will frequently use the following notation.

\b{notation}\label{standardnotation}
 A higher matrix factorization 
$$
\Bigl(d: \   A_{1}\to A_{0},\ h:\oplus_{p=1}^{c}A_{0}(p) \to A_{1}\,\Bigr)
$$
with respect to $f_{1}, \dots, f_{c}$ as in Definition~\ref{matrixfiltration}
involves the following data: 
\b{mylist}
\item[$\bullet$] a ring $S$ over which $A_0$ and $A_1$ are free modules;
\item[$\bullet$]  for $1\le p\le c$, the  rings $R(p):= S/(f_{1}, \dots, f_{p})$, and in particular $R=R(c)$;
\item[$\bullet$]  for $s=0,1$, the filtrations  $0= A_{s}(0) \subseteq \dots\subseteq A_{s}(c) = A_{s}$, preserved by $d$;
\item[$\bullet$]  the induced maps 
$$
 A_{0}(p) \ar{h_p} A_{1}(p) \ar{d_p} A_{0}(p);
$$
\item[$\bullet$]  the quotients  $B_{s}(p) = A_{s}(p)/A_{s}(p-1)$ and the projections $\pi_{p}: A_{1}(p) \to B_{1}(p)$;
\item[$\bullet$] 
the two-term complexes induced by $d$:
\b{align*}
&\A(p): \ A_{1}(p) \ar{d_{p}} A_{0}(p)\\ 
&\B(p): \   B_{1}(p) \ar{b_{p}} B_{0}(p)
\end{align*}
\item[$\bullet$]  the modules 
$$
M(p) = \coker \Big(R(p)\otimes d_p: \ R(p)\otimes A_{1}(p) \to R(p) \otimes A_{0}(p)\,\Big),
$$
and in particular, the HMF module $M = M(c)$ of $(d,h)$.
\end{mylist}

\noindent
We sometimes write $h = (h_{1}|\cdots |h_{c})$.
We say that the higher matrix factorization is {\it trivial} if $A_1= A_0=0$.

If $1\le p\le c$ then $d_p$
together with the maps
$h_q$ for $q\leq p$, is a higher  matrix factorization with respect to $f_1,\dots, f_p$; we write it as
$(d_p,h(p))$, where $h(p) = (h_{1}|\cdots|h_{p})$.  We call $(d_1, h_1)$ the {\it codimension $1$ part} of the higher matrix factorization; $(d_1,h_1)$ is
 a hypersurface matrix factorization for $f_1$  over  $S$ (it could be trivial).
If $q\ge 1$ is the smallest number such that $A(q)\not= 0$ and $R'=S/(f_1,\dots ,f_{q-1})$, then writing ${-}'$ for $R'\otimes -$,
the maps 
$$
b_{q}': \ B_{1}(q)'\to B_{0}(q)'\ \ \hbox{and}\ \ h_q': \  B_{0}(q)'\to B_{1}(q)'
$$ 
form a
hypersurface matrix factorization for the element
$f_{q}\in R'$. We call it the {\it top nonzero part} of the higher matrix factorization $(d,h)$.

For each $0\le p\le c$ set  $R(p):=S/(f_{1}, \dots, f_{p})$. 
The HMF module $$M(p)=\coker (R(p)\otimes d_p)$$ is an $R(p)$-module. 
\end{notation}

\myspace
Next, we make some conventions about complexes.

We write $\U[-a]$ for the shifted complex, with
$\U[-a]_i = \U_{i+a}$ and differential $(-1)^{a}d$.

Let $(\W,\partial^{W})$ and $({\bf Y},\partial^{Y})$ be complexes. The complex ${\bf W}\otimes {\bf Y}$ has differential 
$$\partial^{W\otimes Y}_{q} = \sum_{i+j=q}\, \Big((-1)^j\partial^{W}_i\otimes \id+
\id\otimes \partial^{Y}_j\Big)
\, .$$ A map of complexes $\gamma: \W[a]\to \Y$ 
is homotopic to 0 if there exists a map $\alpha: \W[a+1]\to \Y$ such that 
$$
\gamma = \partial^{\Y}\alpha - \alpha \partial^{\W[a+1]} = 
\partial^{\Y}\alpha -(-1)^{a+1} \alpha \partial^{\W}.
$$
If $\varphi: \   {\bf W}[-1] \to {\bf Y}$ is a map of complexes, so that $-\varphi\partial^{W} = \partial^{Y}\varphi$, then the {\it mapping cone}
 $\MM(\varphi)$ is the complex $\MM(\varphi)={\bf Y}\oplus {\bf W}$ with modules $\MM(\varphi)_i=Y_i\oplus W_i$ and differential 
 $$
 \bordermatrix{~ 
 & Y_i &W_i\cr 
 Y_{i-1} & \partial_i^Y & \varphi_{i-1}\cr
 W_{i-1} & 0 &\partial_i^{W}}\, .
 $$

If $f$ is an element in a ring $S$ then we write ${\bf K}(f)$ for the two-term Koszul complex $f: eS\to S$, where we think of $e$ as an exterior variable. If $({\bf W},\partial)$ is any complex of $S$-modules we write
 ${\bf K}(f)\otimes {\bf W}= e{\bf W}\oplus {\bf W} $; it is the mapping cone of the map ${\bf W}\to {\bf W} $
that is $(-1)^if:\ W_i\to W_i$.

\section{The minimal $S$-free resolution of a higher matrix factorization module} \label{Sfinitematr}

We will use the notation in \ref{standardnotation} throughout this section.
Suppose that $M$ is the HMF module of a  higher matrix factorization
$(d,h)$ with respect to a regular sequence
$f_{1},\dots,f_{c}$ in a local ring $S$. Theorem~\ref{fresthm} expresses the  minimal 
$S$-free resolution of $M$  as an iterated mapping cone of {Koszul extensions}, which we will
now define in \ref{mcone}. We say that a complex $(\U,d)$ is a {\it left complex} if  $U_j=0$ for $j<0$; thus for example the free
resolution of a module  is a left complex.

\b{definition}\label{mcone} 
Let $S$ be a  ring.
Let $\B$ and $\L$ be
$S$-free left complexes, and let 
$\psi: \   {\bf B}[-1] \to {\bf L}$ be a map of complexes. 
Note that $\psi$ is zero on $B_{0}$.
Denote ${\bf K} := {\bf K}(f_1,\dots ,f_p)$ the Koszul complex on $f_1,\dots ,f_p\in S$.
An  {\it $(f_{1}\dots,f_{p})$-Koszul extension} of $\psi$ is a map of complexes
$$\Psi:\ \KK\otimes {\bf B}[-1]  \to {\bf L}$$ 
extending 
$$
\KK_{0}\otimes \B[-1] = \B[-1]  \ar{\psi} \L
$$
whose restriction to $\KK\otimes B_{0}$ is  zero.
\end{definition}

The next proposition shows that Koszul extensions exist in the case we will use.

\b{prop} \label{better}  
Let $f_{1}, \dots, f_{p}$ be elements of a  ring $S$.  Let ${\bf L}$ be a free resolution of an $S$-module $N$ annihilated by $f_{1},\dots, f_{p}$.
Let $\psi: \   {\bf B}[-1] \to {\bf L}$ be a map from an $S$-free left complex $\B$.
\b{mylist}
\item[\rm (1)] There exists an $(f_{1}\dots,f_{p})$-Koszul extension of $\psi$.
\item[\rm (2)] If $S$ is local, the elements $f_{i}$ are in the maximal ideal, ${\bf L}$ is  minimal, and the map $\psi$ is minimal, then every Koszul extension of $\psi$ is minimal. 
\end{mylist}
\end{prop}

\proof 
Set $\KK = \KK(f_{1}\dots,f_{p})$, and let
$\varphi: \   \KK \otimes {\bf L} \to {\bf L}$ be any map extending
the identity map $S/(f_{1}, \dots f_{p})\otimes N \to N$. The map
$\varphi$ composed with the tensor product map
$\hbox{\id}_{\KK} \otimes \psi$ is a Koszul extension, proving existence.
For the second statement, note that if $\psi$ is minimal,
then so is the Koszul extension we have constructed. 
Since any two extensions of a map from a free complex to a resolution are homotopic, it follows that every Koszul extension is minimal. 
\qed

\myspace 
  We can now describe our construction of
an $S$-free resolution of an HMF module.

\b{construction}\label{finresconstruction}\label{fres}
 Let $(d,h)$ be a higher matrix factorization with respect to a regular sequence $\ff c $ in a ring $ S$. Using notation as  in~\ref{standardnotation},  we 
choose splittings $A_{s}(p) = A_{s}(p-1)\oplus B_{s}(p)$ for $s=0,1$, 
so $$A_{s}(p) = \oplus_{1\le q\le p}  B_{s}(q)\,$$
and denote by $\psi_p$  the component of $d_p$ mapping $B_{1}(p)$ to  $A_{0}(p-1)$. 

\begin{mylist}
\item[$\bullet$] Set $\L(1) := \B(1)$, a free resolution of $M(1)$ with zero-th term $B_{0}(1) = A_{0}(1)$.
\item[$\bullet$] For $p\geq 2$, suppose that 
$\L(p-1)$ is an $S$-free resolution of $M(p-1)$ with zero-th term 
$L_{0}(p-1) = A_{0}(p-1)$.
Let 
$$
\psi_p': \   \B(p)[-1]\to \L(p-1)
$$
be the map of complexes
induced by $\psi_p: \   B_{1}(p)\to A_{0}(p-1)$, and  let
$$
\Psi_p: \   \KK(f_{1},\dots, f_{p-1})\otimes \B(p) [-1]\to \L(p-1)
$$ be an $(f_1,\dots ,f_{p-1})$-Koszul extension. Set
$
\L(p) = \MM(\Psi_p).
$
\end{mylist}

\noindent The following
theorem implies that $H_{0}(\L(p)) = M(p)$, so that the construction
can be carried through to $\L(c)$.  Note that $\L(c)$ has a filtration with successive quotients of the form
$\KK(f_{1}, \dots, f_{p-1})\otimes \B(p)$. 

\end{construction}

\b{thm} \label{fresthm}
With notation and hypotheses as in {\rm \ref{fres}} the complex $\L(p)$ is an $S$-free resolution of $M(p)$
for $p=1,\dots ,c$.
Moreover, if $S$ is local and $(d,h)$ is minimal, then the resolution $\L(p)$ is minimal.
\end{thm}

\b{remark}\l{koszulexpression}\l{rmvanishing}
The underlying free module of the  Koszul complex ${\bf K}(f_1,\dots ,f_{p-1})$  is the exterior algebra  on generators $e_i$ corresponding to
the $f_i$. Set $B(p) = B_0(p)\oplus B_1(p)\,,$ and thus we get that as an $S$-free module ${\bf L}(p)$ is
$${\bf L}(p)={\bf L}(p-1)\oplus
S\langle e_1,\dots ,e_{p-1}\rangle\otimes_S B(p)\, .$$
 The only non-zero components  of the differential that land in $B_0(p)$ are
  those of the map $d$ and
  $$f_i: \, e_iB_0(p)\to B_0(p)\quad\hbox{ for }\ i<p\, .$$  \e{remark}

\b{example}\label{exampletwo} 
Here is the case of codimension $2$.
After choosing splittings $A_s(2) = B_s(1)\oplus B_s(2)\, ,$ a  higher matrix factorization $(d,h)$  for a regular sequence $f_1,f_2\in S$    is
a diagram of free   $S$-modules

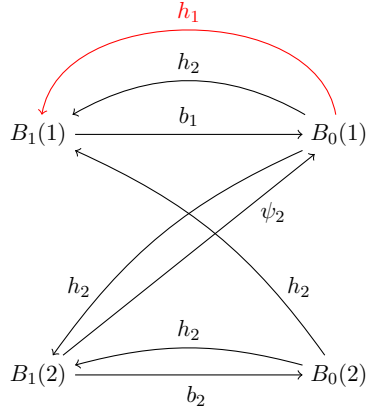
\begin{center}
\begin{tikzpicture}  
 [every node/.style={scale=0.8},  auto]
  \node (TL) {$B_0(1)$};
    \node (TM) [node distance=5cm, left of=TL] {$B_1(1)$};
     \node (M) [node distance=4cm, below of=TM] {$B_1(2)$};
       \node (L) [node distance=4cm, below of=TL] {$B_0(2)$};
   \draw[->] (TM) to node  {$b_1$} (TL);
     \draw[->] (M) to node [below=10pt, right=-5pt] {$b_2$} (L);
      \draw[->] (M) to node [above=20pt, right=30pt]  {$\psi_2$} (TL);
  \draw[->, bend right=30] (TL) to node [swap]{$h_2$} (TM);
   \draw[->, bend right=80, red] (TL) to node [swap]{$h_1$} (TM);
  \draw[->, bend right=15] (L) to node [swap]{$h_2$} (M);
  \draw[->, bend right=15] (TL) to node [below=25pt, left=30pt]{$h_2$} (M);
  \draw[->, bend right=15] (L) to node [below=25pt, right=30pt]{$h_2$} (TM);
\end{tikzpicture}
\end{center}
 \noindent where $d$ has components $b_1,b_2,\psi_2$, and for some $C,D$ we have
\myalign \label{conditions}
\nonumber  b_1h_1&=f_1\Id_{B_0(1)}\ \ \ \hbox{on}\ B_0(1)\\ 
h_1b_1&=f_1\Id_{ B_1(1)}\ \ \   \hbox{on} \ B_1(1)\\ 
\nonumber  dh_2&=f_2\Id+f_1C\ \ \  \hbox{on}\ B_0(1)\oplus B_0(2) \\ 
\nonumber \pi_2 h_2 d_2 &= f_{2}\pi_2+f_1D\pi_2 \ \ \ \hbox{on}\ B_1(1)\oplus B_1(2) \, .
\end{align}
Applying  Theorem~\ref{fresthm}, we may write the $S$-free resolution of the HMF module
$M = \coker(S/(f_1,f_2)\,\otimes d)$ in (\ref{ffour}). The homotopy for $f_1$ is shown with red arrows, and the homotopy for $f_2$ is not shown. 
\eq \label{ffour} \  \end{equation} \vglue -1.4cm \ 
\begin{center}
\begin{tikzpicture}  
 [every node/.style={scale=0.8},  auto]
  \node (TL) {$B_0(1)$};
    \node (TM) [node distance=4cm, left of=TL] {$B_1(1)$};
     \node (M) [node distance=2cm, below of=TM] {$B_1(2)$};
       \node (L) [node distance=2cm, below of=TL] {$B_0(2)$};
       \node (BM) [node distance=2cm, below of=M] {$e_1B_0(2)$};
        \node (BL) [node distance=4cm, left of=BM] {$e_1B_1(2)$};
   \draw[->] (TM) to node  {$b_1$} (TL);
     \draw[->] (M) to node [above=10pt, left=-10pt] {$b_2$} (L);
      \draw[->] (M) to node [above=5pt, right=-20pt]  {$\psi_2$} (TL);
     \draw[->] (BL) to node [swap] {$b_2$} (BM);
  \draw[->] (BL) to node [above=12pt, right=-8pt] {$-f_1$} (M);
   \draw[->] (BL) to node [above=20pt, right=-15pt]  [swap] {$h_1\psi_2$} (TM);
   \draw[<-,  bend right=20, red] (BL) to node [above=15pt, right=25pt] {$h_1=\id$} (M);
 \draw[->, bend right=30, red] (TL) to node [swap]{$h_1$} (TM);
 \draw[<-, bend right=20, red] (BM) to node [above=15pt, right=25pt]{$h_1=\id$} (L);
  \draw[->] (BM) to node [above=10pt, right=-5pt]{$f_1$} (L);
\end{tikzpicture}
  \end{center}
\myspace
\noindent 
\end{example}

Before giving the proof of Theorem~\ref{fresthm} we exhibit  some consequences for the structure of 
modules  that can be expressed as HMF modules. We keep notation
as in \ref{standardnotation}.

\vglue .5cm
\b{cor} \label{IntroSres}
With notation and hypotheses as in  {\rm \ref{fres}}, if in addition 
$S$ is local and the higher matrix factorization is minimal, then
the minimal $S$-free resolution of $M$ has a filtration by minimal $S$-free resolutions of the modules
$M(p):= \coker(S/(\ff{p})\otimes d_p)$, whose 
successive quotients are the complexes
$$
 \KK(f_{1}, \dots, f_{p-1}) \otimes_{S} \B(p).
$$
\vglue -.75cm \rightline{\qed}
\end{cor}
\vskip .3cm

\myspace
\b{cor}\label{bettifinite}
With notation and hypotheses as in  {\rm \ref{fres}}, if in addition $S$ is local and the higher matrix factorization $(d,h)$ is minimal, then
the  Poincar\'e series of the HMF module $M$ of the higher matrix factorization $(d,h)$ is
$$
{\poin}^S_{M}(x)=\sum_{1\le p\le c}\ (1+x)^{p-1} \Big(x\,\rank (B_1(p))+\rank(B_0(p))\Big)\, .$$
\end{cor}

\b{cor} \label{mcm} 
With notation and hypotheses as in {\rm \ref{fres}},  if $M(p)\not= 0$ then its projective
dimension  over $S$ is $p$, and $f_{p+1}$ is a non-zerodivisor on $M(p)$.
If $S$ is  a local Cohen-Macaulay ring
then the module $M(p)$ is a maximal Cohen-Macaulay $R(p)$-module.
\end{cor}

\proof The resolution $\L(p)$ has length $p$, and no module annihilated by a regular
sequence of length $p$ can have projective dimension $<p$. The Cohen-Macaulay
statement  follows from this and the Auslander-Buchsbaum formula.

Suppose that $f_{p+1}$ is a zerodivisor on $M(p)$. Hence, $f_{p+1}$ is contained in a minimal prime $\bf n$ over
$\hbox{ann}_S(M(p))$. Since $f_1,\dots ,f_p$ annihilate $M(p)$, they are
contained in $\bf n$ as well. Therefore, the height of $\bf n$ is $\ge p+1$. The projective dimension of
$M(p)_{\bf n}$ over $S_{\bf n}$ is less or equal to $p$, so it is strictly less than $\dim (S_{\bf n})$.
Thus the minimal $S_{\bf n}$-free resolution of $M(p)_{\bf n}$ is a complex of length $<\dim (S_{\bf n})$
and its homology $M(p)_{\bf n}$ has finite length. This is a contradiction
by the New Intersection Theorem, cf. \cite{PW}.
\qed

\b{cor}\label{filteredFiniteResolution}
With notation and hypotheses as in  {\rm \ref{fres}}, if in addition 
$S$ is local and the higher matrix factorization is minimal, then
$M(p)$ has no $R(p)$-free summands.
\end{cor}

\proof
If $M(p)$ had an $R(p)$-free summand, then with respect to suitable bases the minimal presentation matrix $R(p)\otimes d_p$ of $M(p)$ would have a row of zeros.  
Thus a matrix representing $R(p-1)\otimes d_p$
would have a row of elements divisible by $f_{p}$. Composing with $h_p$ we see that a matrix representing $R(p-1)\otimes d_p h_p$ would have a row of elements in 
$\mm f_{p}$. However $R(p-1)\otimes \left(d_p h_c\right) = f_{p}\Id$, a contradiction.
\qed

\myspace
The following result shows that HMF modules are quite special. Looking ahead to Corollary~\ref{highsyzy}, we see that
it can be applied to {\it any} $S$ module that is a sufficiently high syzygy over $R$.

\b{cor}\label{simplerelations} 
With notation and hypotheses as in  {\rm \ref{fres}}, suppose in addition that $S$ is local and that the higher matrix factorization $(d,h)$ is minimal, and let $n = \sum_{p}\rank B_{0}(p)$, the rank of the target of $d$. In a suitable basis, the minimal presentation matrix of the HMF module $M$ consists of the matrix $d$ concatenated with an $\Big(n\times \sum_p (p-1)\rank B_{0}(p)\Big)$-matrix that is the direct sum of matrices of the form 
$$
\b{pmatrix}
 f_1 & \kern-.5em\dots & \kern-.5em f_{p-1} \\ 
 \end{pmatrix}
\otimes \Id_{B_{0}(p)} =
\b{pmatrix}
f_{1}&\cdots &f_{p-1}&0&\cdots&0&\cdots&0&\cdots&0\\ 
0&\cdots&0&f_{1}&\cdots &f_{p-1}&\cdots&0&\cdots&0\\ 
0&\cdots&0&0&\cdots &0&\cdots&0&\cdots&0\\ 
\vdots&\vdots&\vdots&\vdots&\vdots&\vdots&\vdots&\vdots&\vdots&\vdots\\ 
0&\cdots&0&0&\cdots &0&\cdots &f_{1}&\cdots&f_{p-1} 
\end{pmatrix}
$$
\end{cor}

We remark that a similar property holds for all matrices of the differential in the minimal free resolution of $M$.

\myspace
\proof In the notation of
Construction~\ref{fres}, the given direct sum is the part of the map $\L_{1}(c)\to \L_{0}(c)$ that corresponds to
$$
\oplus_{p} \Big(\K(f_{1}, \dots, f_{p-1})\Big)_{1}\otimes B_{0}(p)\to \oplus_{p}B_{0}(p).
$$
\vglue -.75cm
\ \qed

\myspace\vglue .1cm
Theorem~\ref{fresthm} and Corollary~\ref{filteredFiniteResolution} allow us to express the Betti numbers of an HMF module in terms of the ranks of the modules $B_{s}(p)$. Recall that if $S$
is a local ring with residue field $k$ then the
 {\it Betti numbers} of a  module $N$ over $S$ are
$b_i^S(N)=\dim_k(\tor_i^S(N,k))$. They are often studied via the {\it Poincar\'e series}:
$$
{\poin}^S_N(x)=\sum_{i\ge 0} \,b_i^S(N)\,x^i \, .
$$

Corollary~\ref{bettifinite} makes it worthwhile to ask whether there are interesting restrictions on the ranks of the 
$B_{s}(p)$. Here is a first result in this direction:

\b{cor}\label{nozeros} 
With notation and hypotheses as in  {\rm \ref{fres}}, suppose in addition that $S$ is local and Cohen-Macaulay and that the higher matrix factorization $(d,h)$ is minimal. If $B_{1}(p) = 0$ for some $p$,  then 
$B_{1}(q)=B_{0}(q) = 0$ for all  $q\leq p$.
\end{cor}

\proof Suppose that $B_{1}(p)=0$.
If $B_{0}(p)\neq 0$ then $M(p)$ would have a free summand, 
contradicting Corollary~\ref{filteredFiniteResolution}, so $B_{0}(p) = 0$ as well.
 It 
follows that $h_p$ restricts to a map $A_{0}(p-1)\to A_{1}(p-1)$, and thus
$M(p-1)$ is annihilated by $f_{p}$. However, if $M(p-1)\neq 0$ then
by Corollary~\ref{mcm} it would be
a  maximal Cohen-Macaulay module over the ring $R(p-1)$, and this is 
a contradiction. Thus $M(p-1) = 0$, so $B_{s}(q) = 0$ for $q\leq p$.
\qed

\b{example}\label{unstablemf}
Let $S=k[x,y,z]$ and let $f_{1}, f_{2}$ be the regular sequence $xz, y^{2}$. We give an example of a higher matrix factorization with respect to $f_{1}, f_{2}$ such that
$B_{1}(2) \neq 0,$ but $B_{0}(2) = 0.$ 
If

\myspace
\begin{center}
\begin{tikzpicture}  
 [every node/.style={scale=0.8},  auto]
  \node (TL) {$B_0(1)=S^2$};
    \node (TM) [node distance=6cm, left of=TL] {$B_1(1)=S^2$};
     \node (M) [node distance=2.5cm, below of=TM] {$B_1(2)=S $};
       \node (L) [node distance=2.5cm, below of=TL] {$B_0(2)=0\, ,$};
   \draw[->] (TM) to node  {$\b{pmatrix}z&-y\\  0&x\end{pmatrix}$} (TL);
     \draw[->] (M) to node [above=10pt, right=-1pt] {$0$} (L);
      \draw[->] (M) to node [above=10pt, left=3pt]  {$\b{pmatrix}0\\  y\end{pmatrix}$} (TL);
\end{tikzpicture}
\end{center}
\myspace

\noindent
and 
$$
h_1 = \b{pmatrix} x&y\\  0&z \end{pmatrix} 
\ \  \hbox{and}\ \ 
h_2 = \b{pmatrix}0&0\\  -y&0 \\  x&y\end{pmatrix},
$$
then $(d,h)$ is a higher matrix factorization.
\end{example}

In the case of higher matrix factorizations that come from high syzygies (stable matrix factorizations) Corollary~\ref{nozeros} can be strengthened further:  $B_{0}(p) = 0$ implies $B_{1}(p) = 0$ as well; see Corollary~\ref{zerostop}. This is not the case in general, as the above example shows.

\myspace
\noindent{\prooffont Proof of Theorem~\ref{fresthm}:}
The minimality statement follows at once from the 
construction and  Proposition~\ref{better}(2).
Thus it suffices to prove the first statement.

Note that $d_1=b_1$.
The equations in the definition of a higher matrix factorization
imply in particular that $h_1b_{1} = b_{1}h_1 = f_{1}\id$,
so $b_{1}$ is a monomorphism. Note that  $\coker (d_1)$ is annihilated
by $f_{1}$. Thus $\L(1) = \B(1)$ is an $S$-free resolution of 
$$
M(1) = \coker \Big(R(1)\otimes d_1\Big) = \coker (d_1).
$$

To complete the proof we do induction on $p$. By induction hypothesis
$$
\L(p-1): \quad    \cdots\to L_{1}(p-1) \to L_{0}(p-1)
$$
is a free resolution of $M(p-1)$. Since $L_{0}(p-1) = A_{0}(p-1)$,
 the map $\psi_p$ defines a morphism of complexes
$\psi_p': \   \B(p)[-1] \to \L(p-1)$ and thus a mapping cone
\myspace
\begin{center}
\begin{tikzpicture}  
 [every node/.style={scale=1},  auto]
\node(00){$L_{2}(p-1)$};
\node(pre)[node distance = 1.8cm, left of=00]{};
\node(prepre)[node distance = .5cm, left of=pre]{$\cdots$};
\node(preprepre)[node distance=1.5cm, left of=prepre]{};
\node(preprepre1)[node distance=.7cm, below of=preprepre]{};
\node(preprepre2)[node distance=.7cm, below of=preprepre1]{};
\node(01)[node distance=2.5cm, right of=00]{$L_{1}(p-1)$};
\node(02)[node distance=3.5cm, right of=01]{$L_{0}(p-1)\ $};
\node(11)[node distance=1.4cm, below of=01]{$B_{1}(p)$};
\node(12)[node distance=1.4cm, below of=02]{$B_{0}(p)\,.$};
\draw[->] (pre) to node {} (00);
\draw[->] (00) to node  {} (01);
\draw[->] (01) to node [ pos=.4]{} (02);
\draw[->] (11) to node [swap,pos=.6]{$b_{p}$} (12);
\draw[->] (11) to node [above=3pt, left=3pt] {$\psi_p$} (02);
\end{tikzpicture}
\end{center}
\vglue .1cm

To simplify the notation, denote by $\K$  the Koszul complex $\K(
f_{1},   \dots,  f_{p-1})$
of $f_{1},   \dots, \allowbreak  f_{p-1}$, and write $\kappa_{i}: \   \wedge^{i}S^{p-1}\to \wedge^{i-1}S^{p-1}$
for its differential.
Also, set $B_{s} := B_{s}(p)$ and $\B:\ B_1\ar{b_p} B_0$. 

Since $M(p-1)$ is annihilated by
$(f_{1}, \dots, f_{p-1})$, Proposition~\ref{better} shows that there exists a Koszul extension  $\Psi_p: \   \K\otimes \B[-1]\to \L(p-1)$ of $\psi_p'$. Let $(\L(p),\epsilon )$
be the mapping cone of $\Psi_p$, and note that the zero-th terms of $\L(p)$ is
$L_{0} = L_{0}(p-1)\oplus B_{0} = A_{0}(p)$. 
We will show that $\L(p)$ is a resolution of $M(p)$. 

We first show that $H_{0}(\L(p)) = \coker ( \epsilon_{1}) = M(p).$ 
If we drop the columns corresponding to $B_{1}$
from a matrix for $\epsilon_{1}$ we get a presentation of $M(p-1)\oplus (R(p-1)\otimes B_0(p-1))$, 
so $\coker( \epsilon_{1})$ is annihilated by $(f_{1}, \dots, f_{p-1})$. 
Moreover, the map $h_p: A_{0}(p) \to A_{1}(p) \subset L_{1}(p)$
defines a homotopy for multiplication by $f_{p}$ modulo $(f_{1}, \dots, f_{p-1})$,
so  $\coker( \epsilon_{1})$ is annihilated by $f_{p}$ as well. Thus
$\coker(\epsilon_{1} )= \coker(R(p)\otimes \epsilon_{1}) = M(p)$ as required.

We next analyze the homology of the complex $\K\otimes \B$. It is isomorphic to
$\B\otimes \K$, which is
the mapping cone of the map 
$$(-1)^ib_{p}\otimes \id:\ B_{1}[-1]\otimes K_i\to B_{0}\otimes K_i \,,$$ so there is a long
exact sequence 
$$
\cdots \to H_{i}(\K\otimes B_{1}) \to H_{i}(\K\otimes B_{0})\to  H_{i}(\K\otimes \B) \to H_{i-1}(\K\otimes B_{1}) \to \cdots.
$$
Since $\K\otimes B_{s}$ is a resolution of $R(p-1)\otimes B_s$
we see that $H_{i}(\K\otimes \B) = 0$ for $i>1$ and there
is a four-term exact sequence
$$
0 \to H_{1}(\K\otimes \B) \to R(p-1) \otimes B_{1}
{\mathop{{\mathrel{\smash-}}{\mathrel{\mkern-8mu}}\a\a\a
{\mathrel{\smash-}}{\mathrel{\mkern-8mu}}\rightarrow}\limits^{{R(p-1)\otimes b_{p}}}}
\,  R(p-1) \otimes B_{0} \to H_{0}(\K\otimes \B) \to 0.
$$

Since $\L(p)$ is the mapping cone of $\Psi_p$, we have a long exact sequence in homology
of the form
$$
\cdots  \to H_{i}(\L(p-1)) \to H_{i}(\L(p)) \to H_{i}(\K\otimes \B)\ar{{\Psi_p}_{*}} H_{i-1}(\L(p-1))\to \cdots\, ,
$$
so from the vanishing of the $H_{i}(\K\otimes \B)$ for $i>1$ we see that $H_{i}(\L(p)) = 0$ for $i>1$.

It remains to prove only that $H_{1}(\L(p)) = 0$, or equivalently that the map
$${\Psi_p}_{*}: \   H_{1}(\K\otimes \B)\to H_{0}(\L(p-1)) = M(p-1)$$ is a monomorphism. From the four-term
exact  sequence above we see that 
$$H_{1}(\K\otimes \B) =\ker\Big(R(p-1) \otimes b_{p}\Big)\, . 
$$
Also note that by construction the map 
$$
{\Psi_p}_{*}: \ker\Big(R(p-1) \otimes b_{p}\Big) \to H_{0}(\L(p-1)) = \coker\Big(R(p-1) \otimes d_{p-1}\Big)
$$
is induced by $$\psi_{p}: R(p-1) \otimes B_{1}(p) \to R(p-1) \otimes A_{0}(p-1)\, .$$
Since $\overline{L}_{0}(p-1)=\overline{A}_{0}(p-1)$, the  proof is finished by the next Lemma~\ref{newex}, 
which we will use again in Section~\ref{Stoinf}.
\qed

\b{lemma}\label{newex}
With notation and hypotheses as in Construction~ {\rm \ref{fres}},
$\psi_{p}$ induces a monomorphism from
$\ker\Big(R(p-1) \otimes b_{p}\Big)$
to
$\coker\Big(R(p-1) \otimes d_{p-1}\Big)$.
\end{lemma}

\proof
To simplify notation we write $\overline{-}$ for $R(p-1)\otimes -$. Consider the diagram:

\myspace
\begin{center}
\begin{tikzpicture}  
 [every node/.style={scale=0.9},  auto]
\node(01) {$u\in \overline{A}_{1}(p-1)$};
\node(02)[node distance=6cm, right of=01]{$\overline{A}_{0}(p-1)\ $};
\node(11)[node distance=2cm, below of=01]{$v\in\overline{B}_{1}(p)$};
\node(12)[node distance=2cm, below of=02]{$\overline{B}_{0}(p)\,.$};
\draw[->] (01) to node {$\overline{d}_{p-1}$} (02);
\draw[->] (11) to node [swap]{$\overline{b}_{p}$} (12);
\draw[->] (11) to node [above=5pt, left=5pt] {$\overline{\psi}_p$} (02);
\end{tikzpicture}
\end{center}

\noindent 
We must show that 
if $v\in \ker(\overline{b}_p)$ and 
$\overline{\psi}_p(v)=\overline{d}_{p-1}(u)$ for some $u\in \overline A_{1}(p-1)$,
then $v=0$. 

Write $\overline{\pi}_p$ for the projection of $\overline A_1(p) = \overline A_{1}(p-1)\oplus \overline B_1(p)$ to 
$\overline B_1(p)$,
and note that $\overline d_{p}$ is the sum of the three maps in the diagram above.
Our equations say that $d_{p}(-u,v) = 0$. 
By condition (b) in Definition~\ref{matrixfiltration}, 
$$
f_pv= f_{p}\overline{\pi}_p(-u,v) = \overline{\pi}_p\overline{h\mathstrut}_p\,\overline d_p(-u,v)=0.
$$
Since $f_p$ is a non-zerodivisor in $R(p-1)$, it follows that $v=0$.
\qed

\vglue -.5cm
\section{Resolutions with a surjective CI operator} \label{Scitools}

We begin by recalling the definition of CI operators.
Suppose that $\ff c \in S$ is a regular sequence and $({\bf V},\partial  )$ is a complex of   free modules over $R = S/(\ff c )$. Suppose that $\widetilde\V$ is a lifting of $\V$ to $S$, that is,
a sequence of free modules $\widetilde {V}_{i}$ and maps $\widetilde{\partial}_{i+1}: 
\widetilde V_{i+1}\to \widetilde {V}_{i}$ such that
$\partial = R\otimes \widetilde \partial$. Since $\partial^{2} = 0$ we can choose maps
$\widetilde {t_{j}}: \  \widetilde {V}_{i+1}\to \widetilde {V}_{i-1}$, where $1\le j\le c$,  such that 
$$
\widetilde\partial^{2} = \sum_{j=1}^{c}f_{j}\widetilde{t}_{j}.
$$
We set
\vglue -.5cm
$$
t_{j}:= R\otimes \widetilde {t_{j}}.
$$
Since 
$$
\sum_{j=1}^{c} {f \mathstrut}_{j}\widetilde{t}_{j}\,\widetilde {\partial}= \widetilde {\partial\mathstrut}^{3} 
= \sum_{j=1}^{c} f_{j}\widetilde {\partial}\,\widetilde {t}_{j},
$$ 
and the $f_{i}$ form a regular sequence,
we see that each $t_{j}$ commutes with $\partial$,
and thus the $t_{j}$ define a map of 
complexes $\V[-2]\to \V$,  \cite[1.1]{Ei1}. 
In the case $c=1$, we have $\widetilde{\partial}^2=f_1\widetilde{t}_1$ and we sometimes
write $\widetilde{t}_1={1\over f_1}\widetilde{\partial}^2$ and call it the {\it lifted CI operator}.

\cite[1.2 and 1.5]{Ei1} shows that the operators $t_{j}$ are, up to homotopy,
independent of the choice of liftings.
They are called  the {\it  CI operators} (sometimes called Eisenbud operators)
associated to the sequence $f_{1}, \dots, f_{c}$.

We next recall the  definition of higher homotopies and the Shamash construction. The version for a single element  is due to Shamash \cite{Sh};  \cite{Ei2} treats the more general case of a collection of
elements.

\b{definition}\label{hhomotop} 
Let $f_1,\dots ,f_c\in S$, and ${\bf G}$ be a  free complex of   $S$-modules. We denote ${\bf a}=(a_1,\dots ,a_c)$, where each $a_i\ge 0$ is an integer, and set $|{\bf a}|=\sum_ia_i$.
A  {\it  system of higher homotopies}  $\sigma$  for $f_1,\dots ,f_c$ on ${\bf G}$ is 
a collection of maps $$\sigma_{{\bf a}}: \   \G\to \G[-2|{\bf a}|+1]$$  of the underlying modules such that 
the following three conditions are satisfied:
\b{mylist}
\item[(1)] $ \sigma_{\bf 0}$ is the differential on ${\bf G}$.
\item[(2)] For each $1\le i\le c$, the map $ \sigma_{\bf 0}\sigma_{{\bf e}_i}+ \sigma_{{\bf e}_i}\sigma_{\bf 0}$ is multiplication by $f_i$ on $\G$, where ${\bf e}_i$ is the $i$-th standard vector.
\item[(3)] If ${\bf a}$ is a multi-index with  $|{\bf a}|\ge 2$,  then  $\sum_{{\bf b}+{\bf s}= {\bf a}}  \sigma_{\bf b}\sigma_{{\bf s}}=0$.
\end{mylist}
\noindent A   {system of higher homotopies}  $\sigma$  for one element $f\in S$ on ${\bf G}$ consists 
 of maps $\sigma_{j}: \   \G\to \G[-2j+1]$ for $j=0,1,\dots$, and will be denoted $\{\sigma_j\}$.
\end{definition}

\b{prop}\label{homotopyproof}
{\rm \cite{Ei2, Sh}}
If  ${\bf G}$ is a free resolution of an $S$-module 
annihilated by elements  $f_1,\dots ,f_c\in S$, then there exists 
 a  system of higher homotopies on $\G$ for $f_1,\dots ,f_c$. \end{prop}
 
For the reader's convenience we present a short proof following \cite{Sh}:

\myspace 
 \proof
 It is well-known that homotopies $\sigma_{{\bf e}_i}$ satisfying (2) in Definition~\ref{hhomotop}  exist.
Equation (3) in \ref{hhomotop} can be written as
$$d\sigma_{\bf a}=-\sum_{{\bf b}+{\bf s}={\bf a}\atop {\bf b}\neq {\bf 0}}\  \sigma_{\bf b}\sigma_{\bf s}\,.$$ As $\G$ is a free resolution,  in order to show by induction  on ${\bf a}$ and on the homological degree that the desired
$\sigma_{\bf a}$ exists,
it suffices to show that the right-hand side is annihilated by $d$. Indeed,

\b{align*}
&-\sum_{{\bf b}+{\bf s} ={\bf a}\atop {\bf b}\neq  {\bf 0}}\  (d\sigma_{\bf b})\sigma_{\bf s}
=\sum_{{\bf b}+{\bf s}={\bf a}\atop {\bf b}\neq  {\bf 0}}\, \sum_{{\bf m}+{\bf r}={\bf b}\atop {\bf r}\neq  {\bf 0}}\,
\sigma_{\bf r}\sigma_{\bf m}\sigma_{\bf s}-\sum_{\{i:\,{\bf e_i}<{\bf a}\}}\,f_i\sigma_{{\bf a}-{\bf e_i}}
\\  &=
\sum_{{\bf m}+{\bf r}+{\bf s}={\bf a}\atop {\bf r}\neq  {\bf 0}}\, \sigma_{\bf r}\sigma_{\bf m}\sigma_{\bf s}-\sum_{\{i:\,{\bf e_i}<{\bf a}\}}\,f_i\sigma_{{\bf a}-{\bf e_i}}
\\  &=
-\sum_{\{i:\,{\bf e_i}<{\bf a}\}}\,f_i\sigma_{{\bf a}-{\bf e_i}}+
\sum_{{\bf r}\neq  {\bf 0}}\sigma_{\bf r}\Bigg(\sum_{{\bf m}+{\bf s}={\bf a}-{\bf r}}\sigma_{\bf m}\sigma_{\bf s}\Bigg)\\  &=
\sum_{{\bf r}\neq  {\bf 0}\atop {\bf r}\neq {\bf a}-{\bf e_i}}\sigma_{\bf r}\Bigg(\sum_{{\bf m}+{\bf s}={\bf a}-{\bf r}}\sigma_{\bf m}\sigma_{{\bf s}}\Bigg)
+\sum_{\{i:\,{\bf e_i}<{\bf a}\}}\, \sigma_{{\bf a}-{\bf e_i}}(\sigma_{\bf e_i}\sigma_ {\bf 0}+\sigma_ {\bf 0}\sigma_{\bf e_i}-f_i)
=0\,,
\end{align*}
where the first and the last equalities hold by induction hypothesis.
\qed

\b{construction}\label{standard} (cf. \cite[Section 7]{Ei1})
Suppose that $f_1,\dots ,f_c$ are elements  in a ring $S$, and that
$ {\bf G}$ is a  free complex over $S$  with  a system $\sigma$ of higher homotopies. This gives rise to
a new complex $\std(\G, \s)$. To define it, we will write $S\{y_1,\dots ,y_c\}$ for the divided
power algebra over $S$ on variables $y_1,\dots ,y_c$; thus, 
$$
S\{y_1,\dots ,y_c\} \cong \hom_{\hbox{graded $S$-modules}}(S[t_1,\dots ,t_c], S) =  { \oplus}
\,Sy_1^{(i_1) }\cdots y_c^{(i_c)}
$$
where the $y_1^{(i_1) }\cdots y_c^{(i_c) }$ form the dual basis to the monomial basis  of the polynomial ring
$S[t_1,\dots,t_c]$. We  
will use
the fact that $S\{y_1,\dots ,y_c\}$ is an $S[t_1,\dots, t_c]$-module with action $t_jy_j^{(i)} = y^{(i-1)}_j$ (see
\cite[Appendix 2]{Ei3}).

Set $R=S/(f_1,\dots ,f_c)$.
The graded module 
$$ 
S\{y_1,\dots ,y_c\}\otimes  {\bf G}\otimes R,
$$ 
where  each $y_i$ has degree $2$, becomes a free complex over $R$ when
equipped with the
differential
$$
\delta := \sum
\,t^{\bf a}\otimes  \sigma_{\bf a}\otimes  R.
$$ 
This complex is called the {\it Shamash complex}  and denoted
$\std ({\bf G}, \sigma)$.

 In the case when we consider only one element $f\in S$, we denote the divided power algebra
by $S\{y\}$, where  the $y^{(i) }$ form the dual basis to the basis $t^i$ of the polynomial ring
$S[t]$. 
\end{construction}

\b{prop}\label{standardthm}
 {\rm \cite{Ei1, Sh}} 
Let $f_1,\dots ,f_c$ be a regular sequence in a ring $S$, and let $N$ be a module
over $R := S/(f_1,\dots ,f_c)$. If $\G$ is an $S$-free resolution of $N$ and $\sigma$ is
a system of higher homotopies for $f_1,\dots ,f_c$ on $\G$,
then   $\std ({\bf G},\sigma)$ is an $R$-free resolution of $N$. 
\end{prop}

\begin{construction}\label{operators on Ext}
In \cite[1.2 and 1.5]{Ei1} Eisenbud shows that the CI operators are, up to homotopy,
independent of the choice of liftings, and also that they commute up to homotopy.
If $S$ is local
with maximal ideal $\mm$ and residue field $k$, and  $\V$ is 
 an $R$-free resolution of an $R$-module
$N$, then the CI operators $t_{j}$ induce well-defined, commutative maps $\chi_{j}$
on $\ext_{R}(N,k)$, and
thus make $\ext_{R}(N,k)$ into a module over the polynomial ring
 $\RR:=k[\chi_{1}, \cdots, \chi_{c}]$, where the variables $\chi_{j}$ have degree $2$.
The $\chi_{j}$ are also called
CI operators. By \cite[Proposition 1.2]{Ei1}, the action of $\chi_{j}$ can be defined using any CI operators on any 
$R$-free resolution of $N$.
Because the $\chi_j$ have degree $2$,
we may split  any $\RR$-module into
even degree and odd degree parts; in particular, we write
$$
\ext_{R}(N,k) = \ext_{R}^{even}(N,k) \oplus \ext_{R}^{odd}(N,k).
$$
\end{construction}

A version of the following result was first proved in [Gu] by Gulliksen, who used a different construction of operators on Ext. Other constructions of operators were introduced and used by Avramov \cite{Av}, Avramov-Sun \cite{AS}, Eisenbud \cite{Ei1}, and Mehta \cite{Me}. The relations between these constructions were explained by Avramov and Sun \cite{AS}. 
We will use only the construction from \cite{Ei1} outlined at the beginning of this section.  Using that construction, we provide a new  and short proof of the following result.

\b{thm}\label{extGeneration} {\rm \cite{AS, Ei1, Gu}}
Let $f_{1}, \dots, f_{c}$ be a regular sequence in a
  local ring $S$ with residue field $k$, and set
 $R= S/(f_{1}, \dots, f_{c})$.
If $N$ is an $R$-module with finite projective dimension over $S$, then the action of the
CI operators makes $\ext_{R}(N,k)$ into a finitely generated
$\RR:= k[\chi_{1}, \dots, \chi_{c}]$-module.
\end{thm}

\proof
Let $\G$ be a finite $S$-free resolution of $N$. By Proposition~\ref{homotopyproof}, there exists 
 a  system of higher homotopies on $\G$. Proposition~\ref{standardthm} shows that  $\std ({\bf G},\sigma)$ is an $R$-free resolution of $N$. Consider its dual.
 By \cite[Theorem 7.2]{Ei1} (also see Construction~\ref{standardCI}), the CI operators can be chosen to act on 
 $\std ({\bf G},\sigma)$ as
 multiplication by the variables, and thus they commute. By the construction of the Shamash resolution,
 it is clear that $\hbox{Hom}_R(\std ({\bf G},\sigma),k)  $ is a finitely generated module over $\RR$. As the CI operators commute with the differential, it follows that both the kernel and the image of the differential are submodules, so they are finitely generated as well. Thus, so is the quotient module $\ext_{R}(N,k)$.
 \qed

\myspace
In this paper we will  use higher homotopies and the Shamash construction for one element $f\in S$.
We focus on that case in the rest of the section.

\b{construction}\label{standardCI}
Suppose that $f\in S$, and that
$ ({\bf G},\partial) $ is a  free complex over $S$  with  a system $\sigma$ of higher homotopies.
We use the notation in Construction~\ref{standard}. 
The {\it standard lifting} $\tilde \std(\G,\sigma)$ of the Shamash complex to $S$ is $ S\{y\}\otimes  {\bf G}$ with
the maps
$\widetilde \delta = 
\sum
\,t^{\,j}\otimes  \sigma_j.$
In particular,
$\widetilde \delta \Big\vert_{\G}=\partial$, so of course $\tilde \delta^{2}\Big\vert_{\G}=\partial^{2}=0$.
Moreover, the equations of Definition~\ref{hhomotop} say precisely that, 
$\tilde \delta^{2}$ acts on the complementary summand $\G' = \oplus_{i>0}y^{(i)}\G$
by $ft$; that is, it sends each $y^{(i)}\G$ isomorphically to $fy^{(i-1)}\G$.
Thus
$$\widetilde{ \delta}^{2} = f {t\otimes 1}\, .$$
The {\it standard CI operator} for $f$
on $\std(\G, \s)$ is $t\otimes 1$. 
Note that $t: \  \std(\G,\s) \to  \std(\G,\s)[2]$ is surjective, and
is split by the map sending $y^{(i)}u\in S\{y\}\otimes \G\otimes  S/(f)$ to $y^{(i+1)}u$.
Also, the {\it standard lifted CI operator} $$\tilde t:=t\otimes 1:\ \tilde \std(\G,\sigma) \to \tilde\std(\G,\sigma)$$ commutes with the
lifting  $\tilde \delta = \sum\,t^{\,j}\otimes  \sigma_j$ of the differential $\delta$.
\end{construction}

We will use the  following modified version of Proposition~\ref{standardthm}:

\b{prop} \label{second} 
Let $\widetilde{\G}$ be a complex of $S$-free modules with 
a system of higher homotopies $\sigma$ for a non-zerodivisor $f$ in a ring $S$. 
If $\F=\std(\widetilde{\G},\sigma)$, then 
$H_{j}(\F) = 0$ for all $0<j\leq i$ if and only if  $H_{j}(\widetilde{\G}) = 0$ for all
$j\leq i$. In particular, $\std(\widetilde{\G},\sigma)$ is an $S/(f)$-free resolution of a module $N$ if and only if
$\widetilde{\G}$ is an $S$-free resolution of $N$. 
\end{prop}

\proof 
We first show that (without any exactness hypothesis) $H_{0}(\widetilde{\G}) = H_{0}(\F)$.
Since the standard lifted CI operator $\widetilde t: \   \widetilde{F_{i}}\to \widetilde{F}_{i-2}$ is surjective, $f$ annihilates $N:=\coker\Big( \widetilde{\delta}:  \widetilde{F}_{1}\to \widetilde{F}_{0}\Big)$, and thus $N = \coker\Big( \delta:   F_{1}\to F_{0} \Big)= H_{0}(\F)$.
But for $i\leq 1$ we have $\widetilde{G}_{i} = \widetilde{F_{i}}$, so $H_{0}(\widetilde{\G}) = H_{0}(\F)$ as required.

Set $\overline{\G}=R\otimes \widetilde{\G}$.
We now use the short exact sequences of complexes
\b{align*}
&0\longrightarrow \overline\G \longrightarrow \ \F \ \ar{t}\ 
 \F [2]\longrightarrow  \ 0 \\       
       &0\longrightarrow \widetilde{\G} \ar{f}
\ \widetilde{\G} \longrightarrow \overline \G \longrightarrow \ 0\, ,           
\end{align*}
which yield  long exact sequences  
\eq   \label{eqstar}
\cdots\to H_{j-1}(\F) \to H_{j} (\overline\G)\to H_{j}(\F) \to H_{j-2}(\F) \to H_{j-1} (\overline\G) \to 
\cdots 
\end{equation} 
\eq \label{eqlong}
\cdots \to H_{j+1}( \overline\G)\to H_{j}( \widetilde{\G}) \ar{f} H_{j}(  \widetilde{\G}) \to H_{j}(\overline\G)\to H_{j-1}( \widetilde{\G}) \to \cdots
\end{equation}
respectively.
Since $\sigma_{1}$ is a homotopy for $f$ on $\widetilde{\G}$, the latter sequence breaks up into short
exact sequences
\eq \label{starstar}
0\to H_{j}(\widetilde{\G}) \to H_{j}(\overline\G) \to H_{j-1}(\widetilde{\G}) \to 0\, .
\end{equation}

First, assume that $H_{j}(\F) = 0$ for $1\le j\leq i$.
From the long exact sequence (\ref{eqstar})
 we conclude that $H_{j}(\widetilde{\G})= 0$ for $2\leq j\leq i$,
and then  (\ref{starstar}) implies that $H_{j}(\widetilde \G) = 0$ for $1\leq j\leq i$.

Conversely, suppose that $H_{j}(\widetilde \G) = 0$ for $1\leq j\leq i$.  It is well known that if we apply the Shamash construction to a resolution then we get
 a resolution, but since the bound $i$ is
not usually present we give an argument:

Assume that $H_{j}(\widetilde{\G}) = 0$ for $1\le j\leq i$.
By (\ref{starstar}) it follows that $H_{j}(\overline\G)= 0$ for $2\leq j\leq i$. Applying (\ref{eqstar}), 
 we conclude that
$H_j(\F)\cong H_{j-2}(\F)$ for $3\le j\le s$. Hence, it suffices  to prove that
$H_1(\F)=H_2(\F)=0$.

We will prove that $H_{1}(\F)=0$.
Let $g_1\in \widetilde{G}_1$ be an element that reduces modulo $f$ to $\overline g_{1}$.
We have 
$$
\widetilde{\partial} (g_1)=fg_{0} = \widetilde{\partial } \sigma_1 (g_{0})
$$ 
for some $g_{0}\in G_0$. Thus
$ g_{1}- \sigma_1 (g_{0})\in\ker(\widetilde{\partial} )$ is a cycle in $\widetilde{\G}$. Since $H_{1}(\widetilde{\G}) = 0$, 
we must have 
$g_{1}- \sigma_1 (g_{0})=\widetilde{\partial } (g_{2})$ for some $g_{2}\in \widetilde{G}_2$. Using the
isomorphism $\widetilde {F_{2}} = \widetilde{G}_{2}\oplus \widetilde{G}_{0}$ we see that
$$
g_{1}= \sigma_1 (g_{0})+\widetilde{\partial } (g_{2}) = \widetilde \delta(g_{0}+g_{2}).
$$
  It follows that $\overline g_{1} = \delta(\overline g_{0}+\overline g_{2})$ is a boundary in $\F$, 
as required.

Finally, we show that  $H_{2}(\F)=0$. Part of  
(\ref{eqstar})
is the exact sequence
$$
H_{2} (\overline \G) \to H_{2}(\F) \to H_{0}( \F) \ar{\beta} H_{1} (\overline \G)\to H_1(\F).
$$
Since $H_{2}(\overline \G) = 0$,  it suffices to show that the map marked $\beta$ is
a monomorphism. But we already showed that $H_{1}(\F)=0$, so $\beta$ is an
epimorphism. Since its source
and target are isomorphic finitely generated modules over the  ring $S$, this implies that it is an isomorphism, whence $H_{2}(\F) = 0$.
\qed

\myspace
It follows from Theorem~\ref{extGeneration} that CI operators on the resolutions of high syzygies over complete intersections are often surjective, in a sense we will make precise. To prepare for the 
study of this situation, we consider what can be said when a CI operator is  surjective.

\b{prop}\label{first}
Let $f\in S$ be a non-zerodivisor in a  ring $S$, and let
$$
(\F,\delta):\quad \cdots \to F_i\ar{\delta_i}F_{i-1}\to\ \dots\  \to F_1\ar{\delta_1}F_0
$$ 
be a   complex of free $R:=S/(f)$-modules.
Let $(\widetilde \F, \widetilde \delta)$ be a lifting of $(\F,\delta)$ to $S$. Set
\b{align*}
\widetilde t :&= (1/f)\widetilde \delta^{2}: \   \widetilde \F\to \widetilde \F[2]\, ,
\\ 
\widetilde{\G} &= \ker(\,\widetilde {t}\,)\, .
\end{align*}
Suppose that $\widetilde t$ is surjective. Then: 
\b{mylist}
\item[\rm (1)] {\rm \cite[Theorem 8.1]{Ei1}}
The maps $\widetilde{\delta}: \   \widetilde {F}_{i}\to \widetilde {F}_{i-1}$ induce maps $\widetilde{\partial}: \   \widetilde{G}_{i}\to \widetilde{G}_{i-1}$, and
$$
\widetilde{{\bf G}}: \quad\cdots \to \widetilde{G}_{i+1}\ar{\widetilde{\partial}_{i+1}} \widetilde{G}_{i}\to\cdots  \to \widetilde{G}_{1} \ar{\widetilde{\partial}_1} \widetilde{G}_{0}
$$
is an $S$-free complex. If $S$ is local and $\F$ is minimal, then
so is $\widetilde{\G}$.

\item[\rm (2)] We may write 
$\widetilde {F}_{i} = \oplus_{j\geq 0}\widetilde{G}_{i-2j}$
 in such a way that the lifted CI operator $\widetilde t$ consists of 
the projections 
$$
\widetilde{F}_i=\bigoplus_{0\leq j\leq{i}/{2}}\widetilde{G}_{i-2j} \ar{\widetilde t} \bigoplus_{0\leq j\leq ({i-2})/{2}}\widetilde{G}_{i-2-2j}=\widetilde{F}_{i-2}\, .
$$
If $\sigma_{j}: \   \widetilde{G}_{i-2j}\to \widetilde{G}_{i-1}$ denotes the appropriate component of the 
map $\widetilde{\delta}: \   \widetilde{F}_{i}\to \widetilde{F}_{i-1}$, then $\s = \{\s_{j}\}$
is a system of higher homotopies on $\widetilde{\G}$, and $ \F \cong \std(\widetilde{\G},\s)$.
\end{mylist}
\end{prop}

\proof  
(2): Since the maps $\widetilde t$ are surjective, it follows inductively that
we may write $\widetilde {F_{i}}$ and $\widetilde t$ in the given form. 
The component corresponding to $ \widetilde{G}_{i-2j}\to \widetilde{G}_{i-1}$ in $\widetilde{\delta}: \   \widetilde{F}_{m}\to \widetilde{F}_{m-1}$ is the same  for any $m$ with $m\ge i-2j$ and $m\equiv i\  \mod (2)$ because
$\widetilde \delta$ commutes with $\widetilde t$.
 The condition 
that $\sigma$ is a sequence of higher homotopies is equivalent to the condition
that $\widetilde \delta^{2} = f\widetilde t$, as one sees by direct computation. It is now immediate
that $\F\cong \std( \widetilde{\G},\sigma)$.
\qed

\b{cor}\l{resolvingOverTwoRings} 
With hypotheses and notation as in Proposition~{\rm \ref{first}}, suppose in addition that
$S$ is a local ring and that $(\F,\delta)$ is a minimal $R$-free resolution of  $N$.
The minimal $S$-free resolution of $N$ is $(\widetilde{\G} ,  \widetilde{\partial})= \ker(\widetilde t)$.
If we split the epimorphisms $t: F_{i}\to F_{i-2}$ and correspondingly write
$F_{i} = \overline {G}_{i}\oplus F_{i-2}$ then the differential $\delta: \   F_{i}\to F_{i-1}$
has the form
$$
\delta_{i} = \bordermatrix{~ &  \overline {G}_{i}&F_{i-2}\cr 
 \overline {G}_{i-1}& \overline \partial_{i}&\varphi_{i}\cr  
F_{i-3}&  O&\delta_i}.
$$
\vglue -.5cm
\ \qed
\end{cor}

As an immediate consequence of Propositions~\ref{second}  and \ref{first} we obtain a result of Avramov-Gasharov-Peeva; their proof relies on the
spectral sequence proof of [AGP, Theorem~4.3].

\b{cor}\label{reduction}
 {\rm \cite[Proposition~6.2]{AGP}} Let  $f\in S$ be a non-zerodivisor in a local ring. If $N$ is a  module
over $S/(f)$ then the  CI operator $\chi$ corresponding to $f$ is a non-zerodivisor on 
$\ext_S(N,k)$ if and only if the minimal $S/(f)$-free resolution of $N$ is 
obtained by a Shamash construction applied to the minimal
free resolution of $N$ over $S$. \end{cor}

\proof Nakayama's Lemma shows that the CI operator $t: \ \F[-2]\to\F$ is surjective if and only if the operator
$\chi: \   \ext_{R}(N,k) \to \ext_{R}(N,k)$ is injective.\qed

\section{The minimal $R$-free resolution of a higher matrix factorization module} \label{Stoinf}

Let $(d,h)$ be a higher matrix factorization with respect to a regular sequence $\ff c$ in a ring $S$, 
 and $R= S/(\ff c)$. We will describe an $R$-free resolution of the HMF
module $M$ that is minimal when $S$ is local and $(d,h)$ is minimal.

\b{construction}\label{infresconstruction}\label{inres} 
Let $(d,h)$ be a higher matrix factorization with respect to a regular sequence $\ff c $ in a ring $ S$. Using notation as  in ~\ref{standardnotation}, choose splittings $A_{s}(p) = A_{s}(p-1)\oplus B_{s}(p)$ for $s=0,1$, 
so $$A_{s}(p) = \oplus_{1\le q\le p}  B_{s}(q)\,,$$
 and write $\psi_p$ for the component of $d_p$ mapping $B_{1}(p)$ to  $A_{0}(p-1)$. Set $$\A(p):\ A_1(p)\ar{d_p}A_0(p)\quad\hbox{ and }\quad \B(p):\ B_1(p)\ar{b_p}B_0(p)\, .$$
 \vglue -.5cm\ 
 \b{mylist}
\item[$\bullet$] Set $\U(1)=\B(1)$, and note that $h_1$ is a homotopy for  $f_1$.
Set $$\T(1):= \std(\U(1), h_1)\,.$$ Its beginning is the complex $R(1)\otimes \A(1)$. 
\item[$\bullet$] Given an $R(p-1)$-free resolution
$\T(p-1)$ of $M(p-1)$ with  beginning
$R(p-1)\otimes \A(p-1)$, 
let 
$$\Psi_p: \   R(p-1) \otimes \B(p)[-1] \to \T(p-1)$$ be the 
map of complexes induced by  $\psi_p: \   B_{1}(p)\to A_{0}(p-1)$. Set
$$
\U(p) := \MM\left(\Psi_p \right).
$$
We will show that $\U(p)$ is an $R(p-1)$-free resolution of $M(p)$. 
Thus we can choose
a system of higher homotopies $\sigma(p)$ for $f_{p}$ on $\U(p)$
that begins with $d_p$ (that is, $\sigma(p)_0=d_p$) and
$$
R(p-1) \otimes h_p: \    R(p-1) \otimes A_{0}(p)\to R(p-1) \otimes A_{1}(p).
$$
Set
$$
\T(p) := \std(\U(p), \s(p)).
$$
\end{mylist}

\myspace
\noindent  The underlying graded module of $\T(p)$ is $\U(p) = \MM(\Psi_p)$ tensored with a divided power algebra on a variable $y_{p}$ of degree $2$.
Its first differential is 
$$
R(p)\otimes \A(p): \   R(p)\otimes A_{1}(p)\ 
{\mathop{{\mathrel{\smash-}}{\mathrel{\mkern-8mu}}\a\a
{\mathrel{\smash-}}{\mathrel{\mkern-8mu}}\rightarrow}\limits^{R(p)\otimes d_{p}}}
\  R(p)\otimes A_{0}(p), 
$$
which is the presentation  of $M(p)$.
We see by induction on $p$ that the term $T_{j}(p)$ of homological degree $j$ in
$\T(p)$
is a direct sum of the form
\begin{equation}\label{decomposition}
 T_{j}(p) = \bigoplus y_{q_{1}}^{(a_{1})}\cdots y_{q_{i}}^{(a_{i})} B_{s}(q)\otimes R(p)
\end{equation}
where the sum is over all terms with
\b{align*}
p\geq q_{1} &> q_{2}>\cdots > q_{i}\geq q \geq 1, \\ 
a_{m}&>0\ \hbox{for}\ 1\le m\le i\,,\\ 
j &= s+\sum_{1\le m\le i}2a_{m}.
\end{align*}
We say that an element $y_{q_{1}}^{(a_{1})}\cdots y_{q_{i}}^{(a_{i})} v$ with $v\in B_{s}(q)$
and $a_{1}>0$ is
{\it admissible of weight $q_{1}$}, and we make the convention that the admissible elements in $B_{s}(q)$  have weight $0$.

The complex $\T(c)$ is thus filtered by:
$$\T(0):=0\subseteq
R\otimes \T(1)\subseteq\cdots \subseteq R\otimes \T(p-1)\subseteq
\T(c)\, ,$$
where $R\otimes \T(p)$ is the subcomplex spanned  by elements of weight $\le p$ with 
with $v\in B_{s}(q)$ for $q\le p$.
\end{construction}

\b{thm}\l{inresthm}
With notation and hypotheses as in  {\rm \ref{inres}}:
\b{mylist}
\item[\rm (1)] The complex $\T(p)$ is an $R(p)$-free resolution of $M(p)$ whose
first differential is $R(p)\otimes d_p$ and whose second differential is
$$
R(p)\otimes \left(\Big(\oplus_{q\leq p} \,A_{0}(q)\Big) \ar{h} A_{1}(p)\right),
$$
where the $q$-th component of $h$ is 
$h_q: \   A_{0}(q) \to A_{1}(q)\hookrightarrow A_{1}(p)$.

\item[\rm (2)] If $S$ is local then $\T(p)$ is the minimal free resolution of $M(p)$ if and only
if the higher matrix factorization $\Big(d_p,h(p)=(h_{1}|\cdots|h_{p})\Big)$  (see  {\rm \ref{standardnotation}} for notation) is  minimal.

\end{mylist}
\end{thm}

\noindent{\prooffont Proof of Theorem~\ref{inresthm}(1):}
We do induction on $p$. 
To start the induction, note that $\U(1)$ is the two-term complex
$\A(1) = \B(1)$. By hypothesis, its differential $d_1$  and homotopy
$h_1$ form a hypersurface matrix factorization for $f_{1}$, and 
$\T(1)$ has the form
$$
\T(1): \   R(1)\otimes \left( \cdots \ar{h_1} A_{1}(1)\ar{d_1} A_{0}(1)\ar{h_1} A_{1}(1)\ar{d_1} A_{0}(1)\right).
$$

Inductively, suppose that $p\geq 2$, and that 
$$
\T(p-1): \   \cdots \to T_{2}\to T_{1}\to T_{0}
$$ 
is an $R(p-1)$-free resolution of 
$M(p-1)$ whose first two maps are as claimed. 
We write $\overline -$ for $R(p-1) \otimes -$.
It follows that the first
map of $\U(p)$ is 
$$
\overline{d}(p): \overline{A}_{1}(p) = T_{1}\oplus \overline{B}_{1}(p)\to \overline{A}_{0}(p)= T_{0}\oplus  \overline{B}_{0}(p).
$$
 Since $R(p-1)\otimes \left(d_p h_p\right) = f_{p}\,\Id_{A_{0}(p)}$ we may take $R(p-1)\otimes h_p$ to
be the start of a system of higher homotopies $\sigma(p)$ for $f_{p}$ on $R(p-1)\otimes \U(p)$.
It follows from the definition that the first two maps in 
$\T(p)=\std(\U(p),\sigma(p))$ are as asserted.

By Proposition~\ref{second},
the Shamash construction takes an $R(p-1)$-free resolution to an $R(p)$-free
resolution of the same module. Thus  for the induction
it suffices to show that $\U(p)$ is an $R(p-1)$-free resolution of $M(p)$.
Since the first map of $\U(p)$ is $\overline{d}(p)$, and since $\overline{h}(p)$
is a homotopy for $f_{p}$, we see at once that
$$
H_{0}(\U(p)) = \coker( \overline{d}(p)) = \coker( R(p)\otimes d_p )= M(p).
$$

To prove that $\U(p)$ is a resolution, note first that $\U(p) _{\ge 2}= \T(p-1)_{\geq 2}$, 
and the image of $U(p)_{2} = T(p-1)_{2}$ is contained in the summand $T(p-1)_{1}\subseteq U(p)_{1}$,
so
$H_{i}(\U(p)) = H_{i}(\T(p-1) )= 0$ for $i\geq 2$. Thus it suffices to prove that $H_{1}(\U(p)) = 0$.

Let $(y,v) \in U(p)_{1} = T(p-1)_1 \oplus \overline{ B}(p)_1$ be a cycle in $\U(p)$. Thus, 
$\overline{b}_p(v)=0$ and $\overline{\psi}_p(v)=-\overline{d}_{p-1}(y)$. By Lemma~\ref{newex},
we conclude that $v=0$.
\qed

 \myspace
For the proof of part (2) of Theorem~\ref{inresthm} we will use the form of the resolutions $\T(p)$
to make a special lifting of the differentials to $S$, and thus to produce
especially ``nice'' CI operators. We pause in the proof of Theorem~\ref{inresthm} to describe this construction and deduce
some consequences.

\b{prop}\l{speciallifting1}
With notation and hypotheses as in  {\rm \ref{inres}}, there exists a  lifting
of the filtration $\T(1)\subseteq \cdots \subseteq \T(c)$  to a filtration
$
\tilde\T(1) \subseteq \cdots \subseteq \tilde\T(c)
$ 
over $S$, and a lifting
$\widetilde{\delta}$ of the differential $\delta$ in $\T(c)$ to $S$
with lifted CI operators $\tilde t_1,\dots ,\tilde t_c$  on  $\tilde \T(c)$ such that
for every $1\le p\le c$:
\b{mylist}
\item[\rm (1)]  Both $\widetilde{\delta}$  and  $\widetilde{t}_{p} $ 
 preserve $\tilde \T(p)$, and $\widetilde{t}_{p} \Big\vert_{\widetilde \T(p)}$
 commutes with  $\widetilde{\delta}\Big\vert_{\tilde \T(p)}$
on $\widetilde \T(p)$ .
 \item[\rm (2)] The CI operator ${t}_p$ vanishes on the subcomplex $R\otimes \U(p)$
and induces an isomorphism from $R\otimes T(p)_j/U(p)_j$ to $R\otimes T(p)_{j-2}$ that sends
an admissible element $y_{q_{1}}^{(a_{1})}\cdots y_{q_{i}}^{(a_{i})} v$ with $q_1=p$ to $y_{q_{1}}^{(a_{1}-1)}\cdots y_{q_{i}}^{(a_{i})} v$.
\end{mylist}
\end{prop}

\proof
If $p=1$ the result is obvious. Thus we may assume by induction that liftings 
$$
0\subset \tilde \T(1) \subseteq \cdots \subseteq \tilde\T(p-1), 
$$ 
$\tilde \delta(p-1)$ and $\tilde t_1,\dots ,\tilde t_{p-1}$ on $\tilde \T(p-1)$ satisfying the Proposition have been constructed.
We use the maps $\psi_{p}$ and $b_{p}$ from the definition of the higher matrix factorization to construct a lifting of $\U(p)$ from 
the given lifting of $\T(p-1)$.
In addition, we choose liftings $\widetilde\sigma$ of the maps (other than the differential) in the
system of higher homotopies $\sigma(p)$ for $f_p$ on $\U(p)$.

By construction, $\T(p) = \std(\U(p),\sigma(p))$, so we take the standard lifting to $S$ from \ref{standardCI}, that is,   take $\tilde\T(p)=\oplus_{i\geq 0} \,y_{p}^{(i)}\tilde\U(p)$ 
with
lifting of the differential
$\widetilde \delta = \sum
\,t^{\,j}\otimes  \widetilde{\sigma}_j
$, where $t$ is the dual variable to $y_p$.

By  Construction~\ref{standardCI} it follows that,  modulo $(f_1,\dots ,f_{p-1})$, the map $\widetilde {\delta }^2$ vanishes on $\widetilde{\U}(p)$ and induces
$f_p$ times the projection $  \widetilde{T}_j(p)/ \widetilde{U}_j(p)\to   \widetilde{T}_{j-2}(p)$.

We choose $\widetilde{t}_p$ to be the standard lifted CI operator, which
 vanishes on $\widetilde{\U}(p)$ and is the projection $  \widetilde{T}_j(p)/ \widetilde{U}_j(p)\to   \widetilde{T}_{j-2}(p)$. Then  $\widetilde{\delta}_{i-2}\widetilde{t}_p=
\widetilde{t}_p\widetilde{\delta}_i$ by construction; see \ref{standardCI}.

Recall that $\tilde {\delta}\Big\vert_{\widetilde{\T}(p-1)}$ is the lifting $\tilde \delta(p-1)$ given by induction.
Therefore,  from $\tilde \delta$ we can choose maps $\tilde t_1,\dots ,\tilde t_{p-1}$ on $\widetilde{\T}(p)$
that extend the
maps  $\tilde t_1,\dots ,\tilde t_{p-1}$ given by induction on $\widetilde{\T}(p-1)\subseteq \widetilde{\U}(p)$. 
\qed

\myspace
The  CI operators commute up to homotopy, and
 it is an open conjecture from  \cite{Ei1} (see also \cite[Section~9]{AGP}) that they
 can be chosen to commute when restricted to the minimal free resolution of
 a high syzygy in the local case. 
Proposition~\ref{speciallifting1} allows us to give a partial answer, based on the following general criterion.

\b{prop}\label{comcriterion}
 Let $\ff c$ be a regular sequence in a local ring $S$, and 
let $R = S/(\ff c)$. Suppose that $(\F,\delta)$ is a complex over $R$ with lifting
$(\tilde \F, \tilde \delta)$ to $S$, and let $\tilde t_{1}, \dots, \tilde t_{c}$  on $\tilde\F$  define
CI operators corresponding to $\ff c$.
 If, for some $j$, $\widetilde{t_j}$ commutes with $\tilde \delta^{2}$, then $t_{j}$ commutes with each $t_{i}$.
\end{prop}

\proof
Since $\tilde \delta^{2} = \sum f_{i}\tilde t_{i}$ by definition, we have
$\sum f_{i}\tilde t_{j} \tilde t_{i} = \sum f_{i}\tilde t_{i}\tilde t_{j}$, or equivalently
$\sum f_{i}(\tilde t_{j} \tilde t_{i}-\tilde t_{i}\tilde t_{j}) = 0$. Since $\ff c$ is a regular
sequence it follows that $\tilde t_{j} \tilde t_{i}-\tilde t_{i}\tilde t_{j}$ is zero modulo
$(\ff c)$ for each $i$.\qed

\myspace
As an immediate consequence, we have:

\b{cor}\label{partialcommutativity} 
Suppose that $S$ is local. With CI operators on $\T(p)$ chosen as in 
Proposition~ {\rm \ref{speciallifting1}} the operator $t_{p}$ commutes on $\T(p)$ with
each $t_{i}$ for $i<p$. \qed
\end{cor}

\begin{cor}\label{action on ext} Let $k[\chi_{1},\dots, \chi_{c}]$ act on $\ext_{R}(M,k)$  as in
Construction~{\rm\ref{operators on Ext}}. There is an isomorphism
 $$
\ext_{R}(M,k) \cong \bigoplus_{p=1}^{c} k[ \chi_{p}, \dots, \chi_{c}] \otimes_{k} \hom_S(\B(p),k)
$$
of vector spaces such that, for $i\geq p$,  $\chi_{i}$ preserves the summand 
$$
k[ \chi_{p}, \dots, \chi_{c}] \otimes \hom_S(\B(p),k)
$$
 and acts on it via the action on the first factor.
\end{cor}

\proof Since $\T(c)$ is a minimal free resolution of $M$, the $k[\chi_{1},\dots, \chi_{c}]$-module $\ext_{R}(M,k)$
is isomorphic to
$
\Hom_{R}(\T(c),k).
$
Using the decomposition in (\ref{decomposition}) we see that the underlying graded free module of $\Hom_{R}(\T(c),k)$ is
$$\bigoplus_{p}k[\chi_{p},\dots,\chi_{c}]\otimes_{k}\Hom_{S}(\B(p),k).$$ From part (2) of Proposition~\ref{speciallifting1} we see that,
for $i\geq p$, the action of
$\chi_{i}$ on the summand $k[\chi_{p},\dots,\chi_{c}]\otimes_{k} \Hom_{S}(\B(p),k)$ is via the natural action on the first factor.
\rightline{\qed}

\myspace
Corollary~\ref{action on ext} provides a standard decomposition of $
\ext_{R}(M,k)$ in the sense of \cite{EP}.

We will  complete the proof of Theorem~\ref{inresthm}:

 \myspace
\noindent{\prooffont Proof of Theorem~\ref{inresthm}(2):} We suppose that $S$ is local with maximal ideal $\mm$. If the resolution $\T(p)$ is minimal then it follows at once from the description of the first
two maps that $(d,h)$ is minimal. We  will prove the converse
by induction on $p$.

If $p=1$ then $\T(1)$ is the periodic resolution
$$
\T(1): \ \ \ \ \cdots\ar{h_{1}} A_{1} \ar{d_{1}}A_{0}\ar{h_{1}} A_{1}\ar{d_{1}} A_{0}
$$
and only involves the maps $(d_1, h_1)$; this
 is obviously minimal if and only if $d_1$ and $h_1$ are minimal.

Now suppose that $p>1$ and that $\T(q)$ is minimal for $q<p$. 
Let $\delta_{i}: \   T_{i}(p)\to T_{i-1}(p)$ be the differential of $\T(p)$. 
We will prove minimality of $\delta_{i}$ by a second induction, on $i$, starting with $i=1,2$.

Recall that the underlying graded module of $\T(p) = \std(\U(p), \s)$
is the divided power algebra $S\{y_{p}\} = \sum_{i} Sy_{p}^{(i)}$
tensored with the underlying module of $R(p)\otimes \U(p)$. Thus the beginning
of the resolution $\T(p)$ has the form
$$
\cdots \to R(p)\otimes y_{p}A_{0}(p) \oplus R(p) \otimes T_{2}(p-1) \ar{\delta_{2}}
R(p) \otimes A_{1}(p) \ar{\delta_{1}} R(p) \otimes A_{0}(p).
$$
The map
$\delta_{1}$ is induced by $d_p$, which is minimal by hypothesis.
Further, $\delta_{2} = (h_p, \partial_2)$
where the map $\partial_2$ is the differential of $\T(p-1)$ tensored with $R(p)$.
The map $h_p$ is minimal by hypothesis, 
and $\partial$ is minimal by induction on $p$, so $\delta_{2}$ is minimal as well.

Now suppose that $j\ge 2$ and that $\delta_{i}$ is minimal for $i\leq j$. 
We must show that $\delta_{j+1}$ is minimal, that is, $\delta_{j+1} (w)\in \mm T_{j}(p)$
for any $w\in T_{j+1}(p)$. 
By Construction~\ref{inres}, 
$\delta_{j+1} (w)$ can be written uniquely as a  sum of admissible elements of the form
$$
y_{q_{1}}^{(a_{1})}\cdots y_{q_{i}}^{(a_{i})} v
$$
with $0\neq v\in B_{s}(q)$ and
\b{align*}
p\geq q_{1} &> q_{2}>\cdots > q_{i}\geq q \geq 1, \\ 
a_{m}&>0 \ \ \hbox{for}\ 1\le m\le i,\\ 
j &= s+\sum_{1\le m\le i}2a_{m}\, .
\end{align*}

If $\delta_{j+1}(w)\notin \mm T_{j}(p)$ then there exists a  summand
$
y_{q_{1}}^{(a_{1})}\cdots y_{q_{i}}^{(a_{i})} v
$
 in this expression that is not in $\mm T_{j}(p)$.  
 Since $\delta_{j+1} (w)$ has homological degree $j\geq 2$, the weight of this
 summand must be $>0$, that is, a factor $y_{q_{1}}^{(a_{1})}$ must be present.

Choose such a summand  with weight $q'_{1}$  as large as possible.
We choose $t_{q'_{1}}$ as in Proposition~\ref{speciallifting1}. Then
 $t_{q'_{1}}$ sends every admissible element of weight $< q'_{1}$ to zero. The admissible summands of $\delta_{j+1}(w)$ with weight $>q'_{1}$ can be ignored
 since they are in $\mm T_{j-2}(p)$. By Proposition~\ref{speciallifting1} it follows that
 $t_{q'_{1}}\delta_{j+1}(w)\notin\mm T_{j-2}(p)$. Since
$$
t_{q'_{1}}\delta_{j+1}(w) = \delta_{j-1}t_{q'_{1}}(w)\, ,
$$ 
this contradicts the induction hypothesis.
\qed

\myspace
Gulliksen \cite{Gu} shows  that the  Poincar\'e series of $M$  over $R$  has the form
${\poin}^R_{M}(x)=g(x)(1-x^2)^{-c}$ for some $g(x)\in{\bf Z}[x]$, and his finite generation result
implies
 that the Betti numbers  are eventually given by two 
polynomials of the  same degree.  Avramov  \cite[Theorem~4.1]{Av} showed that they have  the same leading coefficient.
We can make this very explicit.

\b{cor}\label{bettiinf}
With notation and hypotheses as in  {\rm \ref{inres}}, if in addition $S$ is local and the higher matrix factorization $(d,h)$ is minimal, then:
\b{mylist}
\item[\rm (1)]
The  Poincar\'e series of $M$  over $R$  is
$${\poin}^R_{M}(x)=\sum_{1\le p\le c}\ {1\over{(1-x^2)^{c-p+1}}} \Big(x\, \rank(B_1(p))+\rank(B_0(p))\Big)\, .$$ 
\item[\rm (2)]
 The Betti numbers of $M$ over $R$ are given by the following two polynomials in $z$:
\b{align*}
 b_{2z}^R(M)&=\sum_{1\le p\le c}\ {c-p+z\choose c-p}\rank(B_0(p))\, \\ 
  b_{2z+1}^R(M)&=\sum_{1\le p\le c}\ {c-p+z\choose c-p}\rank(B_1(p))\, .
 \,
 \end{align*}
 \end{mylist}
 \end{cor}

\proof
 For (2), recall that
 the Hilbert function of $k[Z_{p},\dots ,Z_c]$ is $ g_p(z)={c-p+z\choose c-p+1}$. 
 \rightline{\qed}

\myspace
Recall that the  {\it complexity} of an $R$-module $N$ is defined to be 
$$
\hbox{cx}_R(N)=\inf\{ q\ge 0\,|\, \hbox{there exists a } w\in {\bf R} \hbox{ such that } b_i^R(N)\le wi^{q-1} \hbox{ for } i\gg 0\,\}\,.
$$
If the complexity of $N$ is $\mu$ then, as noted above, $$\dim_{k}\ext^{2i}_{R}(N,k) = (\beta/(\mu-1)!)i^{\mu-1}+O(i^{\mu-2})$$ for $i\gg 0$. Following \cite[7.3]{AB} $\beta$ is called the {\it Betti degree} of $N$ and denoted $\hbox{Bdeg}(N)$; this is the multiplicity of the module
$\ext_R^{even}(N,L)$, which is equal to the multiplicity of the module
$\ext_R^{odd}(N,L)$.

\b{cor}\label{complexity} 
With notation and hypotheses as in  {\rm \ref{inres}}, suppose  in addition that $S$ is local. Suppose that  $(d,h)$ is a minimal higher matrix factorization, and set
$$
\gamma=\min\{\,p\,|\,B_1(p)\not=0\,\}\, .
$$
The complexity of $M
:=M(c)$ is  
$$
{\rm cx}_R\,M=c-\gamma+1.
$$
Moreover,  $B_0(p)=0$ for $p<\gamma$, and the Betti degree of $M$ is 
$$\hbox{\rm Bdeg}(M)=\rank(B_1(\gamma))=\rank(B_0(\gamma))\, .$$ If in addition $S$ is Cohen-Macaulay, then
 $\rank(B_1(p))>0$ for  every $\gamma\le p\le c$.
 \end{cor}

\proof
By {Corollary}~\ref{filteredFiniteResolution},
$B_1(p)=0$ implies that $B_0(p)=0$.
Hence the Betti degree of $N$ is equal to $\min\{\,p\,|\,B_1(p)\not=0\,\}$ and $B_0(p)=0$ for $p<\gamma$.

The equality
$\rank(B_1(\gamma))=\rank(B_0(\gamma))$ follows 
since $M(\gamma)$  is annihilated by $f_\gamma$ and has minimal free resolution
$B_1(\gamma)\ar{b_\gamma}B_0(\gamma)$ over  
$S/(f_1,\dots ,f_{\gamma-1})$.

Corollary~\ref{nozeros} implies that $\rank(B_1(p))>0$ for  every $\gamma\le p\le c$, when $S$ is Cohen-Macaulay.
\qed

\section{Resolutions over intermediate rings} \label{Sintemediate}

Using a slight extension of the definition of a higher matrix factorization we can describe the resolutions of the modules
$M(p)$ over any of the rings $R(q)$ with $q< p$.

\b{definition}
 A {\it generalized matrix factorization} over a ring $S$ with respect to a regular
sequence $\ff c\in S$ is a pair of maps $(d,h)$ satisfying the definition
of a higher matrix factorization {\it except} that we drop the assumption that  $A(0)=0$, so that we have a map of free modules $A_1(0)\ar{b_0} A_0(1)$. We do {\it not} require the existence of a map $h_0$.
\end{definition}

\b{construction}\label{interres} 
Let $(d,h)$ be a generalized matrix factorization with respect to a regular sequence $\ff c $ in a ring $ S$. Using notation as  in ~\ref{standardnotation},   we 
choose splittings $A_{s}(p) = A_{s}(p-1)\oplus B_{s}(p)$ for $s=0,1$, and write $\psi_p$ for the component of $d_p$ mapping $B_{1}(p)$ to  $A_{0}(p-1)$. 
\b{mylist}
\item[$\bullet$]
Let $\bf V$ be a free resolution of 
the module $\coker (b_0)$ over $S$, and set $\Q(0):={\bf V}$.
\item[$\bullet$] 
Let $$\Psi_1: \   \B(1)[-1]\to \Q(0)$$
be the map of complexes
induced by $\psi_1: \   B_{1}(1)\to A_{0}(0)$, and  set
$$\Q(1) = \MM(\Psi_1).$$
\item[$\bullet$] For $p\geq 2$, suppose that an $S$-free resolution
$\Q(p-1)$ of $M(p-1)$ with first term 
$Q_{0}(p-1) = A_{0}(p-1)$ has been constructed.
Let $$\psi_p': \   \B(p)[-1]\to \L(p-1)$$
be the map of complexes
induced by $\psi_p: \   B_{1}(p)\to A_{0}(p-1)$, and  let
$$
\Psi_p: \   \KK(f_{1},\dots, f_{p-1})\otimes \B(p) [-1]\to \Q(p-1)
$$ be an $(f_1,\dots ,f_{p-1})$-Koszul extension. Set
$\Q(p) = \MM(\Psi_p).$
\end{mylist}
\end{construction}

\noindent
The proof of Theorem~\ref{fresthm} can be applied in this situation and yields the following result.

\b{prop}\label{intermediatelemma}
Let $(d,h)$ be a generalized matrix factorization over  a ring $S$, and
let $\bf V$ be a free resolution of 
the module $\coker (b_0)$ over $S$. 
For each $p$, the complex
$\Q{(p)}$, constructed in  {\rm \ref{interres}},  is an $S$-free resolution of
the module $M(p)$. 
If the ring $S$ is local then the resulting free resolution is minimal if and only if $(d,h)$ and $\bf V$ are minimal. \qed
\end{prop}

\b{thm}\label{intermediate}
Let $(d,h)$ be a  higher matrix factorization. Fix  a number $1\le j\le c-1$. Let $\T(j)$ be the free resolution of $M(j)$ over the ring $R(j)=S/(f_1,\dots ,f_j)$ given by Theorem~ {\rm \ref{inres}}.
Let $(d',h')$ be the generalized matrix factorization 
 over the ring $R(j)$ with  
\b{align*}
A_s(0)&=R(j) \otimes \Big(\oplus_{1\le q\le j}A_s(q) \Big)\quad \hbox{and}\quad d_0'=R(j)\otimes d_j,\\  
\hbox{for} \ p>j,\ \ 
A_s(p)'&=R(j) \otimes A_s(p+j)\quad \hbox{and}\quad d_{p}'=R(j)\otimes d_{p+j}\, , 
 \end{align*}
 for $s=0,1$ and maps induces by $(d,h)$. Then $M'(0)=M(j)$.
\myspace

\b{mylist}
\item[\rm (1)] Construction~ {\rm \ref{interres}},  starting from
the $R(j)$  free resolution $\Q(0):=\T (j)$ of $M'(0)=M(j)$, 
produces a free resolution $\Q(c-j)$ of $M$ over $R(j)$. 
\vglue .1cm
\item[\rm (2)] If $S$ is local and
$(d,h)$ is minimal, then the resolution $\Q(c-j)$ is minimal. In that case,
the Poincar\'e series of $M$  over $R(j)$  is
\b{align*}{\poin}^{R(j)}_{M}(x)=
&\Bigg(\sum_{1\le p\le j}\ {1\over{(1-x^2)^{p-j-1}}} \Big(x\, \rank(B_1(p))+\rank(B_0(p))\Big)\Bigg)\\  &\Bigg(\sum_{j+1\le p\le c}\ (1+x)^{p-j-1} \Big(x\,\rank (B_1(p))+\rank(B_0(p))\Big)\Bigg)
\, .
\end{align*}
\end{mylist}
\end{thm}

\proof
First, we apply Theorem~\ref{inres}, which gives the resolution $\T (j)$
of $M(j)$ over the ring $R(j)$. Then we apply Proposition~\ref{intermediatelemma}.
\qed

\section{Pre-stable Syzygies and Generic CI Operators} \label{Sgeneric}

Our goal in this section and Section~\ref{SmakeMF} is to show that every sufficiently high syzygy over a complete intersection is  an HMF module. In this section we introduce the concepts  of  {\it pre-stable syzygy} and {\it stable syzygy}. We 
will see that any sufficiently high syzygy in a minimal free resolution over a local complete intersection ring  is a stable syzygy. 
In Section~\ref{SmakeMF} we will
show that a pre-stable syzygy  is an  module.

\b{definition}\label{wellb}
Suppose that $f_{1}, \dots, f_{c}$ is a regular sequence in a
 local ring $S$, and set
 $R= S/(f_{1}, \dots, f_{c})$. We define the concept of a pre-stable syzygy recursively:
 We say that an $R$-module $M$ is a {\it pre-stable syzygy}  with respect to $\ff c$ if 
 either $c=0$ and $M=0$, or $c\ge 1$ and the following conditions are satisfied:
 \vglue .1cm
 \b{mylist}
 \item[(1)] There exists a minimal $R$-free resolution $(\F,\delta)$ of an $R$-module of finite projective dimension over $S$ with a surjective
 CI operator $t_{c}$ on $\F$ and such that $M=\ker(\delta_1)$;
 \vglue .1cm
 \item[(2)] If $\tilde \delta_{1}$ is a lifting of $\delta_{1}$ to $\tilde R := S/(\ff{c-1})$, then  $\tilde M := \ker(\tilde \delta_{1})$ is a pre-stable syzygy with respect to $\ff{c-1}$.
 \end{mylist}
 \noindent We say that a pre-stable syzygy is {\it stable} if the module resolved by $\F$
 in Condition (1) in \ref{wellb} is maximal Cohen-Macaulay and the module $\tilde M$ in Condition (2)
 is a stable syzygy.

 \end{definition}

\b{remark}\label{stableremark}
 The property of being pre-stable is independent of choices: Condition (1) of the definition is independent of the choice of $t_c$ because $t_c$ is uniquely defined up to homotopy, and $\F$ is assumed minimal. Condition (2) is independent of the choice of the lifting of $\delta_{1}$ because, if we write $L$ for the module
resolved by $\F$, then $\ker(\tilde \delta_{1})$  is the second syzygy of $L$ over $\tilde R$
by Propositions~\ref{second} and \ref{first}. 

Note that if $M$ is a pre-stable syzygy, then by (1) it follows that $M$ has finite projective dimension over $S$.
\end{remark}

  The property described in Definition~\ref{wellb} is preserved under taking syzygies:

\b{prop}\label{higherSyzygies}
Suppose that $f_{1}, \dots, f_{c}$ is a regular sequence in a
 local ring $S$, and set
 $R= S/(f_{1}, \dots, f_{c})$. If $M$ is  a pre-stable syzygy over $R$, then 
$\syz_{1}^{R}(M)$ is pre-stable as well.
If $M$ is a stable syzygy over $R$, then so is $\syz_{1}^{R}(M)$.
\end{prop}

\proof
Let $(\F, \delta)$ be a minimal  $R$-free resolution of a module $L$ such that $M = \ker ( \delta_{1})$ and the conditions  in Definition~\ref{wellb} are satisfied. 
 Lifting $\F$ to $\tilde \F$ over $\tilde R:= S/(\ff{c-1})$ and using the hypothesis that $S$ is local, we see that
 the lifted CI operator $\tilde t_{c}$  is surjective on $ \tilde F$.
By Propositions~\ref{second} and \ref{first}, $\tilde \G := \ker (\tilde t_{c})$ is the minimal free resolution of the module $L$ over $\widetilde R$.

 Let $M'=\syz_{1}^R(M)$ and let $L'=\syz_{1}^R(L)$, so that
 $\F'=\F_{\ge 1}[-1]$ is the minimal  free resolution of $L'$. Clearly $t_{c}\Big\vert_{ \F'}$ is surjective.
  The shifted truncation $\tilde \F':= \tilde \F_{\geq 1}[-1]$ is a lifting of $\F'$, and  
$\tilde \G' := \ker \Big(\tilde{t}_c\Big\vert_{\tilde {\F}'}\Big)$
 is a minimal free resolution of $L'$ over $\tilde R$. The complex $
\tilde \G'_{\ge 2} $
agrees (up to the sign of the differential) with $\tilde \G[-1]_{\geq 2}$: 
\myalign\label{syzygies}
 \tilde \G:\quad\quad \dots\to \tilde G_4\to\tilde G_3\to\tilde G_2\to&\tilde F_1\ar{\delta_1}\tilde F_0\\ 
 \nonumber
  \tilde \G':\quad\quad \dots\to \tilde G_4\to\tilde G_3\to\tilde F_2\ar{\delta_2}&\tilde F_1\, ,
  \end{align}
Thus
  $\ker(\widetilde{\delta}_{2})=
 \syz_{1}^{\widetilde R}(\ker (\widetilde{\delta}_{1})).$
  Since $\ker (\widetilde{\delta}_{1})$ is a pre-stable syzygy, we can apply the induction hypothesis to conclude that $\ker(\widetilde{\delta}_{2})$ is pre-stable.
  
  The last statement in the proposition follows from the observation that if $L$ is a maximal Cohen-Macaulay $R$-module, then so is $L'$.
  \qed

\myspace
The next result shows that in the codimension $1$ case, pre-stable syzygies are the same as codimension $1$ matrix factorizations.

\b{prop}\label{stableCodimOne} 
Let $f\in S$ be a non-zerodivisor in a local ring and set
$R=S/(f)$. The following conditions on an $R$-module $M$ are equivalent:
\b{mylist}
\item[\rm (1)] $M$ is a pre-stable syzygy with respect to $f$.
\item[\rm (2)] $M$ has projective dimension $1$ as an $S$-module.
\item[\rm (3)] The minimal $R$-free resolution of $M$ comes from a codimension $1$ matrix factorization of $f$ over $S$.
\end{mylist}
\end{prop}

\proof $(1)\Rightarrow(2)$: Let $\F$ be a minimal free resolution 
satisfying condition (1) in Definition~\ref{wellb}.  
By Proposition~\ref{higherSyzygies} and its proof and notation,
$\syz_2^R(M)$ is a pre-stable syzygy, and thus the free resolution $$
\tilde \G':\ \ \dots \to \widetilde{G}_4\to \widetilde{F}_3\to\widetilde{F}_2$$
(which is the kernel of the lifting of the CI operator $t_c$ on  the minimal free resolution $\F_{\ge 2}$
of $M$)
is zero in degrees $\geq 4$. Since $\tilde \G'$ is the minimal free resolution (up to a shift) of 
 $M$ over $S$, the projective dimension of $M$ over $S$ is $1$.
\vglue .1cm
\hglue -.3cm $(2)\Rightarrow(3)$: If $M$ has projective dimension $1$ then $M$ is the 
cokernel of a square matrix over $S$, and the homotopy for multiplication
by $f$ defines the matrix factorization.
\vglue .1cm
\hglue -.3cm  $(3)\Rightarrow(1)$: Continuing the periodic free resolution of $M$
as an $R$ module two steps to the right we get a minimal free resolution $\F$ of a module $L \cong M$ on which the CI operators are surjective, and also injective on $\F_{\ge 2}$.
It follows that $\ker (\widetilde{\delta}_1)=0$ in the notation of Definition~\ref{wellb}, so it is pre-stable.
\qed

\myspace
We now return to the situation of Theorem~\ref{extGeneration}:
Let $N$ be an $R$-module with finite projective dimension over $S$.
We  regard $E:=\ext_{R}(N,k)$ as a module over $\RR = k[\chi_{1}, \dots, \chi_{c}]$, where $\chi_j$ have degree $2$. Since we think of degrees
in $E$ as cohomological degrees, we write $E[a]$ for the shifted module
whose degree $i$ component is $E^{i+a} = \ext_{R}^{i+a}(N,k)$. 
If $M$ is the $r$-th syzygy module of $N$ then
$\ext_{R}(M,k) = \ext_{R}^{\geq r}(N,k)[-r]$.

Recall that the {\it Castelnuovo-Mumford regularity} 
$\reg\, E$
is defined as
$$
\reg\, E = \max_{0\leq i\leq c}\Big\{i + \{\max\{j \mid H_{(\chi_{1},\dots, \chi_{c})}^{i}(E)^{j} \neq 0\}\} \Big\}.
$$
Since the generators of $\RR$ have degree
$2$, some care is necessary. 
Note that if $\ext^{odd}_{R}(N,k)\allowbreak\neq 0$ then $E= \ext_{R}(N,k)$ can never have regularity $\leq 0$, since it is generated in
degrees $\geq 0$ and the odd part cannot be generated by the even part. Thus we will often have recourse to the condition $\reg\, \ext_{R}(N,k)=1$.
On the other hand, many things work as usual. If we split $E$ into
even and odd parts, $E = E^{even}\oplus E^{odd}$ we have
 $\reg\, E = \max(\reg\, E^{even}, \reg\, E^{odd})$ as usual. Also, 
if $\chi_{c}$ is a non-zerodivisor on $E$ then 
$\reg(E/\chi_{c}E) = \reg \,E$.

\b{thm}\label{epitwo} \label{highsyzprop} 
Suppose that $f_{1}, \dots, f_{c}$ is a regular sequence in a
 local ring $S$ with infinite residue field $k$, and set
 $R= S/(f_{1}, \dots, f_{c})$.  Let $N$ be an $R$-module with finite projective dimension over $S$,
 and let $\L$ be the minimal $R$-free resolution of $N$.
 There exists a non-empty Zariski open dense set $\mathcal Z$ of upper-triangular matrices $(\alpha_{i,j})$
with entries in $k$, such that 
for every
$$
r\ge 2c-1+\reg( \ext_{R}(N,k))\,
$$
the syzygy module $\syz_r^R(N)$ is  pre-stable with respect to
the regular sequence $f'_{1}, \dots,\allowbreak  f'_{c}$ with $f'_{i} = f_{i}+\sum_{j>i}\alpha_{i,j}f_j$.
\end{thm}

 To prepare for the proof of Theorem~\ref{highsyzprop} we will explain the property of the regular
sequence $f'_{1}, \dots, f'_{c}$ that we will use.
Recall that a sequence of elements $\chi'_{c}, \chi'_{c-1},\dots,\chi_{1}'\in \RR$
is said to be an {\it almost regular sequence} on a graded module $E$ if, for $q = c,\dots,1$,
the submodule of elements of $E/(\chi'_{q+1},\dots,\chi'_{c})E$ annihilated by
$\chi'_{q}$ is of finite length. 

\myspace
We will use the following lemma with
 $E=\ext_{R}(N,k)$. 

\b{lemma}\label{genericMonotwo}\label{negreg}
Suppose that $E$ is a  non-zero graded module  of regularity $\le 1$ over ${\mathcal R}=k[\chi_1,\dots ,\chi_c]$. The element $\chi_{c}$ is  almost regular on $E$ if and only if  $\chi_{c}$ is a non-zerodivisor on $E^{\ge 2}[2]$ (equivalently, $\chi_{c}$ is a non-zerodivisor on $E^{\ge 2}$). 
 
More generally, if we set $E(c)=E$ and $$E(j-1)=E(j)^{\ge 2}[2]/\chi_{j}E(j)^{\ge 2}[2]$$ for $j\le c$, then
the sequence $\chi_{c}, \dots, \chi_{1}$ is almost regular on $E$ if and only if  
$\chi_{j}$ is a non-zerodivisor on 
$E(j)^{\ge 2}[2]$ for every $j$.
In that case $\reg\,E(i)\leq 1$.
\end{lemma}

\proof 
By definition the element $\chi_{c}$ is almost regular on $E$ if the submodule $P$
of $E$ of elements annihilated by $\chi_{c}$ has finite length. Since 
$\reg(E)\leq 1$, all such elements must be contained in $E^{\le 1}$.
Hence, $\chi_{c}$ is a non-zerodivisor on $E^{\geq 2}$.

Conversely, if $\chi_{c}$ is a non-zerodivisor on $E^{\geq 2}$ then $P\subseteq E^{\le 1}$
so $P$ has finite length. Therefore, $\chi_{c}$ is almost regular on $E$.

Thus $\chi_{c}$ is almost regular if and only if it is a non-zerodivisor on $E^{\geq 2}$ as claimed.

If $\chi_{c}$ is a non-zerodivisor on $E^{\geq 2}$, then
$$
\reg(E^{\geq 2}/\chi_{c}E^{\geq 2}) = \reg(E^{\geq 2})\leq 3,
$$
whence $\reg(E(c-1))\leq 1$. By induction, $\chi_{c-1}, \dots, \chi_{1}$
is an almost regular sequence on $E(c-1)$ if and only if 
$\chi_{j}$
is a non-zerodivisor on 
$E(j)^{\ge 2}[2]$ for every $j<c$, as claimed.
\qed

\myspace
The following result is a well-known consequence of the ``Prime Avoidance Lemma'' (see for example \cite[Lemma 3.3]{Ei3} for Prime Avoidance):

\b{lemma}\label{primeavoidance}
If $k$ is an infinite field and $E$ is a  graded module over the polynomial ring
$\RR = k[\chi_{1}, \dots, \chi_{c}]$, then 
there exists a non-empty Zariski open dense set $\mathcal Y$ of lower-triangular matrices $(\nu_{i,j})$
with entries in $k$, such that the sequence of elements
$\chi'_{c}, \dots, \chi'_{1}$ with $\chi'_{i} = \chi_{i}+\sum_{j<i}\nu_{i,j}\chi_j$
is almost regular on $E$. \qed
\end{lemma}

Again let  $f_{1}, \dots, f_{c}$ be a regular sequence in a
 local ring $S$ with infinite residue field $k$ and maximal ideal $\mm$, and set
 $R= S/(f_{1}, \dots, f_{c})$.  Let $N$ be an $R$-module with finite projective dimension over $S$,
 and let $\L$ be the minimal $R$-free resolution of $N$.
 Suppose we have CI operators defined by a  lifting $\widetilde{\L}$.
If we make a
change of generators  of $(\ff c)$ using
an invertible matrix $\alpha$ and  $f_{i} '= \sum_{j} \alpha _{i,j}f_{j}$ with $ \alpha _{i,j}\in S$, then
the lifted CI operators on the lifting $\widetilde{\L}$ change as follows:
$$
\widetilde \partial^{2} = \sum_{i}f_{i}'\widetilde t_{i}'
=\sum_{i}\Bigg(\sum_{j} \alpha _{i,j}f_{j}\Bigg)\widetilde t_{i}' = \sum_{j}f_j\Bigg(\sum_{i} \alpha _{i,j}
\widetilde t_{i}'\Bigg)\, .
$$
So the CI operators corresponding to the sequence $f_{1}, \dots, f_{c}$  are expressed as
$t_{j} = \sum_{i}\alpha _{i,j}t_{i}'$. Thus, if we make a
change of generators  of the ideal $(f_{1}, \dots, f_{c})$ using
a matrix $\alpha$ then
the CI operators transform by
the inverse of the transpose of $\alpha$.  Another way to see this is  from the fact that
${\mathcal R}=k[\chi_1,\dots ,\chi_c]$ can be identified with the 
 symmetric algebra of the dual of the vector space $(\ff c)/\mm(\ff c)$. 

In view of this observation, Lemmas~\ref{genericMonotwo} and 
\ref{primeavoidance} can be translated  as follows:

\b{prop}\label{whygeneric}
Let $\ff c\in S$ be a regular sequence in a local ring with infinite residue field $k$, and set $R:=S/(\ff c)$. Let $N$ be an $R$-module of finite projective dimension over $S$, and set $E:= \ext_R(N,k)$.
\b{mylist}
\item[\rm (1)] {\rm \cite{Av, Ei1}} 
There exists a non-empty Zariski open dense set $\mathcal Z$ of upper-triangular matrices $\overline\alpha=(\overline\alpha_{i,j})$
with entries in $k$, such that if $\alpha=(\alpha_{i,j})$ is any matrix over $S$ that
reduces to $\overline\alpha$ modulo the maximal ideal of $S$, then
the sequence $f'_{1}, \dots, f'_{c}$ with $f'_{i} = f_{i}+\sum_{j>i}\alpha_{i,j}f_j$ corresponds to a sequence
of $CI$ operators $\chi'_{c}, \dots, \chi'_{1}$ that is almost regular on $E$. 
\item[\rm (2)] Furthermore, for such $\chi'_{i}$ we have the following property.
Set $E(c)=E$ and $$E(i-1)=E(i)^{\ge 2}[2]/\chi'_{i}E(i)^{\ge 2}[2]$$ for $i\le c$.
Suppose $\reg(E)\leq 1$. Set $\nu=(\alpha^{\vee})^{-1}$. Then $\chi'_{c}$ is a non-zerodivisor on $\ext_R^{\ge 2}(N,k)$,
and more generally
 $\chi'_{i}=\sum_{j}\nu_{i,j}\chi_{j}$ is a non-zerodivisor on 
$E(i)^{\ge 2}[2]$ for every $i$. \qed
\end{mylist}
\end{prop}

We say that $f'_{1}, \dots, f'_{c}$ with $f'_{i} = f_{i}+\sum_{j>i}\alpha_{i,j}f_j$ are {\it generic}
 for $N$ if $(\alpha_{i,j})\in{\mathcal Z}$ in the sense above. 

\myspace  
\noindent {\prooffont Proof of Theorem~\ref{epitwo}:}
To simplify the notation, we may begin by replacing 
$N$ by its $(\reg(\ext_{R}(N,k)) -1)$-st syzygy, and assume that $\reg (\ext_{R}(N,k)) = 1$.
After a general change of  $\ff c$ we may also assume,
by Lemma~\ref{negreg}, that $\chi_{c}, \dots, \chi_{1}$ is an almost regular sequence on $\ext_{R}(N,k)$.
By Proposition~\ref{higherSyzygies} it suffices to treat the case
$r=2c$. Set $M=\syz_{2c}^R(N)$.

Let $(\F, \delta)$ be the minimal free resolution
of $N':=\syz_2^R(N)$, so that $M = \ker(\delta_{2c-3})$.
Since
 $N$ has finite projective dimension over $S$,
the module $N'$ also has finite projective dimension over $S$.

Let $(\tilde \F, \tilde\delta)$ be a lifting of $\F$ to 
$\tilde R:= S/(\ff {c-1})$, and let $\tilde t_{c}$ be the lifted CI operator. 
Set $(\widetilde{\G} ,\widetilde{\delta})= \ker( \tilde t_{c})$. 
By Proposition~\ref{whygeneric},  $\chi_{c}$ is a monomorphism on $\ext_R(N', k)=\ext_{R}^{\geq 2}(N,\allowbreak k)[2]$.
 Since $\chi_{c}$ is induced by $t_{c}$, 
Nakayama's Lemma implies
that $t_{c}$ is surjective, so
in particular $\F_{\geq 2c-2}\to \F_{\geq 2c-4}$ is surjective,
as required for Condition (1) in \ref{wellb} for $c>1$.

Using Nakayama's Lemma again, we see that the lifted CI operator 
$\widetilde{t}_{c}$ is also 
an epimorphism.
Propositions~\ref{second} and \ref{first} show that
$\widetilde{\G}$ is a minimal free resolution of $N'$ over $\tilde R$, and $\F$ is obtained from  
$\widetilde{\G}$ by the Shamash construction \ref{standard}. Hence
  $$ \ext_{\widetilde R}(N',k) = \ext_{R}(N',k)/\chi_{c}\ext_{R}(N',k)\, ,$$ and therefore
$$
\ext_{\widetilde R}(N',k) = \Big(\ext_{R}^{\geq 2}(N,k)\Big/\chi_{c}\ext_{R}^{\geq 2}(N,k)\Big)[2].
$$ 
By
Proposition~\ref{whygeneric}  we conclude  that
$\ext^{\geq 2}_{\tilde R}(N',k)$ has regularity $\leq 1$ over $k[\chi_{1},\dots,\allowbreak  \chi_{c-1}]$.

Suppose now that $c=1$, so that $M = N'$ is the second syzygy of $N$. 
 In this case $\tilde R=S$, and by hypothesis $M = N'$ has finite projective dimension over $S$. Therefore, $\ext_{S}(M,k)$ is a module of finite length. Since it has regularity $\le 1$ (as a module over $k$), it follows that it is
zero except in degrees $\le 1$, that is, the projective dimension of $M$ over $\tilde R$ is $\le 1$. By Proposition~\ref{stableCodimOne}, $M$ is a pre-stable syzygy.

Next suppose that $c>1$.
We have $\ker (\widetilde{\delta}_{2c-3})=\syz^{\tilde R}_{2(c-1)}(N')$,
and by induction on $c$ this is a pre-stable syzygy, verifying Condition (2) in \ref{wellb}. Thus $M$ is
a pre-stable syzygy.
\qed

\b{remark}\label{resolvable}
Some caution is necessary if we wish to work in the graded case (see for example \cite{Pe} for graded resolutions).
Suppose that $S = k[x_{1},\dots, x_{n}]$ is
a standard graded polynomial ring with generators $x_{i}$ in degree 1. Let $f_{1}, \dots, f_{c}$ be a homogeneous regular
sequence, and  set $R=S/(f_{1}, \dots, f_{c})$. Let $N$ be a finitely generated graded $R$-module. When all the $f_{i}$ have the same degree, so that a general linear scalar combination of them is still homogeneous, then   Proposition~\ref{whygeneric} and {Theorem}~\ref{epitwo}
hold for $E=\ext_{R}(N,k)$ verbatim, without first localizing at the maximal ideal. 
However when the $f_1, \dots f_c$ have distinct degrees, there may be no homogeneous linear combination of the $f_{j}$ that corresponds to an eventually surjective CI operator, as
can be seen from the following example.
Let
$R = k[x,y]/(x^2, y^3)$ and consider the module 
$N= R/x \oplus R/y$. Over the local ring
$S_{(x,y)}/(x^2, y^3)$ the CI operator corresponding to $x^{2}+y^{3}$ is eventually surjective. However, the minimal $R$-free resolution of $N$ is the direct sum of  the free resolutions of $R/x$ and $ R/y.$
The  CI operator corresponding to $x^2$ vanishes on the minimal free resolution of $R/y$.
The  CI operator corresponding to $y^3$ vanishes on 
the minimal free resolution of $R/x$, and thus the CI operator
corresponding to
$y^{3}+ax^{3}+bx^{2}y$, for any $a,b$,  does too. 
\end{remark}

\section{The Box complex}\label{Sbox}

Suppose that $f\in S$ is a non-zerodivisor. Given an $S$-free
resolution of an $S/(f)$-module $L$ and a homotopy for $f$, we will construct an $S$-free resolution of the
second syzygy $\syz_2^{S/(f)}(L)$ of $L$ as an $S/(f)$-module, and also a homotopy for $f$ on it.

\theoremstyle{definition}
\newtheorem{bconstruction}[thm]{Box Construction}  
\b{bconstruction}\label{boxc}
Suppose that $f\in S$ is a non-zerodivisor, and that

\myspace
\eq \label{fseven} \  \end{equation} \vglue -1.5cm \ 
\begin{center}
\begin{tikzpicture}  
 [every node/.style={scale=0.9},  auto]
\node(00) {$ Y_{4}$};
\node(01)[node distance=2cm, right of=00]{$Y_{3}$};
\node(02)[node distance=2cm, right of=01]{$Y_{2}$};
\node(03)[node distance=2cm, right of=02]{$Y_{1}$};
\node(04)[node distance=2cm, right of=03]{$Y_{0}$};
\node(05)[node distance=1cm, left of=00] {$\cdots \to$};
\node(06)[node distance=1.5cm, left of=05] {$\Y:\ $};
\draw[->] (00) to node [swap] {$\partial_4$} (01);
\draw[->] (01) to node [swap]{$\partial_3$} (02);
\draw[->] (02) to node [swap]{$\partial_2$} (03);
\draw[->] (03) to node [swap]{$\partial_1$} (04);
\draw[->,  bend right=45] (03) to node[swap] {$\theta_1$} (02);
\draw[->, bend right=45] (04) to node[swap] {$\theta_0$} (03);
\draw[->,  bend right=45] (02) to node[swap] {$\theta_2$} (01);
\draw[->, bend right=45] (01) to node[swap] {$\theta_3$} (00);
\draw[->,  bend left=45] (04) to node {$ \tau_0$} (01);
\draw[->,  bend left=45] (03) to node {$ \tau_1$} (00);
\end{tikzpicture}
\end{center}
\noindent is an $S$-free resolution of a module $L$ annihilated by $f$, with  homotopies $\{\theta_i:\ Y_i\to Y_{i+1}\}_{i\ge 0}$ and  higher homotopies $\tau_0:\ Y_0\to Y_3$ and $\tau_1: Y_1\to Y_4$ for $f$,
so that $\partial_{3}\tau_0 + \theta_1 \theta_0=0$ and $\tau_0 \partial_1+\theta_2\theta_1+\partial_4\tau_1=0$.
We call
the mapping cone 
\vglue .1cm
\eq \label{feight}
\ \end{equation}
\vglue -1.5cm \ 
\begin{center}
\begin{tikzpicture}  
 [every node/.style={scale=0.9},  auto]
\node(pre)[node distance = 1cm, left of=00]{};
\node(00){$\Box(\Y):\quad \quad\to Y_{4}$};
\node(01)[node distance=3.8cm, right of=00]{$Y_{3}$};
\node(02)[node distance=2.5cm, right of=01]{$Y_{2}$};
\node(11a)[node distance=.7cm, below of=01]{$\oplus$};
\node(12a)[node distance=.7cm, below of=02]{$\oplus$};
\node(11)[node distance=1.4cm, below of=01]{$Y_{1}$};
\node(12)[node distance=1.4cm, below of=02]{$Y_{0}$};
\draw[->] (00) to node  {$\partial_4$} (01);
\draw[->] (01) to node [ pos=.4]{$\partial_3$} (02);
\draw[->] (11) to node [pos=.6]{$\partial_1$} (12);
\draw[->] (11) to node [above=2pt, left=1pt] {$\psi$} (02);
\end{tikzpicture}
\end{center}
\vglue .1cm
\noindent of the map $\psi:=\theta_1: \   \Y_{\le 1}[1]\to \YA_{\ge 2}$ the {\it box complex} and denote it $\Box(\Y)$. 
\end{bconstruction}

\theoremstyle{plain}
\newtheorem{bprop}[thm]{Box Proposition} 

\b{bprop}\label{boxthm}
With notation as above, the box complex $\Box(\Y)$
is an $S$-free resolution of the module $\ker\Big( S/(f)\otimes Y_{1} \,
{\mathop{{\mathrel{\smash-}}{\mathrel{\mkern-8mu}}\a\a
{\mathrel{\smash-}}{\mathrel{\mkern-8mu}}\rightarrow}\limits^{S/(f)\otimes \partial_1}}
\,S/(f)\otimes Y_{0}\Big)$, the second $S/(f)$-syzygy of $L$.
Moreover, the maps
\eq   \label{hobox}
\b{pmatrix}
\theta_2&\tau_0\\ 
\partial_2 & \theta_0\end{pmatrix},\ 
 (\theta_3,\tau_1), \ \theta_4,\ \dots
\end{equation}
 give a homotopy for multiplication by $f$ on $\Box(\Y)$ as shown in diagram {\rm (\ref{fnine})}:

\eq \label{fnine} 
\  \end{equation} \vglue -1.5cm \ 
\begin{center}
\begin{tikzpicture}  
 [every node/.style={scale=0.9},  auto]
\node(00){$Y_{4}$};
\node(pre)[node distance = 2cm, left of=00]{};
\node(01)[node distance=3cm, right of=00]{$Y_{3}$};
\node(02)[node distance=5cm, right of=01]{$Y_{2}$};
\node(11a)[node distance=.7cm, below of=01]{$\oplus$};
\node(12a)[node distance=.7cm, below of=02]{$\oplus$};
\node(11)[node distance=1.4cm, below of=01]{$Y_{1}$};
\node(12)[node distance=1.4cm, below of=02]{$Y_{0}\, .$};
\draw[->] (pre) to node {} (00);
\draw[->] (00) to node {$\partial_4$} (01);
\draw[->] (01)[bend left = 10] to node {$\partial_3$} (02);
\draw[->] (11)[bend right=10] to node [swap]{$\partial_1$} (12);
\draw[->] (11)[bend right=20] to node[ pos=.2, anchor=south west] {$\psi$} (02);
\draw[->, red] (11)  to node {$\tau_1$} (00);
\draw[->, bend right=45, red] (02)  to node [swap] {$\theta_2$} (01);
\draw[->, bend left=45, red] (12)  to node  {$\theta_0$} (11);
\draw[->, red] (12) to node [right=50pt, above=-15pt] {$\tau_0$} (01);
\draw[->, bend right=45, red] (01)  to node [anchor=south,pos=.5] {$\theta_3$} (00);
\draw[->, red, bend right=20] (02) to node[right=15pt, below=-5pt] {$\partial_2$} (11);
\end{tikzpicture}
\end{center}
\end{bprop}

A similar formula yields a full system of higher homotopies on $\Box(\Y)$ from higher
homotopies on $\Y$, but we will not need this.

\myspace
\proof 
The following straightforward computation shows that the maps in (\ref{hobox})
 are homotopies for $f$ on $\Box(\Y)$:
\myalign\label{boxhomotopy}
{\b{pmatrix}
\partial_3&\theta_1\\ 
0 & \partial_1\end{pmatrix}
\b{pmatrix}
\theta_2&\tau_0\\ 
\partial_2 & \theta_0\end{pmatrix}
}&{=\b{pmatrix}
\partial_3\theta_2+\theta_1\partial_2&\partial_3\tau_0+\theta_1\theta_0\\ 
\partial_1\partial_2 &\partial_1 \theta_0\end{pmatrix}=
\b{pmatrix}
f&0\\ 
0 &f\end{pmatrix}}
\\ \nonumber
{\b{pmatrix}
\partial_1&\theta_1\\ 
0 & \partial_1
\end{pmatrix}
+\b{pmatrix}\partial_4\theta_3& \partial_4\tau_1\\  0&0\end{pmatrix}
}&{=\b{pmatrix}
\theta_2\partial_3+\partial_4\theta_3&\theta_2\theta_1+\tau_0 \partial_1+\partial_4\tau_1\\ 
\partial_2 \partial_3& \partial_2\theta_1+\theta_0 \partial_1\end{pmatrix}
=\b{pmatrix}
f&0\\ 
0 &f\end{pmatrix}\, .
}
\end{align}

\myspace
Next we will prove
that $\Box(\Y)$ is a resolution.
There is a short exact sequence of complexes
$$
0\to \Y_{\ge 2}\to \Box(\Y)\to {\Y}_{\le 1}\to 0\, ,
$$
 so $H_{i}(\Box(\Y)) = H_{i}(\Y_{\ge 2}) = 0$ for $i\ge 2$ since ${\Y}_{\le 1}$ is a two-term complex.
If $(v,w)\in Y_3\oplus Y_1$ is a cycle, then applying the homotopy maps
in (\ref{hobox}) we get $$(fv,fw)=(\partial_4\theta_3(v)+\partial_4\tau_1(w),0)\,.$$ Since $f$ is a non-zerodivisor,
it follows that $w=0$. Thus $v$ is a cycle in $\Y_{\ge 2}$, which is acyclic, so $v$ is a boundary in $\Y_{\ge 2}$. Hence, the complex $\Box(\Y)$ is acyclic.

To simplify notation, we write  $\overline{-}$  for the functor 
$S/(f)\otimes -$ and set $\psi=\theta_1$. To complete the proof we will show that $H_{0}(\Box(\Y)) = \ker (\overline {\partial_1}: \   \overline Y_{1}\to \overline Y_{0})$.
Since we have a homotopy for $f$ on $\Y$, we see that $f$ annihilates the module resolved by $\Y$.
Therefore, $H_0(\Y)=H_1(\overline{\Y})$.
The complex $\overline{\Box(\Y)}$ is the mapping cone $\MM(\overline{ \psi}\otimes S/(f))$, where $\bar \psi=\psi\otimes S/(f)$, so there is an exact sequence of complexes
$$
0\to \overline{\Y}_{\ge 2}\to \overline{\Box(\Y)}\to \overline{ {\Y}}_{\le 1}\to 0\, .
$$

Since $\Y$ is a resolution, $H_0(\Y_{\ge 2})$ is contained in the free $S$-module $Y_{1}$. 
Thus
$f$ is a non-zerodivisor on $H_0(\Y_{\ge 2})$ and  $\overline{\Y}_{\ge 2}$ is acyclic.
Therefore, the long exact sequence for the mapping cone yields
$$0\to H_1(\MM(\overline{\psi}))\to H_1(\overline{{\Y}}_{\le 1})\ar{\bar\psi} H_0(\overline{ \Y}_{\ge 2})\, .$$
 It suffices to prove that the map induced on homology by
$\overline {\psi}$ is $0$. Let $u\in Y_{1}$ be such that $\overline u\in \ker (\overline \partial_1)$, so
$ \partial_1 (u) = fy$ for some $y\in Y_{0}$.
We also have $fy=  \partial_1\theta_0 (y)$, so $u-\theta_0 (y) \in\ker(\partial_1)$. Since $\Y$ is acyclic
$u= \theta_0 (y)+\partial_2 (z)$ for
some $z\in Y_{2}$. Applying $\psi$ we get
\b{align*}
\psi (u) &= \theta_1\theta_0 (y) +\theta_1\partial_2 (z) \\ 
&= -\partial_3\tau_0( y)+(fz-\partial_3\theta_2(z))\\ 
&=-\partial_{3}\Big(\tau_0(y)+\theta_2(z)\Big)+fz,
\end{align*}
so the map induced on homology by
$\overline {\psi}$ is $0$ as desired.
\qed

\myspace
Proposition~\ref{boxthm} has a partial converse that we will use in the proof of
Theorem~\ref{stable}.

\b{prop}\label{stableExactness} 
Let $f\in S$ be a non-zerodivisor and set $R = S/(f)$.
Let

\myspace
\begin{center}
\begin{tikzpicture}  
 [every node/.style={scale=0.9},  auto]
\node(prepre)[node distance = 7cm, left of=11a]{$\Box(\Y):$};
\node(pre)[node distance = 1cm, left of=00]{$\cdots$};
\node(00){$Y_{4}$};
\node(01)[node distance=2.5cm, right of=00]{$Y_{3}$};
\node(02)[node distance=2.5cm, right of=01]{$Y_{2}$};
\node(11a)[node distance=.7cm, below of=01]{$\oplus$};
\node(12a)[node distance=.7cm, below of=02]{$\oplus$};
\node(11)[node distance=1.4cm, below of=01]{$Y_{1}$};
\node(12)[node distance=1.4cm, below of=02]{$Y_{0}$};
\draw[->] (pre) to node {} (00);
\draw[->] (00) to node  {$\partial_4$} (01);
\draw[->] (01) to node [ pos=.4]{$\partial_3$} (02);
\draw[->] (11) to node [pos=.6]{$\partial_1$} (12);
\draw[->] (11) to node [above=2pt, left=1pt] {$\psi$} (02);
\end{tikzpicture}
\end{center}
\myspace

\noindent be an $S$-free resolution of a module annihilated by $f$. Set
$\theta_1:=\psi$, and with notation as in diagram {\rm (\ref{fnine})}, suppose that 
$$
\b{pmatrix}
\theta_2&\tau_0\\ 
\partial_2 & \theta_{0}\end{pmatrix}
$$
is the first map of a homotopy for 
multiplication by $f$ on $\Box(\Y)$. If the cokernels of $\partial_2$ and of $\partial_3$ are $f$-torsion free,
then the following complex is exact:
\eq \label{fnineseq}
  \dots \to\,Y_4\ar{\partial_4} Y_3 \ar{\partial_3} Y_2  \ar{\partial_2} Y_{1} \ar{\partial_1} Y_{0}\, ,
\end{equation}
and there are homotopies for $f$ as in {\rm (\ref{fseven})}.
\end{prop}

\proof We first show that the sequence is a complex.
 The equation $\partial_3\partial_4 = 0$ follows from our hypothesis. Let $(\theta_3, \tau_1): \   Y_3\oplus Y_1\to Y_2$
be the next map in the homotopy for $f$. 
To show that
$\partial_2\partial_3 = 0$ and $\partial_1\partial_2 = 0$, use the homotopy equations
\b{align*}
0\theta_3+\partial_2\partial_3=0&: \ Y_3\to Y_1\\ 
\partial_1\partial_2=0&:\ Y_2\to Y_0\, .
\end{align*}

The equalities in (\ref{boxhomotopy})  imply that
$\theta_0: \   Y_0\to Y_1$, $\psi=\theta_1: \   Y_1\to Y_2$,  $\theta_2: \   Y_2\to Y_3$, and $\theta_3: \   Y_3\to Y_4$
 form the beginning
of a homotopy for $f$ on (8.9). Thus (8.9) becomes exact after inverting $f$. 
The exactness of (8.9) is equivalent to the statement that the induced maps
$\coker(\partial_3) \to Y_1$ and $\coker(\partial_2) \to Y_2$ are monomorphisms.
Since this is true after inverting $f$, and since the cokernels are $f$-torsion free by
hypothesis, exactness holds before inverting $f$ as well.
\qed

\section{From Syzygies to Higher Matrix Factorizations}\label{SmakeMF}

 Higher matrix factorizations arising from pre-stable syzygies have an additional property. We introduce the
concept of a pre-stable matrix factorizations, which captures that property.

\b{definition}\label{stableMFdef}
 A   higher matrix
factorization $(d,h)$  is a {\it pre-stable matrix factorization} if, in the notation of~\ref{standardnotation},  for each $p=1,\dots, c$ the element $f_{p}$ is a non-zerodivisor on the cokernel of the composite map 
$$
R(p-1)\otimes A_{0}(p-1)  \hookrightarrow R(p-1)\otimes A_{0}(p) \ar{h_p} R(p-1)\otimes A_{1}(p) \ar{\pi_p} R(p-1)\otimes B_{1}(p).
$$
If $S$ is Cohen-Macaulay then we say that the higher matrix factorization $(d,h)$ is a  {\it stable matrix factorization} if the cokernel of the composite map above
is a maximal Cohen-Macaulay $R(p-1)$-module.
\end{definition}

\myspace
The advantage of stable matrix factorizations over pre-stable matrix factorizations is that if $g\in S$ is an element such that $g,f_{1}, \dots, f_{c}$ is a regular sequence and $(d,h)$ is a stable matrix factorization, then $\Big(S/(g) \otimes d, \ S/(g)\otimes h\Big)$ is again a stable matrix factorization. We do not know of pre-stable matrix factorizations that are not stable.

\b{thm}\label{eventually}
Suppose that $f_{1}, \dots, f_{c}$ is a regular sequence in a
 local  ring $S$, and set
 $R= S/(f_{1}, \dots, f_{c})$. If $M$ is a pre-stable  syzygy over $R$ with respect to
 $\ff c$, then $M$ is the HMF module of a minimal pre-stable  matrix factorization 
 $(d,h)$ such that $d$ and $h$ are liftings to $S$ of the first two differentials   in the minimal $R$-free resolution of $M$. If $M$ is a stable syzygy, then $(d,h)$ is stable as well.
 \end{thm}

Combining Theorem~\ref{eventually} and Theorem~\ref{epitwo}
we obtain the following more precise version of Theorem~\ref{highsyzthm}  in the introduction.

\b{cor}\label{highsyzy}
Suppose that $f_{1}, \dots, f_{c}$ is a regular sequence in a
 local ring $S$ with infinite residue field $k$, and set
 $R= S/(f_{1}, \dots, f_{c})$.  Let $N$ be an $R$-module with finite projective dimension over $S$. 
 There exists a non-empty Zariski open dense set $\mathcal Z$ of matrices $(\alpha_{i,j})$
with entries in $k$ such that 
for every $$
r\ge 2c-1+ \reg(\ext_{R}(N,k))\,
$$
the syzygy $\syz_r^R(N)$  is the module of a minimal pre-stable matrix factorization with respect to
 the regular sequence $\{\,f'_{i} = \sum_{j}\alpha_{i,j}f_j\,\}$.\qed
\end{cor}

\noindent{\prooffont Proof of Theorem~\ref{eventually}:}
The proof is by induction on $c$. If $c=0$, then $M=0$ so we are done.

Suppose $c\ge 1$. We use the notation of Definition~\ref{wellb}. By assumption,  the CI operator $t_c$ is surjective on a minimal $R$-free resolution $(\F,\delta)$ of
a module $L$ of which $M$ is the second syzygy. Let $(\tilde\F, \tilde\delta)$ be a lifting of $(\F,\delta)$ to $R'=S/(f_1,\dots ,f_{c-1})$. Since $S$ is local, the lifted CI operator
  $\widetilde{t}_{c}:=(1/f_{c})\tilde\delta^{2}$ is also surjective, and we set $(\widetilde{\G},\tilde\partial):=\ker (\tilde t_c)$.
  By Propositions~\ref{first}, $\F$ is the result of applying the Shamash construction to  $\tilde \G$.
Let $\tilde B _1(c)$ and $\tilde B_0(c)$ be 
the liftings to $R'$ of $F_{1}$ and $F_{0}$ respectively.
By  Propositions~\ref{second} and \ref{first} the minimal $R'$-free resolution of $L$ has the form
\eq \label{sequence}
\dots \to \tilde{G}_{4}\ar{\widetilde{\partial}_4} \tilde A_1(c-1):=\tilde G_{3}\ar{\widetilde{\partial}_3}\tilde A_0(c-1):=\tilde G_{2}\ar{\tilde\partial_2} \tilde B_{1}(c)\ar{\tilde b} \tilde B_{0}(c)
\,,\end{equation}
where ${\tilde b}:=\tilde\partial_1 $, $\widetilde{\partial}_2$, $\tilde \partial_3$, $\widetilde{\partial}_4$ are the  liftings of the differential in $\F$.

Since $L$ is annihilated by $f_{c}$ there exist homotopy maps $\tilde\theta_0 ,\tilde \psi:=\tilde\theta_1, \tilde\theta_2$ and a higher homotopy $\tilde\tau_0$ so that on
\myspace
\eq \label{ften} 
\  \end{equation} \vglue -1.6cm \ 
\begin{center}
\begin{tikzpicture}  
 [every node/.style={scale=0.9},  auto]
\node(00) {$\cdots  $};
\node(01)[node distance=2cm, right of=00]{$\tilde{G}_{3}$};
\node(02)[node distance=2cm, right of=01]{$\tilde G_{2}$};
\node(03)[node distance=2cm, right of=02]{$ \tilde B_{1}(c)$};
\node(04)[node distance=2.5cm, right of=03]{$\tilde B_{0}(c)$};
\node(05)[node distance=1cm, left of=00] {};
\node(06)[node distance=1.5cm, left of=05] {};
\draw[->] (00) to node [swap]  {$\widetilde{\partial}_4$ } (01);
\draw[->] (01) to node [swap]{$\widetilde{\partial}_3$} (02);
\draw[->] (02) to node [swap]{$\tilde \partial_2$} (03);
\draw[->] (03) to node [swap]{$\tilde b=\tilde\partial_1$} (04);
\draw[->, red, bend right=45] (03) to node[swap] {$\tilde \psi=\tilde\theta_1$} (02);
\draw[->, red, bend right=45] (04) to node[swap] {$\tilde \theta_0$} (03);
\draw[->, red, bend right=45] (02) to node[swap] {$\tilde\theta_2$} (01);
\draw[->, red, bend left=45] (04) to node[swap] {$\tilde \tau_0$} (01); 
\end{tikzpicture}
\end{center}

\noindent we have
\myalign \label{homeq}
\nonumber \tilde\partial_1\tilde \theta_0&=f_c\id\\ 
 \tilde \partial_2   \tilde\theta_1+\tilde\theta_0 \tilde\partial_1 &=f_{c}\id\\ 
\nonumber \widetilde{\partial}_3  \tilde\theta_2+\tilde\theta_1 \tilde \partial_2&=f_{c}\id\\ 
\nonumber \widetilde{\partial}_3\tilde \tau_0 + \tilde\theta_1\tilde \theta_0&=0
\, .
\end{align}
Proposition~\ref{boxthm} implies that the minimal free resolution of $M$ over $R'$
has the form 
\vglue .05cm
\eq \label{feleven} \  \end{equation} \vglue -1.5cm \ 
\begin{center}
\begin{tikzpicture}  
 [every node/.style={scale=0.9},  auto]
\node(pre)[node distance = 1.5cm, left of=00]{$\cdots$};
\node(00){$\tilde G_{4}$};
\node(01)[node distance=2.5cm, right of=00]{$\tilde G_{3}$};
\node(02)[node distance=2.5cm, right of=01]{$\tilde G_{2}$};
\node(11a)[node distance=.7cm, below of=01]{$\oplus$};
\node(12a)[node distance=.7cm, below of=02]{$\oplus$};
\node(11)[node distance=1.4cm, below of=01]{$\tilde B_{1}(c)$};
\node(12)[node distance=1.4cm, below of=02]{$\tilde B_{0}(c)$};
\draw[->] (pre) to node {} (00);
\draw[->] (00) to node  {$\tilde\partial_4$} (01);
\draw[->] (01) to node [ pos=.4]{$\widetilde{\partial}_3$} (02);
\draw[->] (11) to node [pos=.6]{$\tilde  b $} (12);
\draw[->] (11) to node [above=3pt, left=2pt] {$ \tilde \psi$} (02); 
\end{tikzpicture}
\end{center}
\vglue .05cm
\noindent 
Using this structure we change the lifting of the differential $\delta_{3}$ so that 
$$\tilde \delta_{3}=\b{pmatrix}\widetilde{\partial}_3& \tilde\psi\\  0& \tilde b\end{pmatrix}\,.$$
Note that the differential $\tilde \partial$ on $\tilde {\G}_{\ge 2}$ has not changed.

Set $M'=\coker(\tilde {\G}_{\ge 2})=\syz_2^{R'}(L)$. Since $M'$ is a pre-stable syzygy, the induction hypothesis implies
that $M'$ is the HMF module of a higher matrix factorization $(d',h')$ with respect to $f_1,\dots ,f_{c-1}$
so that the differential $\tilde G_{3} \to \tilde G_{2}$ is $\tilde \partial_3=d'\otimes R'$
and the differential $\tilde G_{4} \to \tilde G_{3}$ is $\tilde \partial_4=h'\otimes R'$.
Thus, there exist free $S$-modules  $A_1'(c-1)$ and $A_0'(c-1)$ with filtrations so that
$$\tilde G_{3}=A_1'(c-1)\otimes R'\quad\hbox{and}\quad \tilde G_{2}=A_0'(c-1)\otimes R'\,.$$

We can now define a higher matrix factorization for $M$.
Let $B_1(c)$ and $B_0(c)$ be free $S$-modules such that $\tilde B_0(c)=B_0(c)\otimes R'$ and
$\tilde B_1(c)=B_1(c)\otimes R'$. For $s=0,1$, we consider free $S$-modules $A_1$ and $A_0$ with filtrations  such that $A_s(p)=A_s'(p)$ for $1\le p\le c-1$ and
$$A_s(c)=A'_s(c-1)\oplus B_s(c)\,.$$  
We define the map $d: \  A_{1}\to A_{0}$ 
to be 
\eq \label{definedif}
A_1(c) = A_{1}(c-1) \oplus B_{1}(c)\ 
{\mathop{{\mathrel{\smash-}}{\mathrel{\mkern-8mu}}\a\a\a
{\mathrel{\smash-}}{\mathrel{\mkern-8mu}}\rightarrow}\limits^{\b{pmatrix}
d'&\psi_c\\  
0&b_{c}\end{pmatrix}
}}\ 
A_{0}(c-1) \oplus B_{0}(c)=A_0(c)
\end{equation}
\vglue .1cm
where $b_{c}$ and $\psi_c$ are arbitrary lifts to $S$ of $\tilde b$ and $\tilde \psi$.
For every $1\le p\le c-1$, we set $h_p=h_p'$.
Furthermore,
we define $$h_c: \   A_{0}(c) = A_{0}\to A_{1}(c) = A_{1}$$ to be 
\eq \label{defineho}
A_0(c) = A_{0}(c-1) \oplus B_{0}(c) \ 
{\mathop{{\mathrel{\smash-}}{\mathrel{\mkern-8mu}}\a\a\a
{\mathrel{\smash-}}{\mathrel{\mkern-8mu}}\rightarrow}\limits^{
\b{pmatrix}\theta_2 &\tau_0 \\  \partial_2 &\theta_0
\end{pmatrix}
}}\ 
A_{1}(c-1) \oplus B_{1}(c)=A_1(c)
\end{equation}
\vglue .1cm
where $\theta_2, \partial_2,\theta_0, \tau_0$ are arbitrary lifts to $S$ of $\tilde\theta_2, \tilde \partial_2, \tilde \theta_0, \tilde \tau_0$
respectively.

We must verify
conditions (a) and (b) of   Definition~\ref{matrixfiltration}. Since
$(d',h')$ is a higher matrix factorization, we need only check 
\b{align*}
& d h_c \equiv f_{c}\, \Id_{A_{0}(c)} \hbox{ mod}(f_{1}, \dots, f_{c-1})A_0(c)\\ 
&\pi_c h_c d \equiv f_{c}\pi_c\hbox{ mod}(f_{1}, \dots, f_{c-1})B_1(c)\, .
\end{align*}
Condition (a) holds because
\vglue .1cm
$$
\b{pmatrix}d'&\psi\\  0&b_{c}\end{pmatrix} 
\b{pmatrix}\theta_2 &\tau_0 \\  \partial_2 &\theta_0\end{pmatrix} 
= 
\b{pmatrix}
d'\theta_2+\theta_1 \partial_2& d'\tau_0 + \theta_1 \theta_0 \\ 
\partial_1\partial_2 & \partial_1\theta_0
\end{pmatrix} 
\equiv 
\b{pmatrix}
f_{c} &0\\ 
0& f_{c}
\end{pmatrix} 
$$
\vglue .2cm
by (\ref{homeq}).
Similarly, Condition (b) is verified by the computation
\vglue .1cm
$$
\b{pmatrix}\theta_2 &\tau_0 \\  \partial_2 &\theta_0\end{pmatrix} 
\b{pmatrix}d'&\psi\\  0&b_{c}\end{pmatrix}  
= 
\b{pmatrix}
\theta_2 d' & \theta_2 \theta_1 +\tau_0 \partial_1 \\ 
\partial_2 d' & \partial_2\theta_1 +\theta_0 \partial_1
\end{pmatrix} 
\equiv 
\b{pmatrix}
* & *  \\ 
0& f_{c}
\end{pmatrix} .
$$
\vglue .2cm

Next we  show that the  higher matrix factorization that we have constructed is pre-stable.
Consider the complex (\ref{sequence}), which 
is a free resolution of $L$
over $R'$. It follows that $$\coker\Big(\widetilde A_{0}(c-1) \ar{\tilde \partial_2} 
\widetilde B_{1}(c)\Big)\cong \im(\widetilde{\partial}_1)\subset\widetilde{B}_0(c)$$ has no $f_{c}$-torsion, verifying the pre-stability condition.

It remains to show that  $d$ and $h$ are liftings to $S$ of the first
two differentials 
 in the minimal $R$-free resolution of $M$.

By (\ref{ften}) and  Proposition~\ref{boxthm} we have the following homotopies on 
the minimal $R'$-free resolution of $M$:
\begin{center}
\begin{tikzpicture}  
 [every node/.style={scale=0.8},  auto]
\node(pre)[node distance = 2cm, left of=00]{};
\node(00){$\widetilde{G}_{4}$};
\node(01)[node distance=3cm, right of=00]{$\widetilde{G}_{3}$};
\node(02)[node distance=5cm, right of=01]{$\widetilde{G}_{2}$};
\node(11a)[node distance=.7cm, below of=01]{$\oplus$};
\node(12a)[node distance=.7cm, below of=02]{$\oplus$};
\node(11)[node distance=1.4cm, below of=01]{$\widetilde{B}_{1}(c)$};
\node(12)[node distance=1.4cm, below of=02]{$\widetilde{B}_{0}(c)\, .$};
\draw[->] (pre) to node {} (00);
\draw[->] (00) to node {$\widetilde {\partial}_4$} (01);
\draw[->] (01)[bend left = 10] to node {$\widetilde {\partial}_3$} (02);
\draw[->] (11)[bend right=10] to node [swap]{$\widetilde  b=\tilde\partial_1$} (12);
\draw[->] (11)[bend right=20] to node[ pos=.2, anchor=south west] {$\widetilde \psi$} (02);
\draw[->, bend right=45, red] (02)  to node [swap] {$\widetilde \theta_2$} (01);
\draw[->, bend left=45, red] (12) to node  {$\widetilde \theta_0$} (11);
\draw[->, red] (12) to node [right=50pt, above=-15pt] {$\widetilde \tau_0$} (01);
\draw[->, red, bend right=20] (02) to node [right=15pt, below=-5pt] {$\widetilde \partial_2$} (11);
 \end{tikzpicture}
\end{center}
\vglue -3cm
\eq \label{ftwelve} \  \end{equation} 
\vglue 1.5cm
\noindent The minimal $R$-free resolution of $M$ is obtained from the resolution above by applying the Shamash construction. Hence,
the first two differentials are
$$R\otimes 
 \b{pmatrix}\widetilde{\partial}_3&\widetilde{\psi}\\   
0&\widetilde{b}\end{pmatrix} 
\quad\hbox{and}\quad R\otimes
 \b{pmatrix}\widetilde{\partial}_4&\widetilde{\theta_2} &\tau_0 \\  0& \widetilde{\partial_2} &\widetilde{\theta_0}\end{pmatrix}
 \, .
$$

By induction hypothesis  $\widetilde{\partial}_3=R'\otimes d_{c-1}$ and  $\widetilde{\partial}_4=R'\otimes h(c-1)$.
By the construction of $d$ and $h$ in (\ref{definedif}), (\ref{defineho}) we see that  $R\otimes d$ and $R\otimes h$
are the first two differentials in the minimal $R$-free resolution of $M$.

Finally, we will prove that if $M$ is a stable syzygy, then $(d,h)$ is stable as well.
The map  $\partial_2$ is the composite map
$$
\ A_{0}(p-1)  \hookrightarrow A_{0}(p) \ar{h_p} A_{1}(p) \ar{\pi_p} B_{1}(p)\, 
$$
 by construction (\ref{defineho}).
By (\ref{sequence}) it follows that if $L$ is a maximal Cohen-Macaulay $R$-module, then
$\coker (\widetilde{\partial}_2)$ is a maximal Cohen-Macaulay $R'$-module,
 verifying the stability condition
for a higher matrix factorization over $R(p-1)$. By induction, it follows that $(d,h)$ is stable.
\qed

\b{remark}
In order to capture structure when minimality is not present,
Definition~\ref{wellb} can be modified as follows.
We extend the definition of syzygies to non-minimal free resolutions: if $({\bf F},\delta)$ is an $R$-free resolution of an $R$-module $P$, then we define
$\syz_{i,{\bf F}}(P)=\im (\delta_{i})$. 
Suppose that $f_{1}, \dots, f_{c}$ is a regular sequence in a
 local ring $S$, and set
 $R= S/(f_{1}, \dots, f_{c})$. Let $(\F, \delta)$ be an $R$-free resolution, and let $M = \im ( \delta_{r})$ for a fixed  $r\geq 2c$. 
 \end{remark}

 We say that $M$ is a {\it pre-stable syzygy in $\F$} with respect to $\ff c$ if either $c=0$ and $M=0$, or $c\ge 1$ and there exists a  lifting $(\tilde\F, \tilde\delta)$ of $(F,\delta)$ to $R'=S/(f_1,\dots ,f_{c-1})$ such that the CI operator $\widetilde{t}_{c}:=(1/f_{c})\tilde\delta^{2}$ is surjective and,
 setting $(\widetilde{\G},\tilde\partial):=\ker (\tilde t_c)$, the module $\im( \widetilde\partial_{r})$ is pre-stable in $\widetilde{\G}_{\ge 2}$ with respect to $\ff {c-1}$. 

 With minor modifications, the proof of Theorem~\ref{eventually} yields the following result:
 Let $\F$ be an $R$-free resolution. 
 If $M$ is a pre-stable $r$-th syzygy in $\F$ with respect to
 $\ff c$ then $M$ is the HMF module of a pre-stable matrix factorization 
 $(d,h)$ such that $d$ and $h$ are liftings to $S$ of the consecutive differentials $\delta_{r+1}$ and $\delta_{r+2}$ in $\F$.
 If $\F$ is minimal then the higher matrix factorization is minimal.

\myspace
We can use the concept of pre-stable syzygy and Proposition~\ref{stableExactness} in order to build the minimal free resolutions of
the modules $\coker (R(p-1)\otimes b_p)$:

   \b{prop}\label{prestable resolutions}
    Let $(d,h)$ be a minimal pre-stable matrix factorization for a regular sequence $f_{1},\dots, f_{c}$ in a local ring $S$, and use  the notation  of {\rm \ref{standardnotation}}. For every $p\le c$, set
$R(p)=S/(f_1,\dots ,f_{p})$ and    $D(p)=\coker (R(p-1)\otimes b_p)\,.$ Then  $$D(p) = \coker(R(p)\otimes b_{p})\, .$$
Let $\T(p)$ be the minimal $R(p)$-free resolution of $M(p)$ from Construction~ {\rm \ref{inres}} and 
Theorem~{\rm \ref{inresthm}}.
The minimal $R(p-1)$-free resolution of $D(p)$ is
$$\V(p-1):\quad \T(p-1)
 \to R(p-1)\otimes B_{1}(p) \,
 {\mathop{{\mathrel{\smash-}}{\mathrel{\mkern-8mu}}\a\a\a
{\mathrel{\smash-}}{\mathrel{\mkern-8mu}}\rightarrow}\limits^{R(p-1)\otimes b_p}}
\, R(p-1)\otimes B_{0}(p)\, ,
$$
where the second differential is
induced by the composite map
$$
\delta:\ A_{0}(p-1)  \hookrightarrow A_{0}(p) \ar{h_p} A_{1}(p) \ar{\pi_p} B_{1}(p)\, .
$$
The minimal $R(p)$-free resolution of $D(p)$
is 
$$\W(p):\quad \T(p)\to R(p)\otimes B_{1}(p) \,
 {\mathop{{\mathrel{\smash-}}{\mathrel{\mkern-8mu}}\a\a\a
{\mathrel{\smash-}}{\mathrel{\mkern-8mu}}\rightarrow}\limits^{R(p-1)\otimes b_p}}
\,R(p)\otimes B_{0}(p)\, ,
$$
where the second differential is given by 
the Shamash construction applied to $\V(p-1)_{\le 3}$.
\end{prop}

\proof   By Theorem~\ref{inresthm} (using the notation in that theorem)
the complex $\T(p)$ is an $R(p)$-free resolution of $M(p)$.
By  Theorem~\ref{intermediate}, the minimal $R(p-1)$-free resolution of $M(p)$ is
\vglue .1cm
\begin{center}
\begin{tikzpicture}  
 [every node/.style={scale=0.9},  auto]
\node(pre)[node distance = 1.8cm, left of=00]{};
\node(00){$T(p-1)_{2}$};
\node(01)[node distance=2.5cm, right of=00]{$T(p-1)_{1}$};
\node(02)[node distance=6cm, right of=01]{$T(p-1)_{0}$};
\node(11a)[node distance=1cm, below of=01]{$\oplus$};
\node(12a)[node distance=1cm, below of=02]{$\oplus$};
\node(11)[node distance=2cm, below of=01]{$R(p-1)\otimes B_{1}(p)$};
\node(12)[node distance=2cm, below of=02]{$R(p-1)\otimes B_{0}(p)\,.$};
\draw[->] (pre) to node {} (00);
\draw[->] (00) to node  {} (01);
\draw[->] (01) to node [ pos=.4]{$R(p-1)\otimes d_{p-1}$} (02);
\draw[->] (11) to node [pos=.6]{$R(p-1)\otimes b_p$} (12);
\draw[->] (11) to node [pos=.6,above=5pt, left=1pt] {$R(p-1)\otimes \psi_p$} (02);
\end{tikzpicture}
\end{center}
\vglue .1cm
\noindent Since $f_p$ is a non-zerodivisor on $M(p-1)$ by Corollary~\ref{mcm} and since the matrix factorization is pre-stable,
we can apply Proposition~\ref{stableExactness},
where the homotopies 
$\theta_{i}$ and $\tau_{i}$ for $f_p$  are chosen 
to be the appropriate components of the map $R(p-1)\otimes h_p$.
We get 
the minimal $R(p-1)$-free resolution  $$\V(p-1):\quad \T(p-1)
 \to R(p-1)\otimes B_{1}(p) \,
 {\mathop{{\mathrel{\smash-}}{\mathrel{\mkern-8mu}}\a\a\a
{\mathrel{\smash-}}{\mathrel{\mkern-8mu}}\rightarrow}\limits^{R(p-1)\otimes b_p}}
\, R(p-1)\otimes B_{0}(p)\, ,
$$
where the second differential is
induced by the composite map
$$
\delta:\ A_{0}(p-1)  \hookrightarrow A_{0}(p) \ar{h_p} A_{1}(p) \ar{\pi_p} B_{1}(p)\, .
$$
Since we have a homotopy for $f_p$ on $R(p-1)\otimes B_{1}(p)\to
R(p-1)\otimes B_{0}(p)\, $ it follows that $D(p) = \coker(R(p)\otimes b_{p})\, .$

We next apply the Shamash construction to the following diagram with homotopies:
\begin{center}
\hglue -2cm\vbox{
\begin{tikzpicture}  
 [every node/.style={scale=0.9},  auto]
\node(00) {$\V(p-1)_{\le 3}:\ \ $};
\node(01)[node distance=2.5cm, right of=00]{$A_1(p-1)'$};
\node(02)[node distance=3cm, right of=01]{$A_0(p-1)'$};
\node(03)[node distance=3cm, right of=02]{$ B_1(p)'$};
\node(04)[node distance=3.5cm, right of=03]{$B_0(p)'\,,$};
\node(05)[node distance=1cm, left of=00] {};
\node(06)[node distance=1.5cm, left of=05] {};
\draw[->] (00) to node  { } (01);
\draw[->] (01) to node [swap]{$d_{p-1}'$} (02);
\draw[->] (02) to node [swap]{$\partial_2$} (03);
\draw[->] (03) to node [swap]{$\partial_1=b_p'$} (04);
\draw[->,  bend right=45] (03) to node[swap] {$\theta_1:=\psi_p'$} (02);
\draw[->,  bend right=45] (04) to node[swap] {$ \theta_0$} (03);
\draw[->,  bend right=45] (02) to node[swap] {$\theta_2$} (01);
\draw[->,  bend left=30] (04) to node[swap] {$ \tau_0$} (01);
\end{tikzpicture}}
\end{center}
\noindent where $-'$ stands for $R(p-1)\otimes -$.
By Proposition~\ref{second} we obtain an exact sequence
$$
R(p)\otimes A_{1}(p) \to R(p)\otimes A_{0}(p)  \to R(p)\otimes B_{1}(p) \to R(p)\otimes B_{0}(p)\, .
$$
It is minimal since $\theta_0$ is induced by $h_p$.
The leftmost differential $$R(p)\otimes A_{1}(p) \ 
{\mathop{{\mathrel{\smash-}}{\mathrel{\mkern-8mu}}\a\a
{\mathrel{\smash-}}{\mathrel{\mkern-8mu}}\rightarrow}\limits^{R(p)\otimes b_p}}
\  R(p)\otimes A_{0}(p)$$ coincides with the first differential in $\T(p)$.
\qed

\myspace

The following result (stated somewhat differently) and the idea of the proof are from \cite[Theorem 7.3]{AGP}. We will use it in Corollary~\ref{zerostop} in order to obtain numerical information about pre-stable matrix factorizations.

\b{prop}\label{strictlyincreasing} 
Let  $f\in S$ be a non-zerodivisor
in a local ring $S$, and let
$\bf F$ be a  minimal free resolution of a nonzero module over $S/(f)$.
 If the CI operator 
$t: F_2\to F_0$ corresponding to $f$ is surjective, then 
$\rank (F_{1})\geq \rank (F_{0})$,
and if equality holds then $\bf F$ is periodic of period $2$ (that is, $\syz_2^{S/(f)}(L)\cong 
L$ where $L=H_0(\F)$). In the latter case, the ranks of the free modules $F_i$ are constant.
\end{prop}

\proof We lift the first two steps of $\bf F$ to $S$ as
$
\widetilde F_{2}\ar{\widetilde \delta_{2}} \widetilde F_{1}\ar{\widetilde \delta_{1}} \widetilde F_{0},
$
so that $\widetilde \delta_{1}\widetilde \delta_{2} = f\widetilde t$. Since $t$ is surjective
and $f$ is in the maximal ideal, $\widetilde t$ is surjective.
Thus the image of $\widetilde \delta_{1}$ contains $f\widetilde F_{0}$, and it follows
 that  $\rank (\widetilde \delta_{1})=\rank (\widetilde F_{0})$. In particular,
$\rank (F_{1})\geq \rank (F_{0})$. 
In case of equality $\widetilde \delta_{1}$ is a monomorphism,
and we can factor the multiplication by $f$ on $\widetilde  F_{0}$ as 
$\widetilde \delta_{1}\widetilde u_{1}$ for some $u_{1}$\ ---\ a matrix factorization of $f$.
Thus the cokernel of $\delta_{1}$ is resolved by the periodic resolution coming from this
matrix factorization, so $\bf F$ is periodic. Then the ranks of the free modules $F_i$ are constant
by \cite[Proposition~5.3]{Ei1}. \qed

\myspace
Using Proposition~\ref{strictlyincreasing},  we get a stronger version of Corollary~\ref{nozeros}  for  pre-stable matrix factorizations.

\b{cor}\label{zerostop}
Let $(d,h)$ be a minimal pre-stable matrix factorization, and use  the notation  of {\rm \ref{standardnotation}}.
Let  $\gamma$ be the minimal number such that $A(\gamma)\not=0$. Then $\hbox{\rm cx}_R(M)=c-\gamma+1$ and 
\b{align*}
\rank(B_1(p))&=\rank(B_0(p))=0\quad \hbox{for  every}\ 1\le p\le \gamma-1\\ 
\rank(B_1(\gamma))&=\rank(B_0(\gamma))>0\\ 
\rank(B_1(p))&>\rank(B_0(p))>0\quad \hbox{for  every}\ \gamma+1\le p\le c\, .
\end{align*}
The multiplicity of $\ext^{even}$ (equal to the  multiplicity of $\ext^{odd}$ and called the Betti degree)
is the 
size of the hypersurface matrix factorization that is the top non-zero part of the higher matrix factorization $(d,h)$. 

For every $p\le \gamma-1$, the projective dimension of $M$ over $R(p)$ is finite and  we have the equality of Poincar\'e series
$$
{\poin}^{R(p)}_{M}(x)=(1+x)^{p}\,{\poin}^{S}_{M}(x)\, .$$
\end{cor}

\proof 
By Corollary~\ref{bettiinf}(2) it follows that $\hbox{\rm cx}_R(M)=c-\gamma+1$.
The definition of a higher matrix factorization shows that $(d_{\gamma}, h_{\gamma})$ is a codimension $1$ matrix factorization, and hence $\rank(B_1(\gamma))=\rank(B_0(\gamma))>0$.

For  $\gamma+1\le p\le c$, we apply Proposition~\ref{prestable resolutions}.
Since  $\V(p-1)$ is a free resolution, it follows that $B_{0}(p)=0$ if $B_{1}(p)=0$.
But  $B_{0}(p)=B_{1}(p)=0$ implies that $M(p)=M(p-1)$ which contradicts to the fact that
$f_p$ annihilates $M(p)$ and $f_p$ is a non-zerodivisor on $M(p-1)$ by Corollary~\ref{mcm}.
Hence,  $B_{1}(p)\not=0$. 
Since the map $h_p$ is minimal and the free resolution $\V(p-1)$ is minimal, it follows that $B_{0}(p)\not=0$. 
The inequality $\rank(B_1(p))>\rank(B_0(p))$ follows from Proposition~\ref{strictlyincreasing} since
the free resolution $\W(p)$ has a surjective CI operator.
\qed
\myspace
It follows at once that the higher matrix factorization in Example~\ref{unstablemf} is not pre-stable. But in fact Corollary~\ref{zerostop} implies stronger restrictions on the Betti numbers in the finite resolution of modules that are pre-stable syzygies:

\begin{cor}\label{syzbound}
If $M$ is a pre-stable syzygy of complexity $\zeta$ with respect to the regular sequence $\ff c$ in a local ring $ S$ and
$b_{i}^{S}(M)$ denotes the $i$-th Betti number of $M$ as an $S$-module,  then
\begin{align*}
b_{0}^{S}(M)&\geq \zeta\\
b_{1}^{S}(M) &\geq (c-\zeta+1)b_{0}^{S}(M) + \frac{\zeta(\zeta+1)}{2}-1
\end{align*}
\end{cor}

\proof 
Set $\gamma = c-\zeta+1$. 
By Theorem~\ref{fresthm}, Theorem~\ref{inresthm}, and Corollary~\ref{zerostop} we get
$$b_{0}^{S}(M) = b_{0}^{R(\gamma)}(M)=\sum_{p=\gamma}^c\, \rank\,B_0(p)\ge \zeta$$ and
\b{align*}
b_{1}^{S}(M) &= b_{1}^{R(\gamma)}(M) + (c-\zeta)b_{0}^{R(\gamma)}(M)=b_{1}^{R(\gamma)}(M) + (c-\zeta)b_{0}^{S}(M)\cr
b_{1}^{R(\gamma)}(M) &=\sum_{p=\gamma}^c\, \rank\,B_1(p)+\sum_{p=\gamma}^{c-1}\, (c-p) \rank\,B_0(p)\cr
&\ge\bigg(c-\gamma+1-1+\sum_{p=\gamma}^c\, \rank\,B_0(p)\bigg)+\sum_{p=\gamma}^{c-1}\, (c-p) \rank\,B_0(p)\cr
&= \zeta-1+ b_0^S(M)+\sum_{p=\gamma}^{c-1}\, (c-p) \rank\,B_0(p)\, . 
\end{align*}
Therefore, 
\begin{align*}
 b_{1}^{S}(M) =
 &(c-\zeta+1)b_{0}^{S}(M)+\zeta-1+\sum_{p=\gamma}^{c-1}\, (c-p) \rank\,B_0(p)
 \cr
& \ge (c-\zeta+1)b_{0}^{S}(M)+ \zeta-1+{\zeta\choose 2}\cr 
 &=  (c-\zeta+1) b_{0}^{S}(M)+\frac{\zeta(\zeta+1)}{2}-1.
\end{align*}
\vglue -.8cm \ \qed

\myspace\myspace
For example, a pre-stable syzygy module of complexity $\geq 2$ cannot be cyclic and cannot have $b_{1}^{S}(M) = b_{0}^S(M)+1$.

We close this section by a remark on  the graded case. We use the formulas in \ref{epstwo} to study quadratic complete intersections in \cite{EPS2}.

\b{cor}\label{epstwo}
Let $k$ be an infinite field, $S=k[x_1,\dots ,x_n]$ be standard graded with $\deg(x_i)=1$ for each $i$, and $I$ be an ideal generated by a regular sequence of $c$ homogeneous elements of the same degree $q$.
Set $R = S/I$, and suppose that $N$ is a finitely generated  graded $R$-module. Let
$f_1,\dots ,f_c$ be a generic for $N$ regular sequence of forms  minimally generating $I$.
If   $M$
 is a sufficiently high graded syzygy of $N$ over $R$, then $M$ is the  module of a minimal higher  matrix  factorization
 $(d,h)$ with respect to $f_{1}, \dots, f_{c}$;
 it involves modules $B_s(p)$
 for $s=0,1$ and $1\le p\le c$. Denote $b_{i,j}^R(M)=\dim\, \Big(\tor^R_i(M,k)_j\Big)$ the graded Betti numbers of $M$ over $R$.
The graded Poincar\'e series ${\poin}^R_{M}(x,z)=\sum_{i\ge 0}\, b_{i,j}^R(M)x^iz^j$ of $M$ over $R$ is
\eq
 {\poin}^R_{M}(x,z)=\sum_{1\le p\le c}\ {1\over{(1-x^2z^q)^{c-p+1}}} \Big(x\, m_{p;1}(z)+m_{p;0}(z)\Big)\, ,
 \l{series}
 \end{equation} 
where for each $s=0,1$ and $1\le p\le c$ we  use the polynomial $$m_{p;s}(z):=\sum_{j\ge 0}\, b_{s,j}^S\Big(B_s(p)\Big)z^j$$ such that its coefficient
 $b_{s,j}^S\Big(B_s(p)\Big)$ is the number of minimal generators of degree $j$ of
 the $S$-free module $B_s(p)$.
\e{cor}

\proof
By Remark~{\rm \ref{resolvable}}, it follows that Corollary~{\rm \ref{highsyzy}} holds verbatim, without first localizing at the maximal ideal. 
Note that (\ref{series})  is a refined version of the formula in  Corollary~\ref{bettiinf}(1). 
The CI operators $t_i$ on the minimal $R$-free resolution $\bf T$, constructed in \ref{inres}, of the  HMF module $M$
can be taken  homogeneous. Since they are projections by Proposition~\ref{speciallifting1}(2), it follows that they have degree $0$.
The lifted (to $S$) CI operators $\tilde t_i$ satisfy $$\tilde d^{\,2}=f_1\tilde t_1+\cdots +f_c\tilde t_c\, .$$ 
Therefore, $\deg(t_i)=-q$ for every $i$.  \qed

\section{Stable Syzygies  in the Local Gorenstein case}\label{SGor}

 In this section $S$ will denote a local Gorenstein ring. We write $\ff c$ for a regular sequence in $S$ and $R = S/(\ff c)$. Thus $R$ is also a Gorenstein ring.
In this setting
matters are simplified by the fact that a maximal Cohen-Macaulay module is, in a canonical way, an $m$-th syzygy for any $m$. 
 
When $M$ is a maximal Cohen-Macaulay $S$-module
we let $\cosyz^{S}_{j}(M)$ be the dual of the $j$-th syzygy
of $M^{*} := \Hom_{S}(M,S)$.
When we speak of syzygies  or cosyzygies, we will implicitly suppose that they are taken with respect to a minimal resolution.
The following result is well-known.

\theoremstyle{plain}
\newtheorem{clemma}[thm]{Cosyzygy Lemma}

\b{clemma}\label{cosyzygy}
 Let $S$ be a local Gorenstein ring.
\b{mylist}
\item[\rm (1)] If $M$ is a maximal Cohen-Macaulay $S$-module, 
then $M^{*}$ is
 a maximal Cohen-Macaulay $S$-module, $M$ is reflexive, and $\ext_{S}^{i}(M,S) = 0$ for all $i>0$. 
\vglue .1cm
\item[\rm (2)] If $M$ is the first syzygy module
in a minimal free resolution of a maximal Cohen-Macaulay $S$-module, then $M$ has no free summands.
\vglue .1cm
\item[\rm (3)] If $M$ is a maximal Cohen-Macaulay module without free summands, then $$M\cong \syz^{S}_{j}(\cosyz^{S}_{j}(M ))\cong \cosyz^{S}_{j}(\syz^{S}_{j}(M))$$
 for every $j\geq 0$, and $N := \cosyz^{S}_{j}(M)$ is the unique maximal Cohen-Macaulay $S$-module $N$ without free
 summands such that $M$ is  isomorphic to $\syz^{S}_{j}(N)$.
\end{mylist}

\end{clemma}

\noindent {\prooffont Proof Sketch:\ } 
 After replacing $S$ by its completion we may choose a regular local ring $S'\subseteq S$ over
which $S$ is finite, and we have $\ext_{S}(M,S) = \ext_{S'}(M,S')$, and 
$M$ is free over $S'$. Part (2) is obvious over an artinian ring, and the general case follows by factoring out a maximal regular sequence. The first statement of (3) follows from the vanishing
of $\ext$, and the second part follows from the first.
\qed

\myspace
When $M$ is a maximal Cohen-Macaulay module over the Gorenstein ring $S$, we define the Tate resolution of $M$ to be the doubly infinite free complex $\T$ without homology that results from splicing the minimal free resolution of $M$ with the dual of the minimal free resolution of $M^{*}$. If $N$ is also an $S$-module then the {\it stable Ext} is by definition the collection of functors $\extt^{j}(M,N)$, the $j$-th homology of $\hom(\T,N)$; here $j$ can be any integer.

Let $\ff c$ be a regular sequence in a  Gorenstein local ring $S$ with maximal ideal $\mm$ and
residue field $k$. Set $R = S/(f_{1}, \dots, f_{c})$.
Let $M$ be a maximal Cohen-Macaulay $R$-module with no free summands and finite projective dimension over $S$.
If
$\T$ is the Tate resolution of $M$ over $R$, the CI operators corresponding to $\ff c$ are defined on all of $\T$, so that $\extt_{R}(M,k):= \oplus_{i}\extt_{R}^{i}(M,k)$ becomes a graded
module over  the ring $\RR = k[\chi_{1}, \dots, \chi_{c}]$. Then $$\extt^{\geq j}_{R}(M,k) = \ext_{R}(\cosyz_{j}^{R}(M), k)[j]$$ is a finitely generated module
 over $\RR$ for any integer $j$. 
In this case the definition of a stable syzygy (Definition~\ref{wellb}) takes a particularly canonical form:

\b{prop}\label{highsyzdef} 
With hypotheses as above,
$M$ is stable with respect to
$\ff c$ if and only if either $c=0$ and $M=0$, or the following two conditions are satisfied:
\b{mylist}
\item[\rm (1)] $\chi_c$ is 
a non-zerodivisor on $\extt^{\geq -2}_R(M,\allowbreak  k)$. 

\item[\rm (2)]
$
\syz^{R'}_{2}(\cosyz^{R}_{2}(M))
$
is a \highsyz with respect to $f_{1},\dots, f_{c-1}\in S$, where $R'=S/(f_{1},\dots, f_{c-1})$.
\end{mylist}
\end{prop}

\proof $\extt^{\geq -2}_R(M,\allowbreak  k)$ is, up to a shift in grading,
 the same as
$\ext_{R}(\cosyz_{2}^{R}(M), k)$,
and $\cosyz_{2}^{R}(M)$ is the only maximal Cohen-Macaulay module of which
$M$ could be the second syzygy.
\qed

\myspace
We will show that stable syzygies all come from stable matrix factorizations.

\b{thm}\label{stable} 
Let $f_{1},\dots, f_{c}$ be a regular sequence in 
a Gorenstein local ring $S$, and set $R=S/(f_{1}, \dots, f_{c})$. An $R$-module $M$ is a stable syzygy if and only if it is the module of a  minimal stable matrix factorization with respect to $f_{1},\dots, f_{c}$.
\end{thm}

We postpone the proof to give a necessary homological construction:

\b{prop}\label{cosyzres}
Let $f_{1},\dots, f_{c}$ be a regular sequence in 
a Gorenstein local ring $S$, and set $R=S/(f_{1}, \dots, f_{c})$. 
Let $M$ be the HMF module of a minimal stable matrix factorization $(d,h)$. Then
$$\cosyz_2^{R(p)}M(p)=\coker(R(p)\otimes b_{p}) =\coker(R(p-1)\otimes b_{p}) \, .$$
In the notation of Proposition~{\rm\ref{prestable resolutions}},
the minimal $R(p-1)$-free resolution of the module $\cosyz_2^{R(p)}M(p)$ is
$\V(p-1)$, and 
the minimal $R(p)$-free resolution of  the module $\cosyz_2^{R(p)}M(p)$
is 
$\W(p)$.
\end{prop}

\proof
We apply Proposition~{\rm\ref{prestable resolutions}}.
As the higher matrix factorization is stable, we conclude that the depth of the $R(p-1)$-module $\coker(R(p-1)\otimes b_{p})$ is one less than that of a maximal Cohen-Macaulay $R(p-1)$-module.
Therefore, it is a maximal Cohen-Macaulay $R(p)$-module.
The free resolution $\W$ implies that $\cosyz_2^{R(p)}M(p)=\coker(R(p)\otimes b_{p})$.
\ \qed

\b{cor}\label{cosyzcor}
 Let $f_{1},\dots, f_{c}$ be a regular sequence in 
a Gorenstein local ring $S$, and set $R=S/(f_{1}, \dots, f_{c})$. If $M$  is the module of a  minimal stable matrix factorization with respect to $f_{1}, \dots, f_{c}$, then
$$ M(p-1)\cong\syz_{2}^{R(p-1)}\Big(\cosyz_{2}^{R(p)}\Big( M(p)\Big)\Big)\, .$$
\end{cor}

\proof
For
each $p=1,\dots, c$, by Proposition~\ref{cosyzres} we have
\b{align*}
M(p-1)&=\coker(R(p-1)\otimes d_{p-1}) = \syz_{2}^{R(p-1)}\Big(\coker(R(p-1)\otimes b_p)\Big)\\ 
&= \syz_{2}^{R(p-1)}\Big(\cosyz_{2}^{R(p)}\Big( \coker ( R(p)\otimes d_p)\Big)\Big)\\ 
&= \syz_{2}^{R(p-1)}\Big(\cosyz_{2}^{R(p)}\Big( M(p)\Big)\Big)
\end{align*}
where as usual $d_p: \   A_{1}(p)\to A_{0}(p)$ denotes the restriction of $d: \   A_{1}\to A_{0}$.
\qed

\myspace
\noindent{\prooffont Proof of Theorem~\ref{stable}:} 
Theorem~\ref{eventually} shows that
a stable syzygy yields a stable matrix factorization.

Conversely, let $M$ be the module of a  minimal stable matrix factorization $(d,h)$. Use notation as in ~\ref{standardnotation}.
By Proposition~\ref{cosyzres} and in its notation, $\W(p)$ is the minimal $R$-free resolution of 
$\cosyz_2^{R(p)}(M(p))=\coker(R(p)\otimes b_p)\, .$
We have a surjective CI operator $t_c$ on $\W(p)$ because on the one hand, we have it on $\T(p)$ and 
on the other hand $\W(p)_{\le 3}$ is given by the Shamash construction so we have a surjective standard CI operator on $\W(p)_{\le 3}$.
Furthermore, the standard lifting of $\W(p)$ to $R(p-1)$ starts with $\V(p-1)_{\le 1}$, so in the notation of 
Definition~\ref{wellb} we get $\ker(\widetilde{\delta}_1)=M(p-1)$, which is stable by induction hypothesis.
 \qed

\b{cor}\label{takingsyz}
Let $f_{1},\dots, f_{c}$ be a regular sequence in 
a Gorenstein local ring $S$, and set $R=S/(f_{1}, \dots, f_{c})$. Let $M$ be a stable syzygy with a  minimal stable matrix factorization $(d,h)$.
For every $p=1,\dots ,c$ we have $$
\Big(\syz_1^{R(p)}(M(p))\Big)(p-1)=
\syz_1^{R(p-1)}\Big(M(p-1)\Big)\, .
$$
\end{cor}

\proof
By induction, it will suffice to prove this assertion for $M=M(c)$.

The syzygy module $\syz_1^{R}(M)$ is stable by Proposition~\ref{higherSyzygies}.
Recall the proof of Proposition~\ref{higherSyzygies} with $L=\cosyz_{2}^{R}(M)$.
The first and last equalities below  are from Corollary~\ref{cosyzcor}, and then we apply (\ref{syzygies}) to get
\b{align*}
\Big(\syz_1^{R}(M)\Big)(c-1)&=
\syz_{2}^{R(c-1)}\left(\cosyz_{2}^{R}\Big(\syz_{1}^{R}(M)\Big)\right)\\  &=\im(\widetilde{\delta}_3)=\syz_{3}^{R(c-1)}\Big(\cosyz_{2}^{R}(M)\Big)\\ 
&=\syz_{1}^{R(c-1)}\left(\syz_{2}^{R(c-1)}\Big(\cosyz_{2}^{R}(M)\Big)\right)
\\  &=\syz_1^{R(c-1)}\Big(M(c-1)\Big)\,.
\end{align*}
\vglue -.7cm\vbox{ \rightline{\qed}}

\myspace
Recall that if $E$ is a graded $\RR $-module then we define the
{\it $\hbox{\rm S}2$-ification}\ of $E$, written $\hbox{\rm S}2(E)$, by the formula
$$
\hbox{\rm S}2(E) = \oplus_{j\in {\bf Z}} H^{0}(\widetilde E (j))
$$
where $\widetilde E$ denotes the coherent sheaf on projective space associated
to $E$. 

\b{prop}\label{extOfHighSyz}
Suppose that $R = S/I$, where $S$ is a regular local ring and $I$ is generated by a regular
sequence, and let $M$ be maximal Cohen-Macaulay $R$-module.
\b{mylist}
\item[\rm (1)] If $M$ is a stable syzygy then $M$ has no free summand. 
\item[\rm (2)] Set $E:= \extt_{R}^{\ge -2}(M,k)$. If $M$ is a stable syzygy, then $\reg\, E =-1$, 
and $E$ coincides with $\hbox{\rm S}2(E)$
in degrees $\geq -2$.
\end{mylist}
\end{prop}

We could restate the last condition of (2) in terms of local cohomology by saying that $H^{1}_{\RR _{+}}(E)$ is 0 in degree $\geq -2$.

\myspace
\proof (1): This follows at once from part (2) of Lemma~\ref{cosyzygy}.

\hglue -.3cm (2):  We do induction on $c$. If $c=1$ then $E$ is free and generated in degrees $-2$ and $-1$, so the result is obvious, and we may suppose $c>1$.

From Proposition~\ref{highsyzdef} we see that $\chi_{c}$ is a
non-zerodivisor on $E$, so
$$
\reg (E) = \reg(E/\chi_{c}E),
$$
and 
Corollary~\ref{reduction} shows 
$(E/\chi_{c}E)^{\geq 0} = \extt^{\geq 0}_{R'}(M',k)$,
where $M' = \syz^{R'}_{2}\Big(\cosyz^{R}_{2}\allowbreak (M)\Big)$.

Since $M'$ is stable, $\chi_{c-1}$ is a non-zerodivisor on 
$E':=\extt^{\geq -2}_{R'}(M',k)$, and thus also
on $E'^{\geq 0} = (E/\chi_{c}E)^{\geq 0}$, so 
$$
H_{(\chi_{1}, \dots, \chi_{c})}^{0}\Big((E/\chi_{c}E)^{\geq 0}\Big)=0.
$$
 Since the modules 
$E', \, E'^{\geq 0}$ and $E/\chi_{c}E$  differ
by modules of finite length, they have the same $i$-th local
cohomology for $i\geq 1$. By induction, $\reg(E')= -1$, 
so 
$\reg(E/\chi_{c}E) = -1$ as well, proving that $\reg E = -1$.

Finally we show that $E$ agrees with $S\hbox{\rm 2}(E)$ in degrees $\geq -2$.
 Since $\chi_{c}$ is a non-zerodivisor on $E$, we see that
$E$ is a submodule of $F := \hbox{\rm S}2(E)^{\geq -2}$. Because
$\reg\,E= -1$ the natural map $E \to \hbox{\rm S}2(E)$ is surjective in 
degrees $\geq -1$.

Thus we need only prove that $E \to \hbox{\rm S}2(E)$ is surjective in degree
$-2$.
By induction, $E^{\geq 0}/\chi_{c}E = \ext^{\geq 0}(M', k)$
 has depth at least 1. But from the exact sequence
$$
0 \to \chi_{c}F/\chi_{c}E \to E^{\geq 0}/\chi_{c}E \to  E^{\geq 0}/\chi_{c}F \to 0
$$
we see that the module of finite length 
$\chi_{c}F/\chi_{c}E$ is contained in $E^{\geq 0}/\chi_{c}E$,
so $\chi_{c}F/\chi_{c}E=0$. Since $\chi_{c}$ is a non-zerodivisor on
$E$, and thus also on $F$, this implies that $F/E=0$ as well.
\qed

\section{Syzygies over Intermediate Rings}\label{Sfunct}
\noindent
In this section we suppose that $S$ is a Gorenstein ring. Let $I\subset S$ be an ideal generated by a regular
sequence, and set $R= S/I$.
Let $N'$ be a finitely generated $R$-module of finite projective dimension over $S$. 
If $M = \syz_{i}^{R}(N')$ is a  sufficiently high syzygy, then by Theorem~\ref{highsyzprop} and Theorem~\ref{eventually} 
 $M$ comes from a higher matrix factorization with respect to a generic choice of generators
$f_1,\dots ,f_c$ for the ideal $I$. Set $R(p)=S/(\ff p)$.
The following result identifies the HMF module
$M(p)$ with
the module $\syz_{i}^{R(p)}(N)$, where we have chosen  $N = \cosyz_{c+1}^{R}(M)$.

\b{thm}\label{agreement}\label{functoriality}
Let $\ff c$ be a regular sequence
in a local Gorenstein  ring $S$. 
Set $R(p)=S/(\ff p)$ and $R = R(c)$. 
Suppose that $M$  is a stable syzygy with stable matrix factorization $(d,h)$ 
with respect to  $\ff c$.
Let $N = \cosyz_{c+1}^{R}(M)$, and set $M(0)=0$.
\b{mylist}
\item[\rm (1)]
With notation as in  {\rm \ref{standardnotation}},
$$
\syz_{c+1}^{R(p)}(N) \cong M(p)\ \ \hbox{for}\ p\ge 0.
$$
\item[\rm (2)]
We write $\nu_{p}$ for the map
$$
R(p)\otimes \syz_{i}^{R(p-1)}(N)\  \ar{\nu_p} \ \syz_{i}^{R(p)}(N), 
$$
induced by the comparison map from the minimal $R(p-1)$-free resolution of $N$ to
the minimal $R(p)$-free resolution of $N$ inducing the identity map on $N$ (this comparison map is unique up to 
homotopy).
For each $p$, there is a short exact sequence
$$
0\to R(p)\otimes M(p-1)\ar{\nu_{p}} M(p) \ar{}\cosyz_{2}M(p)\to 0.
$$
\end{mylist}
\e{thm}

\myspace
For the proof of Theorem~\ref{functoriality} we will make use of the following well-known result. For the reader's convenience we sketch the proof.  Write $\mod (R)$ for the category of finitely generated 
$R$-modules and
$\MCM(R(p))$ for the stable category of maximal Cohen-Macaulay 
$R(p)$-modules, where the morphisms are morphisms in $\mod(R(p))$ modulo
those that factor through projectives.  We say that $S$-modules
$M,M'$
{\it have a common syzygy} if there exists a $j$ such that $\syz_j^S(M)\cong \syz_j^S(M')$ in $\MCM(S)$.

\b{lemma}\label{differentSyzygies} \label{differentSyzygiestwo}
Suppose that $S$ is a 
Gorenstein ring and that $M, M'$ are  $S$-modules.
\b{mylist}
\item[\rm (1)]
If $N, N'$ are $S$-modules and 
 there are exact sequences
\b{align*}
0\to &M \to P_{r} \to\cdots\to P_{0}\to N\to 0,\\ 
0\to &M' \to P'_{r} \to\cdots\to P'_{0}\to N'\to 0
\end{align*}
such that each $P_{i}$ and each $P'_{i}$ is a module of finite projective dimension
over $S$, then $M$ and $M'$ have a common syzygy if and only if
$N$ and $N'$ have a common syzygy.

\item[\rm (2)] If $M$ and $M'$ have a common syzygy and are both
maximal Cohen-Macaulay $S$-modules then
$M\cong M'$ in $\MCM(S)$.

\item[\rm (3)] If $M\cong M'$ 
in $\MCM(S)$, the ring $S$ is local, and both $M$ and $M'$ are
maximal Cohen-Macaulay $S$-modules without free
summands, then $M\cong M'$ as $S$-modules.
\end{mylist}
\end{lemma}

\proof \noindent{\rm (1):}
It suffices to do the case $r=0$. 
Let $N_{1}= \ker (P_{0}\to N)$, and let $\V$ be a free resolution of $N_{1}$. The mapping
cone of a map from $\V$ to a finite resolution of $P_{0}$ is a free resolution of $N$,  so that for
$i\gg 0$ we have $\syz_{i}^{S}(N) \cong \syz_{i-1}(N_{1})$ in  ${\bf MCM}(S)$. By induction, for $i\gg 0$
the $(i-1-r)$-th syzygy of $M$ agrees with the $i$-th syzygy of $N$,
and the same is true for $M'$ and $N'$. 
\myspace

\hglue -.3cm {\rm (2):} If $ \syz^{S}_{j}(M) \cong \syz^{S}_{j}(M') \cong N$, then
$M \cong \cosyz^{S}_{j}(N) \cong M'$ in $\MCM(S)$.
\vglue .1cm

\hglue -.3cm{\rm (3):} Let
$M\ar{\alpha} M' \ar{\beta} M$ be inverse isomorphisms in $\MCM(S)$. This means that $\beta\alpha = \id_{M}+\phi\varphi$,
where 
$
M\ar{\varphi} F\ar{\phi} M
$
for some  free module $F$. 
Since $S$ is local and $M$ has no free summand,
$\varphi$ must have image inside
the maximal ideal times $F$, and thus $\phi\varphi$ has image inside the
maximal ideal times $M$. By Nakayama's Lemma, $\beta\alpha$ is an epimorphism, and it follows
that $\beta\alpha$ is an isomorphism. Since the same goes for $\alpha\beta$,
we see that $M\cong M'$.
\qed

\myspace
\noindent{\prooffont Proof of Theorem~\ref{agreement}}:\\
\noindent (1):
By Corollary~\ref{mcm} $M(p)$ is a maximal Cohen-Macaulay $R(p)$-module, and by Corollary~\ref{filteredFiniteResolution} it has no free summand. In particular, $N = \cosyz_{c+1}^{R}(M)$ is well-defined
and has no free summands. It follows that $\syz_{c+1}^{R(p)}(N) $
is a maximal Cohen-Macaulay $R(p)$-module and by the Cosyzygy Lemma~\ref{cosyzygy} it has no free summands.
By  Lemma~\ref{differentSyzygies}(3), it suffices to show
that the maximal Cohen-Macaulay $R(p)$-modules
$M(p)$ and $ \syz_{c+1}^{R(p)}(N)$ have a syzygy in common over $R(p)$.
We will do this by showing that each of these modules has an $R(p)$-syzygy in common with $M$.

 Observe that $R$ has finite projective dimension over  $R(p)$.
Lemma~\ref{differentSyzygies}(1) implies that, indeed, $M=\syz_{c+1}^{R}(N)$ and $\syz_{c+1}^{R(p)}(N)$ have a common syzygy over $R(p)$.

We next compare $M = M(c)$ with $M(p)$. When $p>q$ the module $R(p)$ has finite projective dimension over $R(q)$.
By Corollary~\ref{cosyzcor},  $$M(p-1) = \syz_{2}^{R(p-1)}\Big(\cosyz_{2}^{R(p)}(M(p))\Big)\,.$$
Applying Lemma~\ref{differentSyzygies}(1) to an $R(p-1)$-free resolution of $\cosyz_2^{R(p)}(M(p))$ and to an $R(p)$-free resolution of $\cosyz_2^{R(p)}(M(p))$, we conclude that $M(p-1)$ and $M(p)$ have a common syzygy over each ring $R(q)$ with $q\le p-1$.

\myspace
\noindent (2):
For each $p$, let $\T(p)$ be the minimal $R(p)$-free resolution of $M(p)$ and let $\W(p)$ be the minimal $R(p)$-free resolution of $\cosyz_{2}^{R(p)}(M(p))$. See also Proposition~\ref{cosyzres}.
Since $M(p-1)$ is a maximal
Cohen-Macaulay $R(p-1)$-module by Corollary~\ref{mcm},  the minimal free resolution
of $R(p)\otimes M(p-1)$ as an $R(p)$-module is  $R(p)\otimes \T(p-1)$.

Since $M(p)$ is a stable syzygy, the CI operator 
$t_{p}$ is surjective on
$\W(p)$. Take a lifting $\tilde t_{p}$ acting on a lifting of $\W(p)$ to $R(p-1)$. 
The kernel of $\tilde t_{p}$ is a minimal $R(p-1)$-free resolution $\widetilde{\G}$ of
$\cosyz_{2}^{R(p)}(M(p))$. By Corollary~\ref{cosyzcor}, $\T(p-1)$ is isomorphic to
$\widetilde{\G}_{\ge 2}[-2]$.
Thus we have a short exact sequence of minimal free resolutions
$$
0\to R(p)\otimes \T(p-1) \to \T(p) \ar{t_{p}} \W(p)[-2] \to 0,
$$
and this induces the desired short exact sequence of modules.

The last claim in the theorem follows from  Corollary~\ref{takingsyz}.
\qed

\b{cor}\label{preservationOfCodimOne} 
With hypotheses as in Theorem~ {\rm \ref{functoriality}}, let $M$ be a stable syzygy with respect to $\ff c$, with stable matrix factorization
$(d,h)$. If we denote the codimension $1$ part of $(d,h)$ by $(d_1, h_1)$, then the codimension $1$ part of the higher matrix factorization of
$\syz_{1}^{R}(M)$ is $(h_1, d_1)$.
\end{cor}

\proof
If $(d_1,h_1)$  is non-trivial, then the minimal $R(1)$-free resolution of $M(1) = R(1)\otimes d_1$ is periodic of the form
$$
\cdots \ar{d_1}F_{4}\ar{h_1}F_{3}\ar{d_1}F_{2}\ar{h_1}F_{1}\ar{d_1}F_{0}.
$$
\vglue -.7cm
\ \qed

\b{thm}\label{codimOneEquivalence}\label{apl} 
Suppose that
${\ff c} $ is a regular sequence in a Gorenstein local ring $S$, and set
$R=S/(\ff c)$.
Suppose that $N$ is an $R$-module  of finite projective dimension over $S$.
Assume that ${\ff c}$ are generic with respect to $N$. Denote $\gamma:=c-\hbox{cx}_R(N)+1$, where $\hbox{cx}_R(N$ is 
the complexity of $N$ (see Corollary~ {\rm \ref{complexity}}). Then:
\b{mylist}
\item[\rm (1)]
 The projective dimension of $N$ over $R(p)=S/(f_1,\dots ,f_p)$ is finite for $ p< \gamma$.
\item[\rm (2)]  Choose a $j\ge 1$ large enough so that $M:=\syz_j^R(N)$ is a stable syzygy 
and
$\syz_{j}^{R(p)}(N)$ is a maximal Cohen-Macaulay $R(p)$-module for every $p\le\gamma$.
The hypersurface matrix factorization for the periodic part of the minimal free resolution of  $N$ over $S/(f_1,\dots ,\allowbreak f_{\gamma})$ is
isomorphic to the top non-zero part of the higher matrix factorization of $M$.
\end{mylist}
\end{thm}

\noindent A version of (1) is proved in \cite[Theorem 3.9]{Av}, \cite[5.8 and 5.9]{AGP}.

\myspace
\proof
Choose $M$ as in (2).
 By Corollary~\ref{complexity}, $M(p)=0$ for $ p<\gamma$.
Apply Proposition~\ref{agreement} for $ p\le\gamma$.
The case $p=\gamma$ establishes (2).
\qed

\b{remark}
In particular, the above theorem shows that
 the codimension
$1$ matrix factorization that is obtained from a high $S/(f_{1})$-syzygy
of $N$ agrees with the codimension $1$ part of
the higher matrix factorization for $M$ over $R$, and both codimension
$1$ matrix factorizations are trivial if the complexity of $N$ is $<c$, where
$M$ is a sufficiently high syzygy of $N$ over $R$.
\end{remark}

\section{Functoriality}\label{Sfunctorial}

From Theorem~\ref{agreement} it follows immediately that the higher matrix factorization construction induces functors on the stable module category. In this section we  make the result explicit. 

Let $R(p) = S/(\ff p)$. If $i>\dim\, R$ then the modules $\syz_{i}^{R(p)}(N)$ are maximal Cohen-Macaulay $R(p)$-modules.
We define functors 
 $$
 {\mathcal F}_{i}: \mod(R)\to\prod_{p}\Mor\Big(\MCM(R(p))\Big)
 $$
taking $N$ to the collection of morphisms
$$
R(p)\otimes \syz_{i}^{R(p-1)}(N)\  \ar{\nu_p} \ \syz_{i}^{R(p)}(N), 
$$
where $\nu_{p}$ is the comparison map defined in  Theorem~\ref{agreement}. The map $\nu_{p}$ is unique up to 
homotopy, and thus yields a well-defined morphism in $\MCM(R(p))$.

\b{cor}\label{functcor} With assumptions and notation as in Theorem~{\rm\ref{agreement}}, for each $p=1\dots c-1$ there exists  a triangle
in $\MCM(R(p+1))$ of the form

\vglue .3cm
\begin{center}
\begin{tikzpicture}  
 [every node/.style={scale=0.9},  auto]
   \node (A) {$R(p+1)\otimes M(p) $} ;
  \node (B) [node distance=5cm, right of=A] {$M(p+1)$} ;
  \node (BM) [node distance=2.5cm, right of=A] {} ;
  \node (C) [node distance=2.5cm, below of=BM] {$M(p+1)[-2]:=\cosyz_2^{R(p+1)}(M(p+1))\ .  $} ;
  \draw[->] (A) to node  {$\nu_p$} (B);
  \draw[->] (B) to node  {} (C);
  \draw[<-] (A) to node {$[1]$} (C);
\end{tikzpicture}
\end{center}

\noindent If $M'$ is a first syzygy of $M$,
then the corresponding triangles for $M'$ are obtained from the triangles for $M$ by applying the shift (equivalently, taking  first syzygy) operator to each $M(p)$. \qed
\end{cor}

We remark that  Theorem~\ref{functoriality} implies that  for $i\ge c+3$   we get a triangle

\vglue .3cm\begin{center}
\begin{tikzpicture}  
 [every node/.style={scale=0.9},  auto]
   \node (A) {$R(p)\otimes \syz_{i}^{R(p-1)}(N) $} ;
  \node (B) [node distance=5cm, right of=A] {$\syz_{i}^{R(p)}(N)$} ;
  \node (BM) [node distance=2.5cm, right of=A] {} ;
  \node (C) [node distance=2.5cm, below of=BM] {$\syz_{i-2}^{R(p)}(N)\ .$} ;
  \draw[->] (A) to node  {$\nu_{p}$} (B);
  \draw[->] (B) to node  {} (C);
  \draw[<-] (A) to node {$[1]$} (C);
\end{tikzpicture}
\end{center}

Let $\MF(\ff c)$ be the full subcategory of $\MCM(R)$ whose objects are stable equivalence classes of maximal Cohen-Macaulay
modules that  are stable syzygies with respect to $\ff c$.
We get a functor
$
{\mathcal F}: \   \MF(\ff c) \to {\mathcal C},
$
where 
 an  object $\M$ of $\C$ is a collection of objects
$M(p) \in \MCM(R(p))$ for $p=1,\dots,c$  that fit into
triangles as in Corollary~\ref{functcor} in $\MCM(R(p+1))$ 
and whose morphisms
$\M = \{M(p)\} \to \M' = \{M'(p)\}$ are collections of morphisms 
$\Big\{\Big(M(p)\to M'(p)\Big) \in   \MCM(R(p)\Big\}$ that commute with the morphisms in 
the triangles. Furthermore, if $M'$ is the first syzygy of $M$,
then ${\mathcal F}(M')$ is obtained from ${\mathcal F}(M)$ by applying the shift (equivalently, taking  first syzygy) operator
in $\MCM(R(p))$ to each $M(p)$ and to each triangle.

\section{Morphisms}\label{SMorph}

\noindent In this section we introduce the concept of an  HMF morphism (a morphism of higher matrix factorizations) 
so that it preserves the structures described in Definition~\ref{matrixfiltration}, and then we 
show that any homomorphism of
matrix factorization modules induces an HMF morphism.

\b{definition}\label{mfhom}
 A {\it morphism of matrix factorizations\/} or {\it HMF morphism} $\alpha: (d,h)\to (d',h')$
 is a triple
 of homomorphisms of free modules
\b{align*}
\alpha_{0}: A_{0}&\to A'_{0}\\ 
\alpha_{1}:A_{1}&\to A'_{1}\\ 
\alpha_{2}:\oplus_{p\leq c}A_{0}(p)&\to \oplus_{p\leq c}A'_{0}(p)
\end{align*}
such that, for each $p$:
\b{mylist}
\item[\rm (a)] $\alpha_{s}(A_{s}(p))\subseteq A'_{s}(p)$ for $s=0,1$. We write
$\alpha_{s}(p)$ for the restriction of $\alpha_{s}$ to $A_{s}(p)$.
\item[\rm(b)] $\alpha_{2}\left(\oplus_{q\leq p}A_0(q)\right) \subseteq  \oplus_{q\leq p}A_0'(q)$,
and the component $ A_{0}(p) \to A'_{0}(p)$ of $\alpha_2$ is $\alpha_{0}(p)$.
We write $\alpha_{2}(p)$ for the restriction of $\alpha_{2}$ to $\oplus_{q\leq p}A_0(q)$.
\item[\rm (c)] The diagram 

\begin{center}
\begin{tikzpicture}  
 [every node/.style={scale=1},  auto]
\node(A1){$\oplus_{q\le p}\, A_0(q)$};
\node(A2)[node distance=3cm, right of=A1]{$A_1(p)$};
\node(A3)[node distance=2.5cm, right of=A2]{$A_0(p)$};
 \draw[->] (A1) to node {$h$} (A2);
  \draw[->] (A2) to node {$d_p$} (A3);
\node(B1)[node distance=2cm, below of=A1] {$\oplus_{q\le p}\, A_0'(q)$};
\node(B2)[node distance=2cm, below of=A2]{$A_1'(p)$};
\node(B3)[node distance=2cm, below of=A3]{$A_0'(p)$};
 \draw[->] (B1) to node {$h'$} (B2);
  \draw[->] (B2) to node {$d'_p$} (B3);
  \draw[->] (A1) to node [swap] {$\alpha_2(p)$} (B1);
   \draw[->] (A2) to node   {$\alpha_1(p)$} (B2);
  \draw[->] (A3) to node  {$\alpha_0(p)$} (B3);
\end{tikzpicture}
\end{center}

\noindent commutes  modulo $(\ff{p-1})$.
\end{mylist}
\end{definition}

\b{thm}\label{mfmorph} 
Suppose that
${\ff c} $ is a regular sequence in a Gorenstein local ring $S$, and set
$R=S/(\ff c)$.
 Let  $M$ and $M'$ be
stable syzygies over $R$, and suppose $\zeta: M\to M'$ is a morphisms of $R$-modules.
With notation as in  {\rm \ref{standardnotation}}, let $M$ and $M'$ be HMF modules of
stable matrix factorizations $(d,h)$ and $(d',h')$, respectively.
There exists an HMF morphism 
$$\alpha: (d,h)\to (d',h')$$ such that the map induced  on $$M= \coker (R\otimes d)\to \coker(R\otimes d') = M'$$ is $\zeta$. 
\end{thm}

We first establish a strong functoriality statement for the Shamash construction.
Suppose that $\G$ and $\G'$ are  $S$-free resolutions of $S$-modules $M$ and $M'$ annihilated by a
non-zerodivisor $f$, and $\zeta: M\to M'$ is any homomorphism. If we choose systems of higher
homotopies $\sigma$ and $\sigma'$ for $f$ on $\G$ and $\G'$ respectively, then the Shamash construction yields resolutions $\std(\G,\sigma)$ and $\std(\G',\sigma')$ of $M$ and $M$ over
$R = S/(f)$, and thus there is a morphism of complexes 
$$\tilde \phi: \std(\G,\sigma)\to\std(\G',\sigma')$$
covering $\zeta$. To prove the Theorem we need more: a morphism  defined over $S$ that commutes with the maps in the  ``standard liftings'' $\tilde \std(\G,\sigma)$ and 
$\tilde \std(\G', \sigma')$ (see Construction~\ref{standardCI}) and  respects
the natural filtrations of these modules. The following statement provides the required
morphism.

\b{lemma}\label{higherhomotopieswithmap} 
Let $S$ be a commutative ring, 
and let $\varphi_{0}:(\G,d) \to (\G',d')$ be a map of $S$-free resolutions of modules annihilated by 
an element $f$. Given systems of
higher homotopies $\sigma_j$ and $\sigma'_j$ on $\G$ and $\G'$, respectively, there exists a system of maps $\varphi_{j}$ of degree $2j$
from the underlying free module of $\G$ to that of $\G'$ such that, for every index $m$,
$$
\sum_{i+j=m}(\sigma'_{i}\varphi_{j}-\varphi_{j}\sigma_{i}) = 0.
$$
\end{lemma}

We say that $\{\varphi_{j}\}$ is
 a {\it system of homotopy comparison maps} if they satisfy the conditions in the lemma above.

Recall that a  map of free complexes $ \lambda:\, \U\to\W[-a]$ is a homotopy for
a  map $\rho:\,  \U\to\W[-a+1]$ if $\delta \lambda-(-1)^a \lambda\partial=\rho$, where
$\partial$ and $\delta$ are the differentials in $\U$ and $\W$ respectively.
Since in Lemma~\ref{higherhomotopieswithmap}  $\sigma_{0}$ and $\sigma'_{0} $ are the differentials $d$ and $d'$, the equation above in Lemma~\ref{higherhomotopieswithmap} says that, for each $m$, the map $\varphi_{m}$ is a homotopy for the sum
$$
-\sum_{i+j=m\atop i>0,j>0}(\sigma'_{i}\varphi_{j}-\varphi_{j}\sigma_{i}).
$$

\myspace
\proof
The desired condition on $\varphi_{0}$ is equivalent to the given hypothesis that $\varphi_{0}$ is a map
of complexes. We proceed by induction on $m>0$  and on homological degree to prove the existence of
$\varphi_{m}$. The desired condition can be written as
$$
d'\varphi_m=-\sum_{i+j=m\atop  i\neq 0}\, \sigma'_i
\varphi_j+\sum_{i+j=m}\,\varphi_j\sigma_i\,.
$$
Since $\G$ is a free resolution, 
it suffices to show that the right-hand side is annihilated by $d'$. Indeed,
\b{align*}
&-\sum_{i+j=m\atop i\neq 0}\, (d'\sigma'_i)\varphi_j+
\sum_{i+j=m}\, (d'\varphi_j)\sigma_i\\  
&=
\sum_{i+j=m\atop i\neq 0}\,\sum_{v+w=i\atop v\neq 0}\,\sigma'_v\sigma'_w\varphi_j
-f\varphi_{m-1} -\sum_{i+j=m}\,\sum_{q+u=j\atop q\neq 0}\,\sigma'_q\varphi_u\sigma_i+
\sum_{i+j=m}\,\sum_{u+q=j}\,\varphi_u\sigma_q\sigma_i\\ 
&= \sum_{v+w+j=m\atop v\neq 0}\,\sigma'_v\sigma'_w\varphi_j
-f\varphi_{m-1} -\sum_{i+q+u=m\atop q\neq 0}\,\sigma'_q\varphi_u\sigma_i
+\sum_{i+u+q=m}\,\varphi_u\sigma_q\sigma_i
\, ,
\end{align*}
where the first  equality holds by  (3) in \ref{hhomotop}  and by the induction hypothesis.
Reindexing the first summand by $v=q$, $w=i$ and $j=u$ we get 
\b{align*}
& \sum_{q+i+u=m\atop q\neq 0}\,\sigma'_q\sigma'_i\varphi_u
-f\varphi_{m-1}  -\sum_{i+q+u=m\atop q\neq 0}\,\sigma'_q\varphi_u\sigma_i
+\sum_{i+u+q=m}\,\varphi_u\sigma_q\sigma_i\\ 
&=-f\varphi_{m-1}+ \sum_{q\neq 0}\, \sigma'_u\Bigg(\sum_{i+u=m-q}\,\sigma'_i\varphi_u-\varphi_u\sigma_i\Bigg)+\sum_{u}\varphi_u\Bigg(\sum_{q+i=m-u}\,\sigma_q\sigma_i\Bigg)\\ 
&= -f\varphi_{m-1}  +0+0+\varphi_{m-1}f=                       0\,,
\end{align*}
where the last  equality holds by  (3) in \ref{hhomotop}  and by induction hypothesis.
\qed

\myspace
The next result reinterprets the conditions of Lemma~\ref{higherhomotopieswithmap} as defining a map
between liftings of Shamash resolutions. 

\b{prop}\label{standardmap}
Let $S$ be a commutative ring, 
and let $\G$ and $\G'$ be  $S$-free resolutions with systems of
higher homotopies $\sigma = \{\sigma_j\}$ and $\sigma' = \{\sigma'_j\}$ for $f\in S$, respectively. Suppose that
$\{\varphi_{j}\}$ is a system of homotopy comparison maps for $\sigma$ and $\sigma'$.
We use the {\it standard lifting\/} of the Shamash resolution
defined in  {\rm \ref{standardCI}}, and the notation established there.
Denote by 
$\widetilde{\varphi}$ the map
with components 
$$
\varphi_{i}: y^{(v)}G_{j} \to y^{(v-i)}G_{j+2i}'
$$
from the underlying graded free $S$-module of the
standard lifting $\tilde \std(\G,\sigma)$ of the Shamash  resolution $\std(\G,\sigma)$,
 to the underlying graded free $S$-module of the
standard lifting $\tilde \std(\G',\sigma')$ of the Shamash  resolution $\std(\G',\sigma')$. 
The maps $\widetilde{\varphi}$ satisfy $\widetilde{\delta'}\widetilde{\varphi}=\widetilde{\varphi}\widetilde{\delta}\,,$ where $\widetilde{\delta}$ and $ \widetilde{\delta'}$ are the standard liftings of the differentials
defined in  {\rm \ref{standardCI}}.
\end{prop}

\proof
Fix $a$ and $v$. We must show that the diagram
\begin{center}
\begin{tikzpicture}  
 [every node/.style={scale=0.9},  auto]
\node(A1){$y^{(a)}G_{v}$};
\node(A2)[node distance=6cm, right of=A1]{$\oplus_{0\le i\le a}\, y^{(a-i)}G_{v+2i-1}$};
 \draw[->] (A1) to node {$ \widetilde{\delta}$} (A2);
\node(B1)[node distance=2.5cm, below of=A1] {$\oplus_{0\le j\le a}\, y^{(a-j)}G_{v+2j}'$};
\node(B2)[node distance=2.5cm, below of=A2]{$\oplus_{0\le j\le a\atop 0\le i\le a-j}\, y^{(a-i-j)}G_{v+2i+2j-1}'\, .$};
 \draw[->] (B1) to node [swap] {$ \widetilde{\delta'}$} (B2);
  \draw[->] (A1) to node {$\widetilde{\varphi}$} (B1);
   \draw[->] (A2) to node {$\widetilde{\varphi}$} (B2);
\end{tikzpicture}
\end{center}
\vglue .1cm
\noindent
commutes. Fix   $0\le q\le a$. The map $ \widetilde{\delta'}\widetilde{\varphi}-\widetilde{\varphi} \widetilde{\delta}$
from $y^{(a)}G_{v}$ to $y^{(q)}G_{v+2a-2q-1}'$ 
is equal to
$ \sum_{i+j=a-q}\,(\sigma'_{i}\varphi_{j}-\varphi_{j}\sigma_{i} )$, which vanishes by Lemma~\ref{higherhomotopieswithmap}.
\qed

\b{remark} 
A simple modification of the proof of Lemma~\ref{higherhomotopieswithmap} shows that systems of homotopy comparison maps also exist in the context of systems of higher homotopies for a regular sequence $\ff c$, not just in the case $c=1$ as above, and one can interpret this in terms of Shamash
resolutions as in Proposition~\ref{standardmap} as well, but we do not need these refinements.
\end{remark}

\noindent{\prooffont Proof of Theorem~\ref{mfmorph}:} The result is immediate for $c=1$, so
we proceed by induction on $c>1$. Let $\widetilde R = S/(f_{1},\dots, f_{c-1})$. To simplify the
notation, we will write
$\widetilde -$ for $\widetilde R\otimes_{S} -$, and $\overline -$ for $R\otimes -$. We will make use of
our standard notation~\ref{standardnotation}.

Since $(d,h)$ is stable we can extend the map $\overline d$ to a complex
$$
\overline{A}_{1}(c)\to \overline{A}_{0}(c)\to \overline{B}_{1}(c)\to \overline{B}_{0}(c)
$$
that is the beginning of an $R$-free resolution $\F$ of $\cosyz_{2}^{R}(M)$, and there is
a similar complex that is the beginning of the $R$-free resolution $\F'$ of
$\cosyz_{2}^{R}(M')$. By stability these cosyzygy modules are maximal Cohen-Macaulay
modules, so dualizing these complexes we may use $\zeta(c):=\zeta: M\to M'$ to induce maps
\b{align*}
 \eta:&\ \cosyz_{2}^{R}(M)\to \cosyz_{2}^{R}(M')\\ 
 \lambda:&\ \F\to \F'\, .\, 
\end{align*}
\vglue -.2cm
\noindent
Moving to $\widetilde R$, we have 
\b{align*}
M(c-1) &= \coker\, \widetilde {d}(c-1) = \syz_{2}^{\widetilde R}(\cosyz_{2}^{R}(M))\\ 
M'(c-1) &= \coker\, \widetilde {d}'(c-1) = \syz_{2}^{\widetilde R}(\cosyz_{2}^{R}(M'))\, .\\ 
\end{align*}

We will use the notation and the construction in the proof of Theorem~\ref{eventually}, where we produced an $\widetilde R$-free resolution
$\V$ of $\cosyz_{2}^{R}(M)$, and various homotopies on it. Of course we have a similar resolution $\V'$ of $\cosyz_{2}^{R}(M)$. See the diagram:
\vglue .2cm
\begin{center}
\begin{tikzpicture}  
 [every node/.style={scale=0.8},  auto]
\node(000) {$\V:\ \cdots  $};
\node(00)[node distance=2.3cm, right of=000]{$
{\displaystyle\oplus\atop{q\le c-1}}\widetilde{A}_0(q)$};
\node(01)[node distance=3cm, right of=00]{$\widetilde{A}_1(c-1)$};
\node(02)[node distance=3cm, right of=01]{$\widetilde{A}_0(c-1)$};
\node(03)[node distance=3cm, right of=02]{$ \tilde B_{1}(c)$};
\node(04)[node distance=3cm, right of=03]{$\tilde B_{0}(c)$};
\draw[->] (000) to node  { } (00);
\draw[->] (01) to node [swap]{$\widetilde{d}_{c-1}$} (02);
\draw[->] (00) to node [swap]{$\widetilde{h}$} (01);
\draw[->] (02) to node [swap]{$\tilde \partial_2$} (03);
\draw[->] (03) to node [swap]{$\tilde b$} (04);
\draw[->,  bend right=45] (03) to node[swap] {$\tilde \psi$} (02);
\draw[->,  bend right=45] (04) to node  [above=8pt, right=-20pt]  {$\tilde \theta_0$} (03);
\draw[->,  bend right=45] (02) to node  [above=8pt, right=10pt]  {$\tilde\theta_2$} (01);
\draw[->,  bend right=50] (04) to node [above=1pt, right=50pt] {$\tilde \tau_0$} (01);
\node(A00) [node distance=3.5cm, below of=000] {$\V':\ \cdots  $};
\node(A0)[node distance=2.3cm, right of=A00]{$
{\displaystyle\oplus\atop{q\le c-1}}\widetilde{A}_0'(q)$};
\node(A1)[node distance=3cm, right of=A0]{$\widetilde{A}_1'(c-1)$};
\node(A2)[node distance=3cm, right of=A1]{$\widetilde{A}_0'(c-1)$};
\node(A3)[node distance=3cm, right of=A2]{$ \tilde B_{1}'(c)$};
\node(A4)[node distance=3cm, right of=A3]{$\tilde B_{0}'(c)\,.$};
\draw[->] (A00) to node  { } (A0);
\draw[->] (A1) to node {$\widetilde{d}_{c-1}'$} (A2);
\draw[->] (A2) to node {$\tilde \partial_2'$} (A3);
\draw[->] (A3) to node {$\tilde b'$} (A4);
\draw[->] (A0) to node {$\widetilde{h}'$} (A1);
\draw[->,  bend left=45] (A3) to node{$\tilde \psi'$} (A2);
\draw[->,  bend left=45] (A4) to node [above=-8pt, right=-20pt] {$\tilde \theta_0'$} (A3);
\draw[->,  bend left=45] (A2) to node[above=-8pt, right=10pt]{$\tilde\theta_2'$} (A1);
\draw[->, bend left=50] (A4) to node   [above=-1pt, right=50pt]  {$\tilde \tau_0'$} (A1);
\draw[->, red] (03) to node[above=14pt, right=-15pt] {$\tilde\xi_{1}$} (A3);
\draw[->, red] (01.270) to node[above=14pt, right=-49pt] {$\tilde\alpha_1(c-1)$} (A1.90);
\draw[->, red] (00.270) to node[above=14pt, right=-49pt] {$\tilde\alpha_2(c-1)$} (A0.90);
\draw[->, red] (02.270) to node[above=14pt, right=-49pt] {$\tilde\alpha_0(c-1)$} (A2.90);
\draw[->, red] (04) to node [above=14pt, right=-15pt] {$\tilde\xi_{0}$} (A4);
\draw[->, red] (04) to node [above=-8pt, right=-45pt] {$\tilde\varphi_1$} (A2);
\draw[->, red] (03) to node [above=-8pt, right=-45pt] {$\tilde\varphi_1$} (A1);
\draw[->, red] (02) to node [above=-8pt, right=-45pt] {$\tilde\varphi_1$} (A0);
\draw[->, red]    (04)  edge[out=229,in=270]   node[above=-5pt, right=-85pt]  {$\tilde\varphi_2$} (A0); 
\end{tikzpicture}
\end{center}
\vglue -9.5cm
\eq \label{diagramone} \  \end{equation}
\vglue 6.5cm
\noindent The map $\eta$ induces  $\tilde\xi: \V\to \V'\, ,$
which in turn induces a map
$$\zeta(c-1): M(c-1)\to M'(c-1)\, .$$ 
See diagram  (\ref{diagramone}).

By induction,
the map $\zeta(c-1)$ is induced by an HMF morphism
 with components 
 $$\alpha_{s}(c-1): A_{s}(c-1)\to A'_{s}(c-1)$$ for $s=0,1$
and  
$$\alpha_{2}(c-1): {\oplus\atop{q\le c-1}}\widetilde{A}_0(q)\to {\oplus\atop{q\le c-1}}\widetilde{A}_0(q)'\, .$$

By the conditions in \ref{mfhom}, it follows that the first two squares on the left are commutative; clearly, the last square on the right is commutative as well. Since $\tilde\alpha_{0}(c-1)$ induces the same map on $M(c-1)$ as $\tilde\xi$, 
the remaining square commutes.
Therefore we can apply Lemma~\ref{higherhomotopieswithmap} and conclude that there 
exists a system of homotopy comparison maps, the first few of which are
shown as $\tilde\varphi_{1}$ and $\tilde\varphi_{2}$ in 
diagram (\ref{diagramone}).

With notation as in (\ref{definedif}) and (\ref{defineho})  we may write the first two steps of the minimal $\tilde  R$-free resolution $\U$  of $M$ 
in the form given by the top three rows of  diagram 
(\ref{diagramtwo}), 
where we have used
a splitting $$A_{0}(c) = A_{0}(c-1)\oplus B_{0}(c)$$ to split the left-hand term
$\oplus_{q\leq c}A_{0}(q)$ into three parts, $\oplus_{q\leq c-1}A_{0}(q), \ A_{0}(c-1),$ and $B_{0}(c)$; and similarly for $M'$ and the bottom three rows.
Straightforward computations using the definition of the system of homotopy comparison maps
shows that  diagram 
(\ref{diagramtwo}) 
commutes.

\eq \label{diagramtwo} \  \end{equation} \vglue -1cm \ 
\begin{center}
\hglue -.75cm \vbox{
\begin{tikzpicture}  
[every node/.style={scale=0.85},  auto]
\node(00) {$\widetilde{A}_0(c-1)$};
\node(11)[node distance=1.5cm, above of=00]{$
{\displaystyle\oplus\atop{q\le c-1}}\widetilde{A}_0(q)$}; 
\node(01)[node distance=4cm, right of=00]{$\widetilde{A}_1(c-1)$};
\node(02)[node distance=4cm, right of=01]{$\widetilde{A}_0(c-1)$};
\node(03)[node distance=1.8cm, below of=01]{$ \tilde B_{1}(c)$};
\node(04)[node distance=1.8cm, below of=02]{$\tilde B_{0}(c)$};
\node(05)[node distance=1.8cm, below of=00]{$\widetilde{B}_0(c)$};
\draw[->] (00) to node  [above=14pt, right=-20pt]{$\tilde \partial_2$} (03);
\draw[->] (03) to node [above=10pt, right=0pt] {$\tilde b$} (04);
\draw[->] (03) to node [above=6pt, right=-12pt] {$\tilde \psi$} (02);
\draw[->] (01) to node{$\widetilde{d}_{c-1}$} (02);
\draw[->] (05) to node [above=-10pt, right=0pt]   {$\tilde \theta_0$} (03);
\draw[->] (00) to node  {$\tilde\theta_2$} (01);
\draw[->] (05) to node  [above=-3pt, right=-35pt]{$\tilde \tau_0$} (01);
\node(A0) [node distance=4.5cm, below of=00]{$\widetilde{A}_0'(c-1)$};
\node(A1)[node distance=4cm, right of=A0]{$\widetilde{A}_1'(c-1)$};
\node(A2)[node distance=4cm, right of=A1]{$\widetilde{A}_0'(c-1)$};
\node(A3)[node distance=1.8cm, below of=A1]{$ \tilde B_{1}'(c)$};
\node(A4)[node distance=1.8cm, below of=A2]{$\tilde B_{0}'(c)\,.$};
\node(A5)[node distance=1.8cm, below of=A0]{$\widetilde{B}_0'(c)$};
\node(A11)[node distance=1.5cm, below of=A5]{$
{\displaystyle\oplus\atop{q\le c-1}}\widetilde{A}_0'(q)$}; 
\draw[->] (A0) to node  [above=14pt, right=-20pt]{$\tilde \partial_2'$} (A3);
\draw[->] (A3) to node [above=-10pt, right=0pt] {$\tilde b'$} (A4);
\draw[->] (A3) to node [above=-6pt, right=2pt] {$\tilde \psi'$} (A2);
\draw[->] (A1) to node [above=10pt,right=-3pt] {$\widetilde{d}_{c-1}'$} (A2);
\draw[->] (A5) to node [above=-10pt, right=15pt]  {$\tilde \theta_0'$} (A3);
\draw[->] (A0) to node  {$\tilde\theta_2'$} (A1);
\draw[->] (A5) to node  [above=-3pt, right=-35pt]{$\tilde \tau_0$} (A1);
\draw[->, red] (05) to node [above=-1pt, right=-1pt] {$\tilde\varphi_1$} (A0);
\draw[->, red] (03) to node [above=-1pt, right=-20pt] {$\tilde\varphi_1$} (A1);
\draw[->, red] (04) to node [above=-1pt, right=-20pt] {$\tilde\varphi_1$} (A2);
\draw[->, red, bend right=45] (00) to node  [above=25pt, right=-45pt]{$\tilde\alpha_{0}(c-1)$} (A0);
\draw[->, red, bend right=45] (05) to node [above=-35pt, right=-9pt]{$\tilde\xi_0$} (A5);
\draw[->, red, bend left=45] (02) to node [above=30pt, right=-5pt]{$\tilde\alpha_{0}(c-1)$} (A2);
\draw[->, red, bend left=45] (03) to node [above=-28pt, right=-6pt]{$\tilde\xi_{1}$} (A3);
\draw[->, red, bend left=45] (04) to node [above=-31 pt, right=-5pt]{$\tilde\xi_{0}$} (A4);
\draw[->, red, bend left=45] (01) to node [above=-8pt, right=-1pt]{$\alpha_{1}(c-1)$} (A1);
\draw[->, red ]    (11)  edge[out=180,in=180]  node[right =-49] 
[above=120pt, right=40pt] {$\tilde\alpha_{2}(c-1)$} (A11);
\draw[->, red,  bend right=60] (05) to node [above=10pt, right=-18pt]{$\tilde\varphi_2$} (A11);
\draw[->, red] (00) edge[out=180,in=160 ]  node[above=90pt, right=15pt] {$\tilde\varphi_{1}$} (A11);
\draw[->] (11) to node{$\tilde h$} (01);
\draw[->] (A11) to node [above=-28pt, right=-28pt]{$\tilde h'$} (A1);
\end{tikzpicture}}
\end{center} 

\myspace
Next, we will construct the maps $\alpha_i$.
We construct $\alpha_0$ by extending the map
$\alpha_{0}(c-1)$ already defined
over $S$ by taking $\alpha_{0}|_{B_{0}(c)}$ to have as components  arbitrary liftings to $S$
of $\tilde\xi_{0}$ and $\tilde\varphi_{1}$. Similarly we take
$\alpha_{1}$ to be the extension of $\alpha_{1}(c-1)$ that has arbitrary liftings 
of $\tilde\xi_{1}$ and $\tilde\varphi_{1}$ as components. Finally, we take
$\alpha_{2}$ to agree with $\alpha_{2}(c-1)$ on $\oplus_{q\leq c-1}A_0({q})$ and
 on  the summand
$A_{0}(c) = A_{0}(c-1)\oplus B_0({c})$, to be the map given by 
$$
\alpha_{0}(c-1): A_{0}(c-1) \to A'_{0}(c-1)
$$ 
and arbitrary liftings 
\b{align*}
\varphi_{1}&: A_{0}(c-1) \to \oplus_{q\leq c-1}A'_{0}(q) & 
\varphi_{1}&: B_{0}(c) \to A'_{0}(c-1)\\ 
\xi_{0}&: B_{0}(c) \to B'_{0}(c) & 
\varphi_{2}&: B_{0}(c) \to  \oplus_{q\leq c-1}\,A'_{0}(q)
\end{align*}
 to $S$ of $\tilde\varphi_{1}, \tilde\varphi_{2}$ and $\tilde\xi_{0}$.

It remains to show that $\overline \alpha_{0} = R\otimes_{S}\alpha_{0}$ induces $\zeta: M\to M'$.

By  Proposition~\ref{cosyzres} the minimal $R$-free resolutions ${\bf F}$ and ${\bf F'}$ of $ \cosyz_{2}^{R}(M)$ and $\cosyz_{2}^{R}(M')$ have the form given in the following diagram.  

\myspace
\begin{center}
\begin{tikzpicture}  
 [every node/.style={scale=0.9},  auto]
\node(000) {$\F:\quad\cdots  $};
\node(02)[node distance=3cm, right of=000]{$ \overline{B}_{0}(c)\oplus\overline{A}_0(c-1)$};
\node(03)[node distance=3.5cm, right of=02]{$  \overline{B}_{1}(c)$};
\node(04)[node distance=3cm, right of=03]{$ \overline{B}_{0}(c)$};
\draw[->] (000) to node  { } (02);
\draw[->] (02) to node {$\overline{(\theta}_0, \overline{\partial}_2)$} (03);
\draw[->] (03) to node {$ \overline{b}$} (04);
\node(A00) [node distance=2.5cm, below of=000] {$\F':\quad\cdots  $};
\node(A2)[node distance=3cm, right of=A00]{$ \overline{B}_{0}'(c)\oplus\overline {A}_0'(c-1)$};
\node(A3)[node distance=3.5cm, right of=A2]{$ \overline{ B}_{1}'(c)$};
\node(A4)[node distance=3cm, right of=A3]{$ \overline{B}_{0}'(c)\, ,$};
\draw[->] (A00) to node  { } (A2);
\draw[->] (A2) to node [swap] {$(\overline{\theta}_0', \overline{\partial}_2')$} (A3);
\draw[->] (A3) to node [swap]{$ \overline{b}'$} (A4);
\draw[->] (03.290) to node[above=5pt, right=2pt] {$ \lambda_{1} = \overline \xi_{1}$} (A3.70);
\draw[->] (02.270) to node[above=5pt, right=-19pt] {$\overline{\alpha}_0$} (A2.90);
\draw[->] (02.290) to node[above=5pt, right=2pt] {$ \lambda_{2}$} (A2.70);
\draw[->] (04.290) to node [above=5pt, right=2pt] {$ \lambda_{0}= \overline \xi_{0}$} (A4.70);
\end{tikzpicture}
\end{center}
\vglue .2cm
\noindent  
By definition the map of complexes $ \lambda:\,\F\to \F'$ induces
$\zeta:\, M\to M'$. Using Lemma~\ref{higherhomotopieswithmap}, we see that the left-hand square
of the diagram also commutes if we replace $ \lambda_{2}$ with $\alpha_{0}$, and thus these two
maps induce the same map $M\to M'$, concluding the proof.
\qed

\section{Strong Higher Matrix Factorizations}\label{SHMFstrong}

We introduce a stronger version of Definition~\ref{matrixfiltration} in which we require that the map $h$ is part of a homotopy. In Theorem~\ref{strongmf} we show that an HMF module always has
a strong matrix factorization.

\b{definition}\l{strong}
Let $(d,h)$ be a  higher matrix factorization and $M=\coker(d\otimes R)$ be its module. We say that 
$(d,h)$ is a {\it strong matrix factorization} for $M$ if
for each $p$, the map $h_p$
can be extended to a homotopy $h(f_p) $ for $f_p$ at  $L(p)_0=A_0(p)$ on the $S$-free resolution $\L(p)$ of $M(p)=\coker\Big({d}_p\otimes S/(f_1,\dots ,f_p)\Big)$
constructed in \ref{finresconstruction}.
\e{definition}

 For example, in the codimension $2$ case in Example~\ref{exampletwo},
 the third equation   in
(\ref{conditions}) shows that 
 a  higher matrix factorization satisfies
 $$\rho dh_2\equiv f_2\rho\ \mod (f_1B_{0}(1))\, ,$$
where  we denote $\rho$ the projection from $A_0$ to
$B_0(1)$. A strong matrix factorization must satisfy  the stronger condition $$\rho dh_2=f_2\rho\, .$$
 
\b{thm} \l{strongmf}\ 
\b{mylist}
\item[{\rm (1)}] If $(d,h)$ is a strong
matrix factorization, then it is a  higher matrix factorization.
\item[{\rm (2)}] Let $(d,h')$ be a  higher matrix factorization.
There exists a  strong
matrix factorization
$(d,h)$ with the same filtrations $
0\subseteq A_{s}(1)\subseteq\cdots \subseteq A_{s}(c) = A_{s}$ for $s=0,1$
as   $(d,h')$. Note that $(d,h)$ and $(d,h')$ have the same higher matrix factorization module $M:=\coker(d\otimes R)$.
If the ring $S$ is local and the  higher matrix factorization $(d,h')$ is  minimal, 
then  $(d,h)$ is minimal as well. 
\end{mylist}
\e{thm}

\proof
We will use the following notation:
if $\varphi$ is a map of modules $$\oplus_{1\le i\le s} P_i\to \oplus_{1\le j\le s} Q_j\,,$$ then we denote
$\varphi_{P_i\to Q_j}$  the projection of $\varphi\big\vert_{P_i}$ on $Q_j$ and call it 
the {\it component of $\varphi$ from $P_i$ to $Q_j$}.

Consider   the finite $S$-free resolution  $\bf L$ of $M$
constructed in \ref{fres}, and use the notation in \ref{rmvanishing} and in \ref{strong}.

(1):  We have to show that conditions (a) and (b) in Definition~\ref{matrixfiltration} are satisfied
by $d$ and $h$.
First, we consider (a). For a fixed $p$,  the map $h(f_p)$ has components 
\begin{align*}
h_p&: \ A_1(p)\longleftarrow A_0(p)\\ h(f_p)_{e_iB_0(w) \leftarrow B_0(v)}&:
\ e_iB_0(w)  \leftarrow  B_0(v)\ \hbox{for} \ i<w \, .\end{align*}
For every $q\le p$ we have
\eq
 d_ph_p\big\vert_{B_0(q)}+\sum_{{1\le i< w\le p}} f_ih(f_p)_{e_iB_0(w) \leftarrow B_0(q)}=f_p\id_{B_0(q)}\, .
\l{strongeq}
\end{equation}
This condition is stronger than condition (a) in  Definition~\ref{matrixfiltration}.

We will prove that
(b) holds. 
Fix $p$, and denote $\partial$ the differential in $\L(p)$. Let $\sigma(p)$ be a homotopy for
$f_p$ on $\L(p)$ that extends $h(f_p)$.
The second differential in $\L(p)$ is mapping
\b{align*}
L(p)_2&=  \Big( \oplus_{i<q\le p}\,e_iB_1(q)\Big)\oplus\Big( \oplus_{j<i<q\le p}e_ie_jB_0(q)\Big) \\
 \downarrow \ &\ \\
 L(p)_1&=\Big(\oplus_{q\le p}B_1(q)\Big)\oplus\Big(\oplus_{i<q\le p}e_iB_0(q)\Big)\, .
 \end{align*}
By Remark~\ref{rmvanishing} the only components of the differential that land in $B_1(p)$ are
  $$f_i: \, e_iB_1(p)\to B_1(p)\ \ \hbox{ for }i<p\,.$$
  Therefore,
\b{align*}
f_p\id_{B_1(p)}=
\pi_p\Big(f_p\id_{A_1(p)}\Big)
&=\pi_p\Big(\sigma(f_p)\partial+\partial \sigma(f_p)\Big)\big\vert_{A_1(p)}\\
&=\pi_p\sigma(f_p)\partial\big\vert_{A_1(p)} +\pi_p\partial \sigma(f_p)\big\vert_{A_1(p)}\\
 & =
  \pi_ph_pd_p+\sum_{1\le i<p} f_i \sigma(f_p)_{e_iB_1(p)\leftarrow A_1(p)}  \\
  &\equiv  \pi_ph_pd_p\ \mod(f_1,\dots ,f_{p-1})B_1(p)\, .
\end{align*}
\vglue -.5cm\ 

(2): For each $p$, there exists a homotopy $\sigma(f_p)$  for $f_p$  on the free resolution ${\bf L}(p)$ since the module $M(p)$ is
annihilated by $f_p$.
Let $h: A_0(p) \to A_1(p)$ be the components of $\sigma(f_p)$ from $A_0(p)$ to $A_1(p)$.

Suppose that $S$ is local and the  higher matrix factorization $(d,h')$  is minimal.
We will prove that the map $h$ is minimal.
By Theorem~\ref{inresthm}, $(d,h')$ yields a minimal $R$-free resolution $\T'$ of $M$.
By (1), $(d,h)$ is a  higher matrix factorization, and so by
Theorem~\ref{inresthm} it  yields an $R$-free resolution $\T$ of $M$ over $R$.
Both resolutions have the same ranks of the corresponding free modules in them because
the free modules in the filtrations of $(d,h')$ and $(d,h)$ have the same ranks.
Therefore, the resolution $\T$ is minimal as well. The second differential in $\T$ is
$\overline{h}=h\otimes R$. Hence, the map $h$ is minimal.
\qed

\vglue .7cm
\noindent{\bf Acknowledgements.} 
We are  grateful to  
Jesse Burke and Hailong Dao who read a draft of the paper and made  helpful comments.
We want to express our thanks to  David Morrison for lecturing to us on the applications of  matrix factorizations in physics. The second author thanks Luchezar Avramov for an introduction to this interesting subject; we also
thank him  and 
Ragnar Buchweitz for useful conversations.
We profited from examples computed in the system Macaulay2 [M2], and we wish to express our thanks to its programmers Mike Stillman and Dan Grayson. 

\vglue .2cm 
Eisenbud is partially supported by NSF grant DMS-1001867, and Peeva is partially supported by NSF grants DMS-1100046, DMS-1406062 and by a Simons Fellowship; both authors are partially supported by NSF 
grant 0932078000 while in residence at MSRI.

\vglue .5cm
\makeatletter
\font \sectionfont = cmbx12 
\renewcommand\section{\@startsection{section}%
{1}{0pt}{-1.8\baselineskip}{.8\baselineskip}%
{\sectionfont \center}}
\makeatother

\bibliographystyle{alpha}

\end{document}